\documentclass[1p]{elsarticle}
\usepackage{soul}
\usepackage{amsmath}
\usepackage{graphicx}
\usepackage{subfig}
\usepackage{array}
\usepackage{amsfonts}
\usepackage[dvipsnames]{xcolor}
\usepackage{booktabs}
\usepackage{physics}
\usepackage{cancel}
\usepackage{pifont}
\usepackage{amssymb}
\usepackage{natbib}
\newsavebox{\measurebox}

\journal{Progress report}

\begin{document}

\begin{frontmatter}

\title{An Unstructured Body-of-Revolution Electromagnetic Particle-in-Cell Algorithm with Radial Perfectly Matched Layers and Dual Polarizations}

\author[ad1]{Dong-Yeop Na}
\ead{dyna22@postech.ac.kr}
\author[ad2]{Fernando L Teixeira}
\ead{teixeira.5@osu.edu}
\author[ad3]{Yuri A Omelchenko}
\ead{omelche@gmail.com}

\address[ad1]{Department of Electrical Engineering, Pohang University of Science and Technology, Pohang, 37673, South Korea}
\address[ad2]{ElectroScience Laboratory and Department of Electrical and Computer Engineering, The Ohio State University, Columbus, OH 43212, USA}
\address[ad3]{Trinum Research Inc., San Diego, CA 92126, USA}

\begin{abstract}
A novel electromagnetic particle-in-cell algorithm has been developed for fully kinetic plasma simulations on unstructured (irregular) meshes in complex body-of-revolution geometries.
The algorithm, implemented in the BORPIC++ code, utilizes a set of field scalings and a coordinate mapping, reducing the Maxwell field problem in a cylindrical system to a Cartesian finite element Maxwell solver in the meridian plane.
The latter obviates the cylindrical coordinate singularity in the symmetry axis.
The choice of an unstructured finite element discretization enhances the geometrical flexibility of the BORPIC++ solver compared to the more traditional finite difference solvers.
Symmetries in Maxwell's equations are explored to decompose the problem into two dual polarization states with isomorphic representations that enable code reuse.
The particle-in-cell scatter and gather steps preserve charge-conservation at the discrete level.
Our previous algorithm (BORPIC+) discretized the $\mathbf{E}$ and $\mathbf{B}$ field components of TE\textsuperscript{$\phi$} and TM\textsuperscript{$\phi$} polarizations on the finite element (primal) mesh~\cite{Na2017Axisymmetric,Na2019Finite}.
Here, we employ a new field-update scheme.
Using the same finite element (primal) mesh, this scheme advances two sets of field components independently: (1) $\mathbf{E}$ and $\mathbf{B}$ of TE\textsuperscript{$\phi$} polarized fields, ($E_{z},E_{\rho},B_{\phi}$) and (2) $\mathbf{D}$ and $\mathbf{H}$ of TM\textsuperscript{$\phi$} polarized fields, ($D_{\phi},H_{z},H_{\rho}$).
Since these field updates are not explicitly coupled, the new field solver obviates the coordinate singularity, which otherwise arises at the cylindrical symmetric axis, $\rho=0$ when defining the discrete Hodge matrices (generalized finite element mass matrices).
A cylindrical perfectly matched layer is implemented as a boundary condition in the radial direction to simulate open space problems, with periodic boundary conditions in the axial direction.
We investigate effects of charged particles moving next to the cylindrical perfectly matched layer.
We model azimuthal currents arising from rotational motion of charged rings, which produce TM\textsuperscript{$\phi$} polarized fields.
Several numerical examples are provided to illustrate the first application of the algorithm.
\end{abstract}


\end{frontmatter}


\section{Introduction}

Plasma turbulence is a plasma state where intermittent self-generated electromagnetic and electrostatic fields emerge as a result of nonlinear interactions of charged plasma particles with waves that grow from parent linear modes due to plasma instabilities. Turbulence arises in all large-scale plasma environments, ranging from laboratory and fusion experiments to astrophysical, space, and ionospheric plasmas. Understanding the nature of this ubiquitous phenomenon is important in many respects because turbulence affects the transport of plasma particles, momentum, and energy in these environments. The Sun-Earth system is particularly dominated by the supersonic and super-Alfv\'{e}nic “solar wind” plasma constantly flowing radially from the Sun. This turbulent plasma carries a “frozen-in” interplanetary magnetic field that “tangles” with the Earth’s dipolar field, inducing geomagnetic “storms” and “substorms” in the Earth’s environment~\cite{liu2011dynamic}.

Of special importance to geophysics are extreme space weather events that can be hazardous to space assets. For instance, enhancement of flux in the radiation belts formed by energetic electrons trapped in the Earth's dipolar magnetic field may result in disabling satellites. To scatter energetic electrons from their magnetic field aligned trajectories, various “radiation belt remediation” schemes have been proposed, among which the most promising ones rely on (i) injecting a heavy ion beam across the ambient magnetic field~\cite{ganguli2019understanding} and (ii) deploying high-power antennas or electron beams~\cite{carlsten2019radiation} to directly drive very-low-frequency (VLF) waves. In the former case, VLF whistler waves can be generated through nonlinear (electron induced scattering) conversion of ion beam excited electrostatic lower hybrid oscillations~\cite{ganguli2019understanding}; in the latter schemes, the whistler modes (longer wavelengths) and X-type modes (shorter wavelengths) are proposed to be driven with VLF antennas or field-aligned modulated electron beams~\cite{carlsten2019radiation}.

Nonlinear interaction of plasma wave modes with plasma ions and electrons can be studied self-consistently with particle-in-cell (PIC) or Vlasov simulations in three dimensions (3D).  For instance, the generation of lower hybrid and whistler waves by an ion velocity ring beam was studied with the full-PIC VPIC~\cite{winske2012generation} and hybrid X-HYPERS~\cite{omelchenko2021rate}  codes, while the excitation of whistler waves by loop and dipole antenna was modeled via full-PIC LSP simulations~\cite{main2018excitation}, and pitch-angle scattering by a field-aligned pulsed electron beam was investigated with the Vlasov SPS code~\cite{delzanno2019high}. Since kinetic simulations are generally computationally expensive, they typically need to be scaled with respect to realistic plasma parameters and space resolution in order to extract essential 3D physics~\cite{omelchenko2021rate}.  Physically, however, 3D simulations may often be assumed to be axisymmetric in the Cartesian ($x,y$) space, with the $z$-axis being defined by the direction of external magnetic field. Essentially, the assumption of axisymmetry reduces a 3D simulation model to a computationally two-dimensional (2D) problem. This geometric simplification enables higher resolution and longer simulations of nonlinear plasma interactions.  For instance, the inherent noise driven by finite-size particles in PIC simulations may prevent 3D PIC models from resolving subtle nonlinear effects, such as induced electron scattering, which is fundamental to the nonlinear theory of ion beam evolution~\cite{winske2012generation}. Moreover, computational electromagnetic 3D models may break the exact axisymmetry of the original problem because of discretization effects.

In addition to being reduced in terms of dimensionality, physical parameters, and mesh resolution, full-PIC simulations (e.g.,~\cite{winske2012generation}) may apply periodic boundary conditions which do not account for convection of electromagnetic waves out of the system and therefore limit physical simulation time. Motivated by the need to simulate the aforementioned physical plasmas in more detail, below we describe a novel axisymmetric full PIC unstructured mesh simulation model. Some of the key features of the model are: (1) The use of a radial perfectly matched layer (PML) absorbing boundary condition designed to model open domain problems, (2) the use of an unstructured-mesh finite-element discretization for fields in the meridian ($z \rho$) plane to provide better geometric flexibility and adaptation to complex geometries, (3) the use of judicious operator and field rescalings to map the problem in the $z \rho$ domain to an equivalent problem in the Cartesian $xy$ plane and enable the reuse of simpler Cartesian finite-element codes, and (4) the decomposition of the problem into two polarizations, which removes the coordinate singularity at the axis ($\rho=0$)  and explores a duality between the polarizations to enable the reuse of the same computer code, with only minor adaptations, for both polarizations.

The main motivation behind the development of the present  BORPIC++ algorithm is to more efficiently and accurately simulate nonlinear interactions of charged plasma particles and waves which produce plasma instabilities in the presence of the cylindrical symmetry, particularly, the generation of whistler waves by the ion velocity ring distribution.
The BORPIC++ algorithm can be useful in performing numerical experiments for various radiation belt remediation schemes.

Our previous BORPIC algorithm \cite{Na2017Axisymmetric,Na2019Finite} has two limitations for this purpose: ($i$) it does not account for open radial boundary conditions and ($ii$) it does not account for azimuthal currents resulting from the rotational motion of charged particles.
To resolve such issues, in the present BORPIC++ algorithm incorporates the following features:
\begin{itemize}
\item Azimuthal currents arising from the rotational motion of uniform charged rings are modeled, producing TM\textsuperscript{$\phi$} polarized fields.
\item A novel field discretization scheme is used to obviate the cylindrical coordinate singularity along the symmetry axis (i.e., $\rho=0$). The discretization is based on the discrete exterior calculus of differential forms and uses a field scaling strategy to represent $\mathbf{E}$ and $\mathbf{B}$ (ordinary forms) for TE\textsuperscript{$\phi$} polarized fields and $\mathbf{D}$ and $\mathbf{B}$ (twisted forms) for TM\textsuperscript{$\phi$} polarized fields in a finite element (primal) mesh such that discrete Hodge matrices (generalized finite element mass matrices) for both polarizations are well-conditioned and singularity-free. We note that dual meshes are defined mathematically in the formulation but are not required in the actual numerical implementation.
\item Radial PMLs are implemented to avoid artificial interactions between plasma particles and reflected waves in the radial direction. We investigate the effects of charged particles moving nearby the radial PML.
\end{itemize}

This paper is organized as follows. Section 2 discusses the key physical assumptions built into the model. Section 3 describes the field solver and the finite-element discretization in detail, including the domain mapping and the treatment of the dual field polarizations.  Sections 4, 5, and 6 describe the gather, pusher, and scatter algorithms for the computational superparticles (each representing a charge ring in 3D space) used to model the plasma medium. Section 7 provides a series of numerical results to verify the accuracy and illustrate the basic capabilities of the model. Finally, Section 8 summarizes the main takeaways from this work.

\section{Body-of-Revolution Electromagnetic Particle-in-Cell Algorithm}
\begin{figure}
\centering
\includegraphics[width=1.0\linewidth]{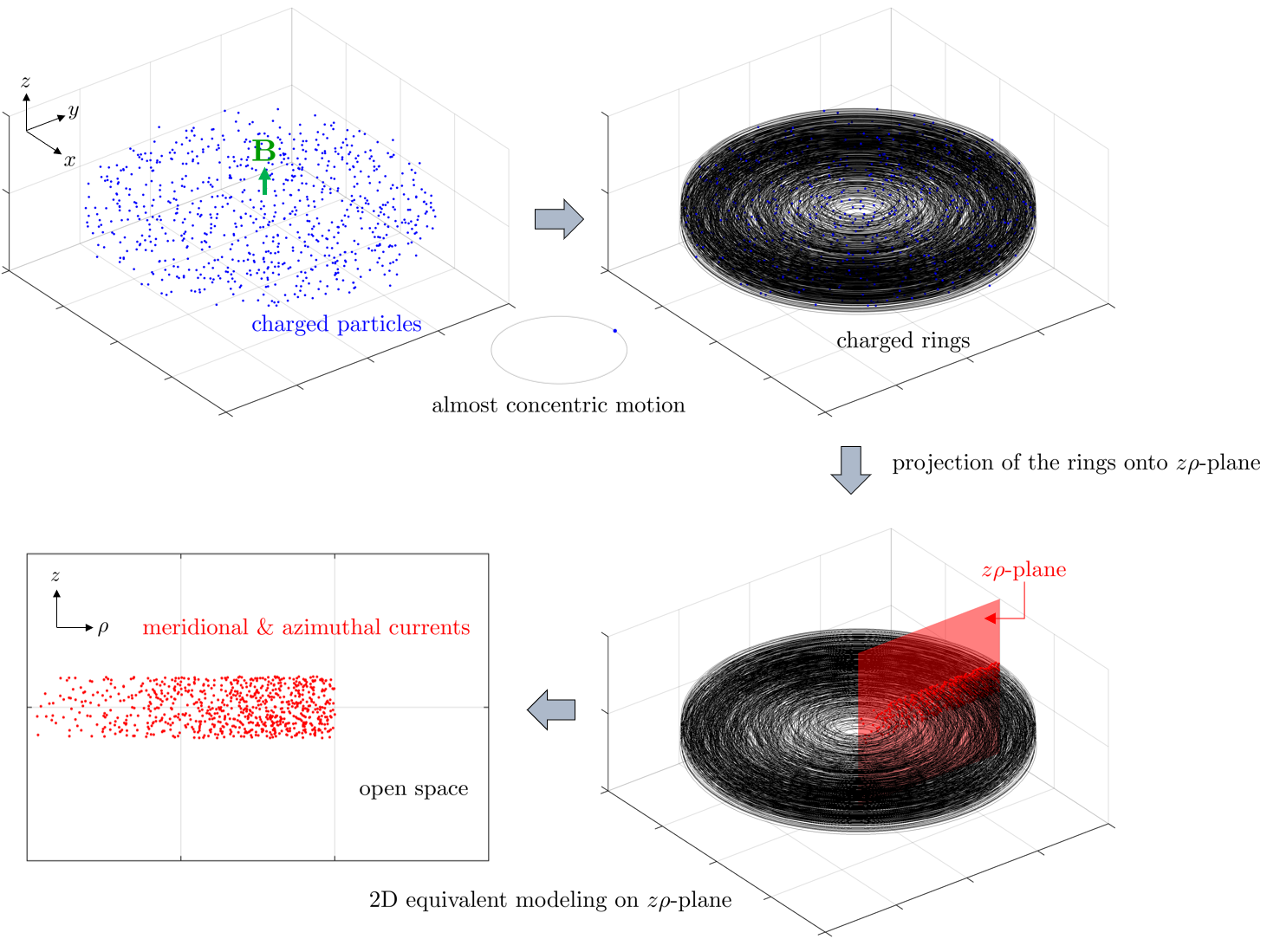}
\caption{Approximation of the motion of charged particles in the presence of a strong axial magnetic field.}
\label{fig:BOR_schematic}
\end{figure}

Body-of-revolution problems consider electromagnetic fields which do not change in the azimuthal direction; hence, for the $\phi$ direction, Fourier eigenmodes with the eigen index $m=0$ is assumed throughout the article.
Fig. \ref{fig:BOR_schematic} illustrates how the actual motion of charged particles in the presence of a strong axial magnetic field is approximated by 2D equivalent modeling on the $z\rho$-plane.

The present body-of-revolution electromagnetic particle-in-cell (BORPIC++) algorithm solves Maxwell-Vlasov equations describing the time evolution of collisionless plasmas on a spatial domain discretized by a finite element mesh. The algorithm computes and updates the interaction between the dynamic electromagnetic field and charged particles during successive time steps. The electromagnetic field may include internal (i.e., the Maxwell field dependent of the particle distribution) and external (independent of the particle distribution) components.  As usual in PIC simulations, coarse-graining of the phase space is assumed whereby each computational (super)particle represents a large number of actual particles (electrons and/or ions) in the plasma medium~\cite{VINCENTI201822}. This decreases the required amount of computational resources. In the BORPIC++ setting, the actual particles associated with a given computational superparticle are uniformly distributed along a circular ring in 3D space. The position of the superparticle corresponds to the intersection of that circular ring with the meridian $z \rho$-plane (more on this below). The interaction between the particles is assumed to be fully mediated through the collective electromagnetic field, with Coulomb-like interactions being negligible due to Debye shielding. The latter assumption is valid as long as the mesh element size is smaller than the Debye length of the plasma medium under consideration.

The BORPIC++ can model physical systems involving a body-of-revolution geometry wherein uniform charged rings can move along the $z$ direction while expanding, shrinking, and rotating.
Such motions of the charged rings can be translated into axial, radial, and azimuthal electric current densities.
Since the charged rings are with a uniform charge density, axial, radial, and azimuthal current densities are independent of azimuthal angle $\phi$; consequently, all terms related to the derivative with respect to $\phi$ in original Maxwell's equations vanish.
This key assumption enables us to decompose total electromagnetic fields into two polarizations: (i) TE\textsuperscript{$\phi$} ($E_{z},E_{\rho},B_{\phi}$) and (ii) TM\textsuperscript{$\phi$} ($E_{\phi},B_{z},B_{\rho}$).
The former is produced by axial and radial motions of the charged rings, and the latter results from the rotational motions.
Since the two sets of polarized fields do not interact with each other, one can handle them separately in the 2D space of $z\rho$-plane.

We should note that a related algorithm (BORPIC) was presented before in~\cite{Na2017Axisymmetric}. However, the BORPIC algorithm modeled the TE\textsuperscript{$\phi$} polarized field (three components) only, as opposed to both TE\textsuperscript{$\phi$} and TM\textsuperscript{$\phi$} fields here (all six electromagnetic field components).
Moreover, the BORPIC algorithm did not incorporate absorbing boundary conditions, being restricted to the analysis of bounded spatial regions.
We have extended our previous BORPIC algorithm to the BORPIC+ algorithm \cite{Na2019Finite} by accounting for both TE\textsuperscript{$\phi$}  and TM\textsuperscript{$\phi$} polarized fields with higher order Fourier azimuthal eigenmodes.
However, in the BORPIC+ algorithm, rotational motions of uniform charged rings are not considered so that the BORPIC+ algorithm is not able to model general motions of the charged rings.
Moreover, the BORPIC+ algorithm does not incorporate open boundary conditions in the radial directions.
More importantly, in the BORPIC+ algorithm, the azimuthal\footnote{Here, $A_{\phi}$ is the azimuthal component of a vector field, $\mathbf{A}$ and $A_{\rho}$ and $A_{z}$ are the transverse (poloidal) components of this vector.} component of $\mathbf{E}$ and transverse components of $\mathbf{B}$ of TM\textsuperscript{$\phi$} polarized fields (instead of the azimuthal component of $\mathbf{D}$ and transverse components of $\mathbf{H}$) were discretized in a primal mesh; consequently, their numerical approximations were not stable close to the $z$-axis because the discrete Hodge matrix encoding the (mapped) permittivity constitutive relation $D_{\phi} = \frac{\epsilon_0}{\rho}E_{\phi}$ exhibits a singularity $\rho=0$.

The new BORPIC++ algorithm resolves the aforementioned issues by: (1) considering axial, radial, and rotational motions of uniform charged rings and taking into account both (TE\textsuperscript{$\phi$} and TM\textsuperscript{$\phi$}) field polarizations, (2) employing absorbing boundary conditions in the radial direction to avoid artificial interactions between plasma particles and reflected waves, and (3) constructing a novel discretization scheme that removes the $1/\rho$ singularity previously present in the finite element mass matrices in the cylindrical system.

The axial singularity, $1/\rho$ in BORPIC+ resulted when building discrete Hodge matrices (generalized finite element mass matrices) using the constitutive relation, $D_{\phi} = \frac{\epsilon_0}{\rho}E_{\phi}$.
This relation was necessary because the BORPIC+ algorithm discretized the $\mathbf{E}$ and $\mathbf{B}$ field components on the finite element (primal) mesh~\cite{Na2017Axisymmetric,Na2019Finite} for both (TE\textsuperscript{$\phi$} and TM\textsuperscript{$\phi$}) polarizations.
To remove this axial singularity in the discrete Hodge matrices, in BORPIC++ we discretize the $\mathbf{D}$ and $\mathbf{H}$ field components of TM\textsuperscript{$\phi$} polarization (instead of the $\mathbf{E}$ and $\mathbf{B}$ components thereof) on the same primal mesh.
To summarize, the new scheme advances two sets of field components independently on the same finite element (primal) mesh: (1) $\mathbf{E}$ and $\mathbf{B}$ of TE\textsuperscript{$\phi$} polarized fields, ($E_{z},E_{\rho},B_{\phi}$) and (2)  $\mathbf{D}$ and $\mathbf{H}$ of TM\textsuperscript{$\phi$} polarized fields, ($D_{\phi},H_{z},H_{\rho}$).
Since these two update schemes are not explicitly coupled to each other, the new field solver obviates the coordinate singularity at the cylindrical symmetry axis, $\rho = 0$ when defining the discrete Hodge matrices (generalized finite element mass matrices). More details on this are presented in Section 3.

During each time step interval, the BORPIC++ algorithm cycles through four stages sequentially: Field Solver,  Gather, Pusher, and Scatter~\cite{vincenti2017pic,10.1063/5.0046842,Na2018Relativistic}. Each of these four stages is described separately in the next Sections.

\section{Finite-Element Field Solver Stage}
\subsection{Cylindrical-to-Cartesian domain mapping}
Consider an object or medium with azimuthal symmetry along the $z$-axis, e.g., a cylindrical or annular waveguide structure or plasma. We denote it as a BOR (body-of-revolution), for short.
The field solver employs a finite element algorithm based on an unstructured mesh with triangular elements.
This provides maximum geometric flexibility in adapting to complex BOR geometries and plasmas.
 We explore the BOR symmetry of the problem by assuming no variation along the $\phi$ direction, which reduces the discretization of the problem to the meridional $z \rho$-plane~\cite{10.1063/1.1384387,mahalingam2010particle,LEHE201666,MASSIMO2016841}.

Although the use of cylindrical coordinates to analyze BOR problems brings computational advantages, vector differential operators such as the curl operator exhibit radial scaling factors not found in the Cartesian coordinate system and in addition a coordinate singularity at $\rho=0$.
These properties can often make certain aspects of the discretization in the  $z \rho$-plane more burdensome. However, by utilizing the transformation-optics principles~\cite{Teixeira1999Lattice,Teixeira1999Differential,Pendry2006Controlling}, it is possible to map Maxwell's equations from the cylindrical system to a Cartesian system where the metric factors are fully embedded into (artificial) constitutive tensors~\cite{Na2019Finite}.
This strategy allows the easy adaptation of finite element codes based on a Cartesian domain to the cylindrical domain.

We start by writing  Maxwell’s equations in the form
\begin{flalign}
\label{eqn:me1}
\nabla\times \mathbf{E} &= -\frac{\partial}{\partial t}\mathbf{B}, \\
\label{eqn:me2}
\nabla\times \mathbf{H} &= \frac{\partial}{\partial t}\mathbf{D}+\mathbf{J},\\
\label{eqn:me3}
\mathbf{D} &=\bar{\bar{\epsilon}}\cdot \mathbf{E}, \\
\label{eqn:me4}
\mathbf{B} &=\bar{\bar{\mu}}\cdot \mathbf{H}.
\end{flalign}
where $\mathbf{E}$ and $\mathbf{B}$ are the electric and magnetic fields in a medium with (generally tensor) electric permittivity $\bar{\bar{\epsilon}}$ and (generally tensor) magnetic permeability $\bar{\bar{\mu}}$, $\mathbf{D}$ and $\mathbf{H}$ are the displacement and magnetizing fields, and $\mathbf{J}$ is the volumetric current density. Assume the above curl operators and vector fields expressed in cylindrical coordinates and components, and that the constitutive parameters of the physical medium,  also expressed in the cylindrical basis, are given by
\begin{flalign}
\bar{\bar{\epsilon}}
=
\left[
\begin{matrix}
\epsilon_{\rho} & 0 & 0
\\
0 & \epsilon_{\phi} & 0
\\
0 & 0 & \epsilon_{z}
\end{matrix}
\right],
\quad
\bar{\bar{\mu}}
=
\left[
\begin{matrix}
\mu_{\rho} & 0 & 0
\\
0 & \mu_{\phi} & 0
\\
0 & 0 & \mu_{z}
\end{matrix}
\right].
\end{flalign}
It should be noted that we explicitly deal with the displacement field $\mathbf{D}$ and magnetizing field $\mathbf{H}$ even in the vacuum background where $\mathbf{D}=\epsilon_0 \mathbf{E}$ and $\mathbf{H}=\mu_0^{-1} \mathbf{B}$ due to numerical purposes as well as the symmetry of the Maxwell equations.
A comprehensive discussion on the numerical motivations can be found in Section \ref{Sec:Motivation_EB_DH}.

Now, let us consider the new (artificial) anisotropic permittivity and permeability tensors $\bar{\bar{\epsilon}}'$ and $\bar{\bar{\mu}}'$ of the form
\begin{flalign}
\bar{\bar{\epsilon}}' = \bar{\bar{\epsilon}}\cdot \bar{\bar{R}}
=
\bar{\bar{\epsilon}}\cdot
\left[
\begin{matrix}
\rho & 0 & 0
\\
0 & \rho^{-1} & 0
\\
0 & 0 & \rho
\end{matrix}
\right],
\label{eqn:eff_permittivity}
\end{flalign}
\begin{flalign}
\bar{\bar{\mu}}' = \bar{\bar{\mu}}\cdot \bar{\bar{R}}
=
\bar{\bar{\mu}}\cdot
\left[
\begin{matrix}
\rho & 0 & 0
\\
0 & \rho^{-1} & 0
\\
0 & 0 & \rho
\end{matrix}
\right].
\label{eqn:eff_permeability}
\end{flalign}
Similarly, it should be mentioned that in plasmas, the original permittivity and permeability tensors ($\bar{\bar{\epsilon}}$ and $\bar{\bar{\mu}}$) are typically assumed to be unit matrices. On the other hand, the primed tensors mentioned above are artificial (geometric) quantities.
We emphasize here that the fields $\mathbf{D}$ and $\mathbf{H}$ (similarly $\mathbf{E}$ and $\mathbf{B}$) are necessary for expressing the different topological roles of the $d/dt$ and curl operators.
However, it is important to clarify that in the plasma problems we consider, the fields $\mathbf{D}$ and $\mathbf{H}$ do not carry separate physical meanings (versus $\mathbf{E}$ and $\mathbf{B}$) unlike their roles in solid materials, for instance.

Using the following rescaling for the fields:
\begin{flalign}
\mathbf{E}' &= \bar{\bar{U}} \cdot \mathbf{E},
\label{eqn:transfE}
\\
\mathbf{B}' &= \bar{\bar{V}}  \cdot \mathbf{B}
\label{eqn:transfB}
\\
\mathbf{H}' &=  \bar{\bar{U}} \cdot \mathbf{H},
\label{eqn:transfH}
\\
\mathbf{D}' &= \bar{\bar{V}} \cdot \mathbf{D}
\label{eqn:transfD}
\\
\mathbf{J}' &= \bar{\bar{V}} \cdot \mathbf{J}
\label{eqn:transfJ}
\end{flalign}
where
\begin{flalign}
\bar{\bar{U}}
=
\left[
\begin{matrix}
1 & 0 & 0
\\
0 & \rho & 0
\\
0 & 0 & 1
\end{matrix}
\right], \quad
\bar{\bar{V}}
=
\rho \,  \bar{\bar{U}}^{-1}
=
\left[
\begin{matrix}
\rho & 0 & 0
\\
0 & 1 & 0
\\
0 & 0 & \rho
\end{matrix}
\right].
\end{flalign}
For example, when writing the above vector fields more explicitly,
\begin{flalign}
\mathbf{A}(\rho,\phi,z) &= \hat{\rho} A_{\rho}(\rho,\phi,z) + \hat{\phi}A_{\phi}(\rho,\phi,z) + \hat{z}A_{z}(\rho,\phi,z), \\
\mathbf{A}'(\rho,\phi,z) &= \hat{\rho} \underbrace{A_{\rho}(\rho,\phi,z)}_{A'_{\rho}(\rho,\phi,z)} + \hat{\phi}\underbrace{\rho A_{\phi}(\rho,\phi,z)}_{A'_{\phi}(\rho,\phi,z)} + \hat{z} \underbrace{A_{z}(\rho,\phi,z)}_{E'_{z}(\rho,\phi,z)},\\
\mathbf{C}(\rho,\phi,z) &= \hat{\rho} C_{\rho}(\rho,\phi,z) + \hat{\phi} C_{\phi}(\rho,\phi,z) + \hat{z} C_{z}(\rho,\phi,z), \\
\mathbf{C}'(\rho,\phi,z) &= \hat{\rho} \underbrace{\rho C_{\rho}(\rho,\phi,z)}_{C'_{\rho}(\rho,\phi,z)} + \hat{\phi}\underbrace{C_{\phi}(\rho,\phi,z)}_{C'_{\phi}(\rho,\phi,z)} + \hat{z} \underbrace{\rho C_{z}(\rho,\phi,z)}_{C'_{z}(\rho,\phi,z)},
\end{flalign}
for $\mathbf{A}=\mathbf{E}$ or $\mathbf{H}$, and $\mathbf{C}=\mathbf{B}$, $\mathbf{D}$, or $\mathbf{J}$.
We can rewrite the original Maxwell's equations using the artificial medium defined above together with the rescaled fields as
\begin{flalign}
\nabla'\times \mathbf{E}' &= -\frac{\partial}{\partial t}\mathbf{B}', \label{eqn:rescale_FL}
\\
\nabla'\times \mathbf{H}' &= \frac{\partial}{\partial t}\mathbf{D}' +\mathbf{J}', \label{eqn:rescale_AL}
\\
\mathbf{D}' &=\bar{\bar{\epsilon}}'\cdot \mathbf{E}', \label{eqn:rescale_GL}\\
\mathbf{B}' &=\bar{\bar{\mu}}'\cdot \mathbf{H}',\label{eqn:rescale_GLM}
\end{flalign}
where
\begin{flalign}
\nabla'\times\mathbf{A}'(\rho,\phi,z)
=
\left|
\begin{matrix}
\hat{\rho} & \hat{\phi} & \hat{z} \\
\frac{\partial}{\partial \rho} &
\frac{\partial}{\partial \phi} &
\frac{\partial}{\partial z} \\
A_{\rho}' & A_{\phi}' & A_{z}'
\end{matrix}
\right|.
\end{flalign}
In other words, the above modified curl operator in the transformed (primed) system is devoid of any radial scaling factor and thus it is locally
isomorphic to the Cartesian curl operator.

We show the equivalence between the original Maxwell's equations in the cylindrical coordinate system represented by \eqref{eqn:me1} to \eqref{eqn:me4} and the transformed one represented by \eqref{eqn:rescale_FL} to \eqref{eqn:rescale_GLM} in \ref{equivalence}.
Furthermore, it can be easily shown that $\bar{\bar{R}}\cdot\bar{\bar{U}}=\bar{\bar{V}}$ and $\bar{\bar{R}}^{-1}\cdot\bar{\bar{V}}=\bar{\bar{U}}$ such that the constitutive relations for the original Maxwellian field variables in \eqref{eqn:me3} and \eqref{eqn:me4} and those for the transformed field variables in \eqref{eqn:rescale_GL} and \eqref{eqn:rescale_GLM} are identical.

\subsection{Field decomposition}
Let us first split each vector field into two polarizations: one parallel and the other normal to the $z \rho$-plane, denoted by the subscripts $\parallel$ and $\perp$, respectively:
\begin{flalign}
\mathbf{E}' &= \mathbf{E}'_{\parallel}+\mathbf{E}'_{\perp},\\
\mathbf{D}' &= \mathbf{D}'_{\parallel}+\mathbf{D}'_{\perp},\\
\mathbf{B}' &= \mathbf{B}'_{\parallel}+\mathbf{B}'_{\perp},\\
\mathbf{H}' &= \mathbf{H}'_{\parallel}+\mathbf{H}'_{\perp},
\end{flalign}
where
\begin{flalign}
\mathbf{E}'_{\parallel}=\hat{z}E'_{z} + \hat{\rho}E'_{\rho}, \quad
\mathbf{E}'_{\perp}=\hat{\phi}E'_{\phi}.
\end{flalign}
and likewise for $\mathbf{D}'$, $\mathbf{B}'$, and $\mathbf{H}'$.

By inserting these expressions into Maxwell's equations, we can effectively separate Maxwell's equations into two set of equations: one involving the TE\textsuperscript{$\phi$} field components
 $\left\{E'_{z},E'_{\rho},B'_{\phi}\right\}$ and  $\left\{D'_{z},D'_{\rho},H'_{\phi}\right\}$ and the other involving the TM\textsuperscript{$\phi$} field components $\left\{E'_{\phi},B'_{z},B'_{\rho}\right\}$ and
  $\left\{D'_{\phi},H'_{z},H'_{\rho}\right\}$.

\subsection{Motivation for $EB$ and $DH$ discretizations}\label{Sec:Motivation_EB_DH}
Now, we are ready to discretize the rescaled EM fields for each polarization. As noted, the discretization uses an unstructured finite-element mesh based on triangular elements.
We discretize $\mathbf{E}'$ and $\mathbf{B}'$ fields for the TE\textsuperscript{$\phi$} polarization and the $\mathbf{D}'$ and $\mathbf{H}'$ fields for the TM\textsuperscript{$\phi$} polarization.
There are two important motivations behind such choice:
\begin{itemize}
\item First, it obviates the singularity present in $\bar{\bar{\epsilon}}'$ and $\bar{\bar{\mu}}'$ for $\rho=0$ along the $z$ axis.
Specifically, from \eqref{eqn:eff_permittivity} and \eqref{eqn:rescale_GL} (and likewise \eqref{eqn:eff_permeability} and \eqref{eqn:rescale_GLM} for the magnetic field case), the rescaled fields obey the relation
\begin{flalign}
D'_{\phi} = \frac{\epsilon_0}{\rho} E'_{\phi}
\label{eqn:DE_phi_eqn}
\end{flalign}
in vacuum.
When solving for $E'_{\phi}$ first and obtaining $D'_{\phi}$ through \eqref{eqn:DE_phi_eqn}, one cannot have an accurate and stable numerical solution of $D’_\phi$ since the discrete Hodge matrix encoding $\epsilon_0/\rho$ has a poor convergence with respect to the mesh refinement.
On the other hand, when solving for $D’_\phi$ first, and then obtaining $E’_\phi$ from $E'_{\phi} = \frac{\rho}{\epsilon_0} D'_{\phi}$, one can obtain a stable numerical solution of $E’_\phi$ since the discrete Hodge matrix involving $\rho/\epsilon_0$ is now well-defined and singularity-free.

\item Second, this choice makes it possible to exploit a duality between the two different polarizations: the discretization of the TE\textsuperscript{$\phi$} polarization case is equivalent to the discretization of the TM\textsuperscript{$\phi$} case with the substitutions below:
\begin{flalign}
\label{eqn:duality}
\mathbf{E} \rightarrow \mathbf{H}, \quad \mathbf{B} \rightarrow \mathbf{D}, \quad \epsilon \rightarrow \mu, \quad \mu \rightarrow \epsilon.
\end{flalign}
This duality makes it possible to reuse the same computer code, with only minor adaptations, for both polarizations.
\end{itemize}

\subsection{TE\textsuperscript{$\phi$} field solver}
In the language of the exterior calculus of differential forms, $\mathbf{E}'_{\parallel}$ and $\mathbf{B}'_{\perp}$ are represented as 1- and 2-forms, respectively~\cite{Teixeira1999Lattice,HE20061}. From a physical standpoint, this difference reflects the distinct type of continuity conditions (tangential versus normal) of these fields.
From a strictly mathematical standpoint, this distinction reflects the different functional spaces these fields belong to. $\mathbf{E}'_{\parallel}$ belongs to the curl-conforming Sobolev space $H(\text{curl})$ and $\mathbf{B}'_{\perp}$ to the div-conforming Sobolev space $H(\text{div})$~\cite{doi:10.1137/08073901X,monkFEM,campos2016constructing,pinto2022semi,Ramachandran2023Review}. On a triangular mesh, 1- and 2-forms are expanded in terms of their natural discrete interpolants: Whitney 1- and 2-forms, respectively~\cite{Bossavit1988Whitney,Teixeira1999Lattice,HE20061}:
\begin{flalign}\label{eqn:exp1}
\mathbf{E}'_{\parallel}(\mathbf{r}_{\parallel},t)
&=
\sum_{i=1}^{N_{1}}
e_{\parallel,i}(t) \,
\mathbf{W}^{(1)}_{i}(\mathbf{r}_{\parallel}),
\end{flalign}
\begin{flalign}\label{eqn:exp2}
\mathbf{B}'_{\perp}(\mathbf{r}_{\parallel},t)
&=
\sum_{i=1}^{N_{2}}
b_{\perp,i}(t) \,
\mathbf{W}^{(2)}_{i}(\mathbf{r}_{\parallel}),
\end{flalign}
where $\mathbf{r}_{\parallel}$ denotes the position vector on the $z \rho$-plane, $N_{1}$ and $N_{2}$ are the number of edges and faces (cells, for a 2D mesh discretizing the $z \rho$ plane) in the finite element mesh, respectively, and
$\mathbf{W}^{(p)}_{i}(\mathbf{r}_{\parallel})$ denotes the vector proxy of the Whitney $p$-form associated with the $i$-th $p$-cell of the primal mesh ($p=1$: edges; $p=2$: faces).
Explicit expressions for $\mathbf{W}^{(p)}_{i}(\mathbf{r}_{\parallel})$ are found in~\cite{Bossavit1988Whitney,Teixeira1999Lattice,HE20061,GILLETTE20111213,LOHI2021113520}, for example.
The expansions in eqs. (\ref{eqn:exp1}) and  (\ref{eqn:exp2}) enforce the appropriate field continuity rules across the finite element mesh elements.

By substituting eq.~(\ref{eqn:exp1}) and eq.~(\ref{eqn:exp2})
into eq.~(\ref{eqn:rescale_FL}), we get
\begin{flalign} \label{eqn:FL}
\sum_{i=1}^{N_{2}}
\left( \frac{\partial  }{\partial t} b_{\perp,i}(t) \right)
\mathbf{W}^{(2)}_{i}=-\nabla' \times \sum_{j=1}^{N_{1}}
e_{\parallel,j}(t) \,
\mathbf{W}^{(1)}_{j}= -\sum_{j=1}^{N_{1}} \sum_{i=1}^{N_{2}} C_{ij}  \, e_{\parallel,j}(t) \, \mathbf{W}^{(2)}_{i}
\end{flalign}
where the last equality is a structural property of the Whitney forms~\cite{Teixeira1999Lattice,LOHI2021113520,10.1007/978-3-319-01601-6,pinto2016electromagnetic} with $C_{ij}$ being the elements of the so called incidence matrix with entries $\{-1,0,+1\}$~\cite{Teixeira1999Lattice}. The incidence matrix encodes the discrete (primed) curl operator distilled from the metric (or more precisely, the coboundary operator on the mesh~\cite{gross2004electromagnetic}).
From eq.~(\ref{eqn:FL}), we immediately obtain
\begin{flalign} \label{eqn:FL2}
\frac{\partial}{\partial t}  b_{\perp,i}(t) = -\sum_{j=1}^{N_{1}} C_{ij}  \, e_{\parallel,j}(t)
\end{flalign}
for $i=1 \ldots N_{2}$. It is convenient to rewrite eq. (\ref{eqn:FL2}) in a more compact form as
\begin{flalign} \label{eqn:FLmatrix}
\frac{\partial  }{\partial t} \mathbf{b_{\perp}}(t) = -\bar{\mathbf{C}} \cdot  \mathbf{e}_{\parallel}(t)
\end{flalign}
where $\bar{\mathbf{C}}$ is the incidence matrix, $\mathbf{b}_{\perp}(t)$ represents the column vector that collects all the degrees of freedom (DoF) of $\mathbf{B}_{\perp}(\mathbf{r}_{\perp},t)$, i.e.,
\begin{flalign}
\mathbf{b}_{\perp}(t)
=
\left[
b_{\perp,1}(t),b_{\perp,2}(t),
\cdots,
b_{\perp,N_{2}}(t)
\right]^{T},
\end{flalign}
and similarly
\begin{flalign}
\mathbf{e}_{\parallel}(t)
=
\left[
e_{\parallel,1}(t),e_{\parallel,2}(t),
\cdots,
e_{\parallel,N_{1}}(t)
\right]^{T}.
\end{flalign}
Next, we rewrite eq.~(\ref{eqn:rescale_AL}) as
\begin{flalign}
\label{eqn:me2rev}
\nabla' \times \left( \bar{\bar{\mu}}'^{-1} \cdot \mathbf{B}'  \right) &= \frac{\partial}{\partial t} \left( \bar{\bar{\epsilon}}' \cdot \mathbf{E}' \right)+\mathbf{J}',
\end{flalign}
and substitute eqs.~(\ref{eqn:exp1}) and (\ref{eqn:exp2})
into eq.~(\ref{eqn:me2rev}) to get
\begin{flalign}
\label{eqn:AL}
\nabla' \times \left( \bar{\bar{\mu}}'^{-1} \cdot
\sum_{i=1}^{N_{2}}
b_{\perp,i}(t) \,
\mathbf{W}^{(2)}_{i} \right) \approx \frac{\partial}{\partial t} \left( \bar{\bar{\epsilon}}' \cdot  \sum_{j=1}^{N_{1}}
e_{\parallel,j}(t) \,
\mathbf{W}^{(1)}_{j} \right)
 +\mathbf{J}',
\end{flalign}
At this point, it should be recognized that the above equation is only an approximation because
$\nabla' \times \mathbf{W}^{(2)}_{i} \not\subset span \{\mathbf{W}^{(1)}_{j}; j=1,\ldots,N_1 \}$,
and therefore we cannot invoke a structural property such as used in eq.~(\ref{eqn:FL}). Instead, an equality is enforced through a Galerkin projection, i.e., by taking the inner product of both sides of eq.~(\ref{eqn:AL}) with Whitney 1-forms $\mathbf{W}^{(1)}_{k}$, $k=1,\ldots,N_1$, so that
\begin{flalign}
\label{eqn:ALweak}
&\iint_{\Omega} \, \mathbf{W}^{(1)}_{k} \cdot \nabla' \times \left( \bar{\bar{\mu}}'^{-1} \cdot \sum_{i=1}^{N_{2}}
b_{\perp,i}(t) \,
\mathbf{W}^{(2)}_{i} \right) dA = \nonumber \\
&
 \iint_{\Omega} \, \mathbf{W}^{(1)}_{k} \cdot \left[ \frac{\partial}{\partial t} \left( \bar{\bar{\epsilon}}' \cdot \sum_{j=1}^{N_{1}}
e_{\parallel,j}(t) \,
\mathbf{W}^{(1)}_{j} \right) \right] dA
 + \iint_{\Omega} \,  \mathbf{W}^{(1)}_{k} \cdot \mathbf{J}' \, dA,
\end{flalign}
for $k=1,\ldots,N_1$, where $dA = dz d\rho$ is the area element in the $z \rho$ plane and $\Omega$ represent the computational domain in that plane.

By using the identity $\nabla \cdot (\mathbf{A} \times \mathbf{B}) = (\nabla \times \mathbf{A} ) \cdot \mathbf{B} - (\nabla \times \mathbf{B} ) \cdot \mathbf{A}$ and assuming that the divergence term integrates to zero\footnote{Note that the associated volume integral is equal, by virtue of the Stokes' theorem, to a surface flux integral on the exterior boundary of the computational domain. This flux integral evaluates to zero given the boundary conditions. In our case, we assume perfect magnetic conductor (PMC) boundary conditions on the z-axis for TE\textsuperscript{$\phi$} polarized fields such that the divergence term including the tangential magnetic field vanishes.
We employ PMC boundary conditions backing the PML in the radial direction, the divergence term goes to zero there as well.
For the lateral boundaries with periodic boundary conditions, the divergence term is zero because the terms on left and right walls cancel out each other in accordance with the Poynting theorem.}, the left-hand side of eq.~(\ref{eqn:ALweak}) can be written as
\begin{flalign}
\label{eqn:ALdis}
&\iint_{\Omega}  \, \left( \nabla' \times \mathbf{W}^{(1)}_{k} \right) \cdot \left( \bar{\bar{\mu}}'^{-1} \cdot
\sum_{i=1}^{N_{2}}
b_{\perp,i}(t) \,
\mathbf{W}^{(2)}_{i} \right)  dA = \nonumber \\
&\iint_{\Omega}  \,
\left( \sum_{j=1}^{N_2} C_{jk}
 \mathbf{W}^{(2)}_{j} \right) \cdot \left( \bar{\bar{\mu}}'^{-1} \cdot
\sum_{i=1}^{N_{2}}
b_{\perp,i}(t) \,
\mathbf{W}^{(2)}_{i} \right)  dA = \nonumber \\
& \sum_{j=1}^{N_2} C_{jk} \sum_{i=1}^{N_{2}}
\left(  \iint_{\Omega}  \,
 \mathbf{W}^{(2)}_{j}  \cdot
 \bar{\bar{\mu}}'^{-1} \cdot
\mathbf{W}^{(2)}_{i} \, dA \right) b_{\perp,i}(t)
\end{flalign}
for $k=1,\ldots,N_1$ and where in the second step we used the same structural property of Whitney forms as in eq.~(\ref{eqn:FL}).
By defining
\begin{flalign}
 \left[\star_{\mu^{-1}}\right]_{i,j} &\equiv  \iint_{\Omega}  \,
 \mathbf{W}^{(2)}_{i}  \cdot
 \bar{\bar{\mu}}'^{-1} \cdot
\mathbf{W}^{(2)}_{j} \, dA
\nonumber \\
&=
\iint_{\Omega}  \,
\left(\hat{\phi}\cdot\mathbf{W}^{(2)}_{i}\right)
{\bar{\bar{\mu}}'}^{-1}_{\phi,\phi}
\left(\hat{\phi}\cdot\mathbf{W}^{(2)}_{j}\right) \, dA,
\label{eqn:TE_Hodge_mu}
\end{flalign}
the last expression in $(\ref{eqn:ALdis})$ can be rewritten in a more compact form as
$ \bar{\mathbf{C}}^T \cdot
 \left[\star_{\mu^{-1}}\right] \cdot
\mathbf{b}_{\perp}(t) $ and by applying similar steps to the right-hand side of eq.~(\ref{eqn:ALweak}), we obtain
\begin{flalign}\label{eqn:ALsemi}
\bar{\mathbf{C}}^T \cdot
 \left[\star_{\mu^{-1}}\right] \cdot
\mathbf{b}_{\perp}(t) = \left[\star_{\epsilon}\right] \cdot \frac{\partial}{\partial t}  \mathbf{e}_{\parallel}(t) + \mathbf{j}_{\parallel}(t),
\end{flalign}
with
\begin{flalign}
\left[\star_{\epsilon}\right]_{i,j}
&\equiv
\iint_{\Omega}\,
\mathbf{W}^{(1)}_{i}
\cdot
\bar{\bar{\epsilon}}' \cdot
\mathbf{W}^{(1)}_{j} \, dA,
\nonumber \\
&=
\iint_{\Omega}\,
\left[
\begin{matrix}
\hat{\rho}\cdot\mathbf{W}^{(1)}_{i} \\ \hat{z}\cdot\mathbf{W}^{(1)}_{i}
\end{matrix}
\right]^{T}
\cdot
\left[
\begin{matrix}
\bar{\bar{\epsilon}}'_{\rho,\rho} & \bar{\bar{\epsilon}}'_{\rho,z} \\
\bar{\bar{\epsilon}}'_{z,\rho} & \bar{\bar{\epsilon}}'_{z,z} \\
\end{matrix}
\right]
\cdot
\left[
\begin{matrix}
\hat{\rho}\cdot\mathbf{W}^{(1)}_{j} \\ \hat{z}\cdot\mathbf{W}^{(1)}_{j}
\end{matrix}
\right]
\, dA,
\label{eqn:TE_Hodge_eps}
\end{flalign}
\begin{flalign}
\label{eqn:Jproj}
{j}_{\parallel,i}(t) \equiv  \iint_{\Omega} \,  \mathbf{W}^{(1)}_{i} \cdot \mathbf{J}'(t) \, dA.
\end{flalign}
As observed in the expressions in \eqref{eqn:TE_Hodge_mu} and \eqref{eqn:TE_Hodge_eps} with \eqref{eqn:eff_permeability} and \eqref{eqn:eff_permittivity}, both discrete Hodge matrices $\left[\star_{\mu}^{-1}\right]$ and $\left[\star_{\epsilon}\right]$ for TE\textsuperscript{$\phi$} polarized fields are devoid from any singularity at $\rho=0$.
Additional details are discussed in \ref{app:singularity}.

Finally, by using a central (leap-frog) finite-difference approximation for the time derivatives in eqs.~(\ref{eqn:FLmatrix}) and (\ref{eqn:ALsemi}), we obtain the following field update equations~\cite{Donderici2008Mixed,Kim2011Parallel}:
\begin{flalign}
\label{eqn:update_TE_1}
\mathbf{b}^{n+\frac{1}{2}}_{\perp}
&= \mathbf{b}^{n-\frac{1}{2}}_{\perp}
- \Delta t\bar{\mathbf{C}}\cdot\mathbf{e}^{n}_{\parallel},
\\
\label{eqn:update_TE_2}
\left[\star_{\epsilon}\right] \cdot \mathbf{e}^{n+1}_{\parallel}
&= \left[\star_{\epsilon}\right] \cdot \mathbf{e}^{n}_{\parallel}
+ \Delta t
\left(\bar{\mathbf{C}}^{T}\cdot\left[\star_{\mu^{-1}}\right]\cdot\mathbf{b}^{n+\frac{1}{2}}_{\perp}
- \mathbf{j}^{n+\frac{1}{2}}_{\parallel}\right).
\end{flalign}
where $\Delta t$ is the time-step increment, $n$ is the time-step index, and
\begin{flalign}
\mathbf{b}_{\perp}^{n+\frac{1}{2}}
&=
\left[
b_{\perp,1}^{n+\frac{1}{2}},b_{\perp,2}^{n+\frac{1}{2}},
\cdots,
b_{\perp,N_{2}}^{n+\frac{1}{2}}
\right]^{T}, \\
\mathbf{e}_{\parallel}^{n}
&=
\left[
e_{\parallel,1}^{n},e_{\parallel,2}^{n},
\cdots,
e_{\parallel,N_{1}}^{n}
\right]^{T}, \\
\mathbf{j}_{\parallel}^{n+\frac{1}{2}}
&=
\left[
j_{\parallel,1}^{n+\frac{1}{2}},j_{\parallel,2}^{n+\frac{1}{2}},
\cdots,
j_{\parallel,N_{1}}^{n+\frac{1}{2}}
\right]^{T}.
\end{flalign}
Recall that, in our context, the current density term $\mathbf{j}^{n+\frac{1}{2}}_{\parallel}$ arises from the movement of charges in the plasma medium.
The current density term is discussed in detail in Section~\ref{Sec:Parallel_Scatter}.
The matrices $\left[\star_{\epsilon}\right]$ and $\left[\star_{\mu^{-1}}\right]$ are symmetric positive-definite and very sparse. In the finite element parlance, they are often termed mass matrices or sometimes discrete Hodge star operators~\cite{He2007Differential,GILLETTE20111213,Teixeira2013Differential}. In that form, eq.~(\ref{eqn:update_TE_2}) requires a sparse (iterative) linear solver at  each time step because $\left[\star_{\epsilon}\right]$ multiplies the unknowns in the l.h.s. of that equation. However, given its structure, $\left[\star_{\epsilon}\right]^{-1}$ can be accurately represented by a sparse approximate inverse (SPAI) computed in embarrassingly parallel fashion~\cite{Kim2011Parallel,Na2016Local}. The SPAI can be obtained once and for all before the start of the time stepping algorithm, thus obviating the need for a linear solver at each time step.

\begin{figure}
\centering
\subfloat[Finite-element (primal) and dual meshes in the $z \rho$-plane.]
{\includegraphics[width=.75\linewidth]{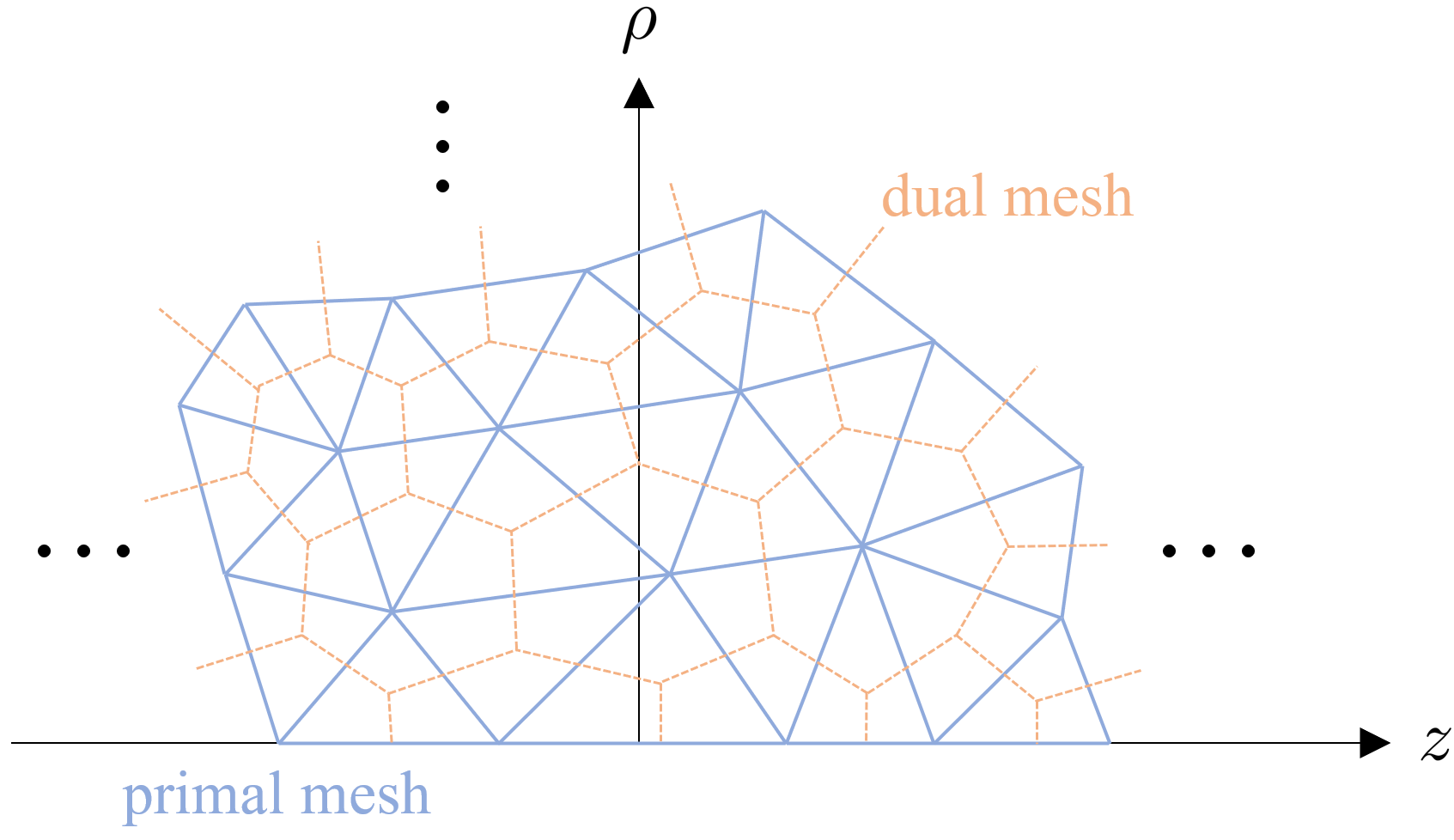}}
\\
\subfloat[Primal mesh discretization of TE\textsuperscript{$\phi$} polarization components.]
{\includegraphics[width=.35\linewidth]{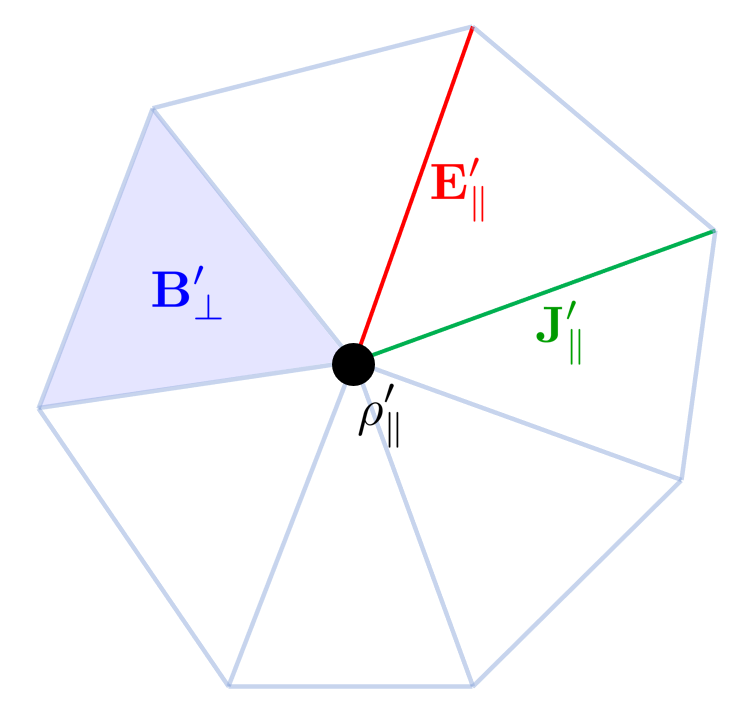}}
\quad\quad\quad
\subfloat[Dual mesh discretization of TE\textsuperscript{$\phi$} polarization components.]
{\includegraphics[width=.35\linewidth]{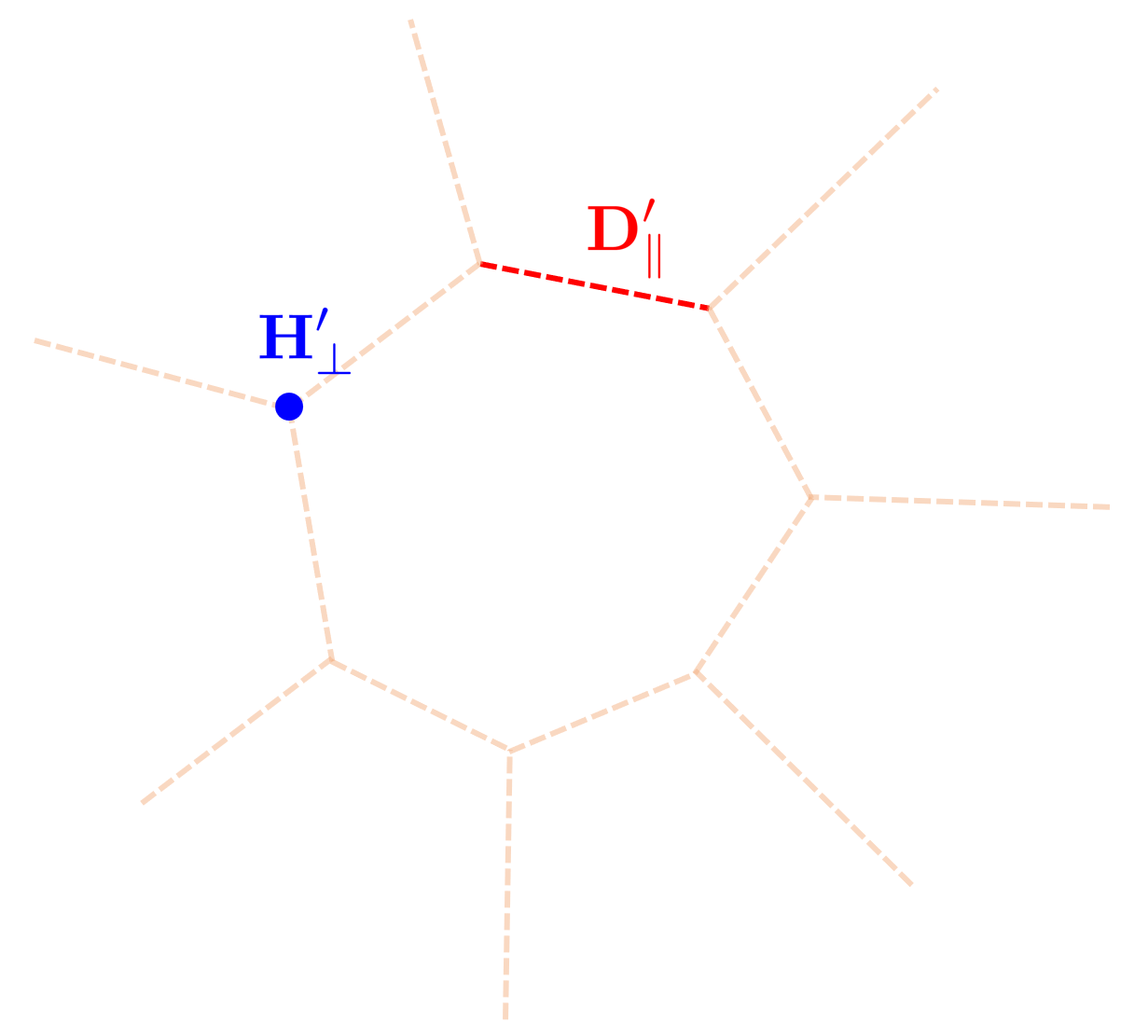}}
\\
\subfloat[Primal mesh discretization of TM\textsuperscript{$\phi$} polarization components.]
{\includegraphics[width=.35\linewidth]{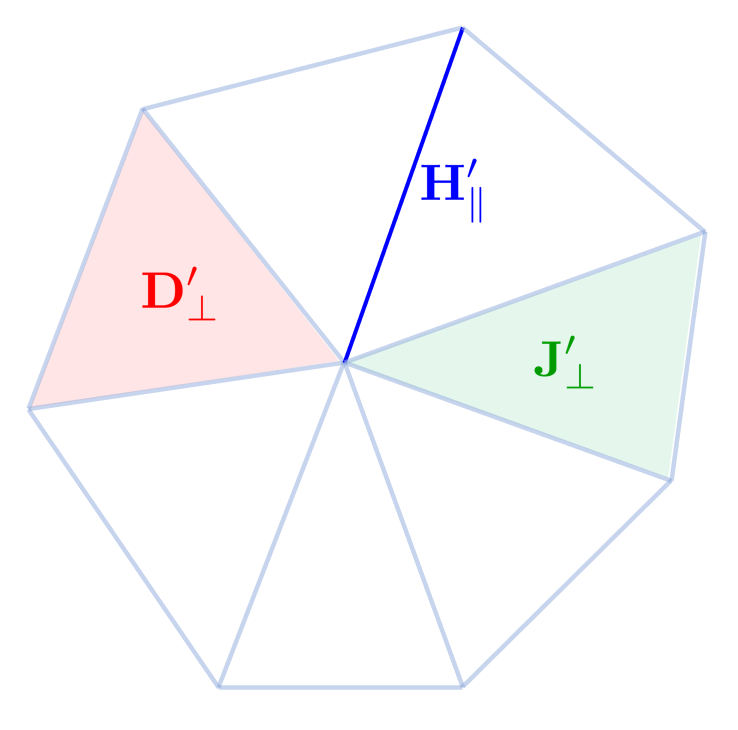}}
\quad\quad\quad
\subfloat[Dual mesh discretization of TM\textsuperscript{$\phi$} polarization components.]
{\includegraphics[width=.35\linewidth]{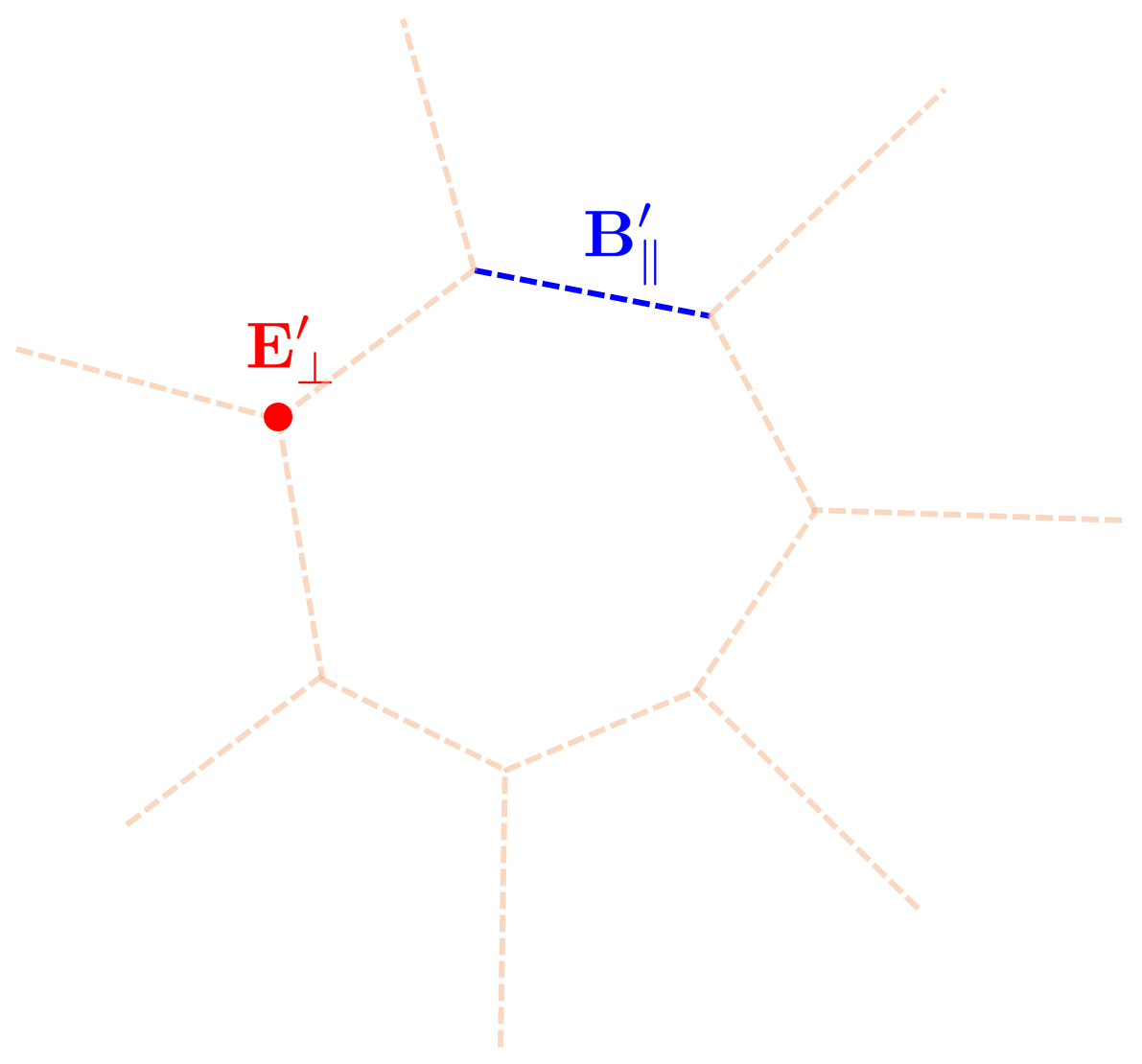}}
\caption{Localization of the field degrees of freedom on the primal and dual meshes.
}
\label{fig:primal_dual_meshes}
\end{figure}

{\it Remark 1:} The coefficients $\mathbf{e}_{\parallel}(t)$ and $\mathbf{b}_{\perp}(t)$ constitute the
the field unknowns and are obtained by the finite element solver. On the unstructured finite-element mesh,
$\mathbf{e}_{\parallel}(t)$
 are defined on the {\it edges} and $\mathbf{b}_{\perp}(t)$ on the {\it faces} (triangular cells in 2D) of the mesh, as illustrated in Fig.~\ref{fig:primal_dual_meshes}. This choice is informed by the exterior calculus of differential forms, where $\mathbf{E}'_{\parallel}$ is represented as a 1-form defined on line integrals (hence on edges as 1-dimensional objects) and $\mathbf{B}'_{\perp}$ is represented as a 2-form defined on surface integrals (hence on faces as 2-dimensional objects). This localization rule is important to ensure a consistent spatial discretization devoid of spurious modes and instabilities~\cite{bochev2006principles,arnold2010finite}. For more details on the use of exterior calculus to inform the consistent (or ``structure-preserving'') spatial finite element discretization of vector and tensor fields, the reader is referred to~\cite{Teixeira1999Lattice,PIER2001,hirani2003discrete,bochev2006principles,Teixeira2013Differential,Teixeira2014Lattice,arnold2018finite,9599150} and references therein. An early reference addressing discrete exterior calculus in the context of electromagnetic PIC algorithms is~\cite{squire2012geometric}.

{\it Remark 2:} Fig.~\ref{fig:primal_dual_meshes} depicts two meshes: the ``primal'' mesh (triangular elements) and the corresponding ``dual'' mesh (polygonal elements).
The primal mesh is the actual mesh used by the finite element solver. The polygonal dual mesh is not used in the actual computations and it may be understood as a device to indicate the localization of the degrees of freedom of the remaining components not used in the field solver stage~\cite{Teixeira1999Lattice,GILLETTE20111213,squire2012geometric,9599150}. For example, the degrees of freedom of $\mathbf{H}'_{\perp}$  (0-form in 2-dimensions) are associated with the nodes of the dual mesh
and those of $\mathbf{D}'_{\parallel}$ (1-form in 2-dimensions) with the edges of the dual mesh. Observe from Fig.~\ref{fig:primal_dual_meshes}(a) that nodes of the dual mesh are in a bijective (1:1) correspondence  with faces (triangular cells) of the primal mesh and the edges of the dual mesh are in a bijective correspondence with the edges of the primal mesh~\cite{PhysRevE.61.3174}. These are instances of the general result that $k$-dimensional cells of the dual mesh are in bijective correspondence with the $(n-k)$-dimensional cells of the primal mesh, where $n$ is the dimension of space. This bijection implies that, in the discrete setting, the same number of unknowns represent
$\mathbf{D}$ and $\mathbf{E}$ or $\mathbf{B}$ and $\mathbf{H}$. This geometric duality is also color-coded in Fig.~\ref{fig:primal_dual_meshes}.

{\it Remark 3:}
It is also to be stressed that we regard transverse components of $\mathbf{E}$ and azimuthal components of $\mathbf{B}$ fields (ordinary forms) as primal quantities for TE\textsuperscript{$\phi$} polarization whereas azimuthal components of $\mathbf{D}$ and transverse components of $\mathbf{H}$ fields (twisted forms) are considered to be primal quantities for TM\textsuperscript{$\phi$} polarization.
Note again that both polarizations are discretized on the finite element mesh.
In typical mixed finite-element time-domain schemes, $\mathbf{E}$ and $\mathbf{B}$ fields are assumed to be primal quantities \cite{HE20061,He2007Differential,Donderici2008Mixed,Kim2011Parallel,Moon2015Exact,Na2016Local,Na2017Axisymmetric,Na2019Finite}.
However, this is not strictly necessarily and one can choose for $\mathbf{D}$ and $\mathbf{H}$ instead to be discretized on the primal mesh \cite{osti_950065,campospinto:hal-01303852,pinto2016electromagnetic,Na2019Polynomial}.
For example, for TE\textsuperscript{$\phi$} polarized fields in $z\rho$ plane, the $\mathbf{D}$ field is represented as a (twisted) 2-form with dergees of freedom associated with the area elements of the primal mesh and expanded using Whitney 2-forms. On the other hand, the $\mathbf{H}$ field is represented as a (twisted) 1-form with degrees of freedom associated with the edges of the primal mesh and expanded using of Whitney 1-forms.
To stress, the reason why we deal with the twisted forms on the primal mesh for the TM\textsuperscript{$\phi$} case is to obtain a singularity-free Hodge matrix elements along the z-axis. Further details on this aspects are presented in Section \ref{Sec.TM_Phi_FETD}.

\begin{figure}
\centering
\subfloat[Primal mesh discretization.]
{\includegraphics[width=.45\linewidth]{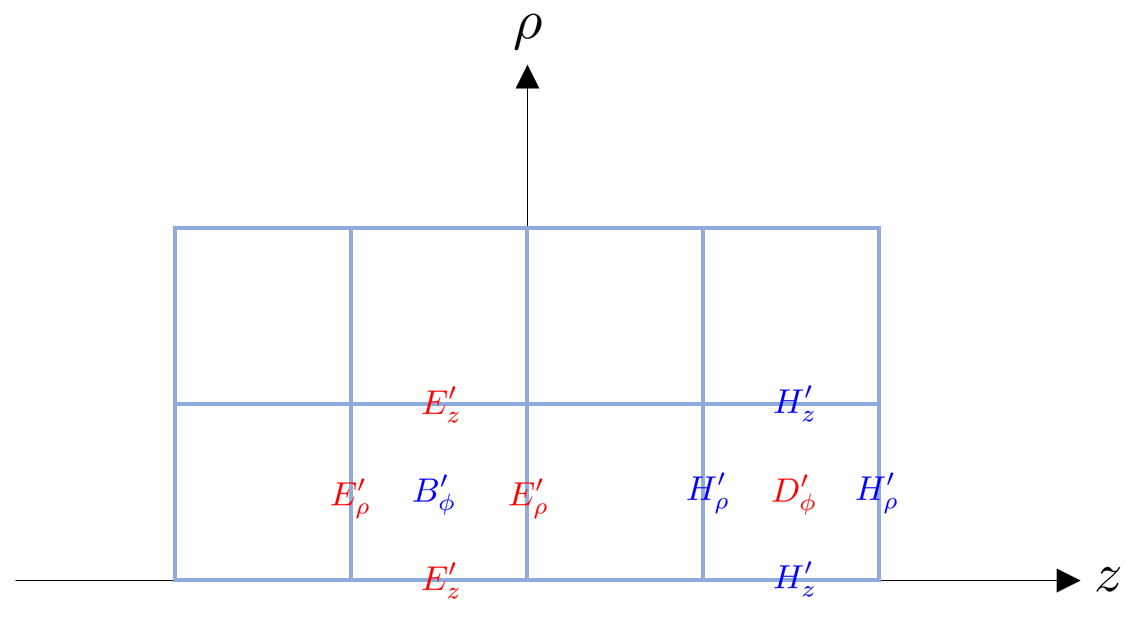}}
\quad
\subfloat[Dual mesh discretization.]
{\includegraphics[width=.45\linewidth]{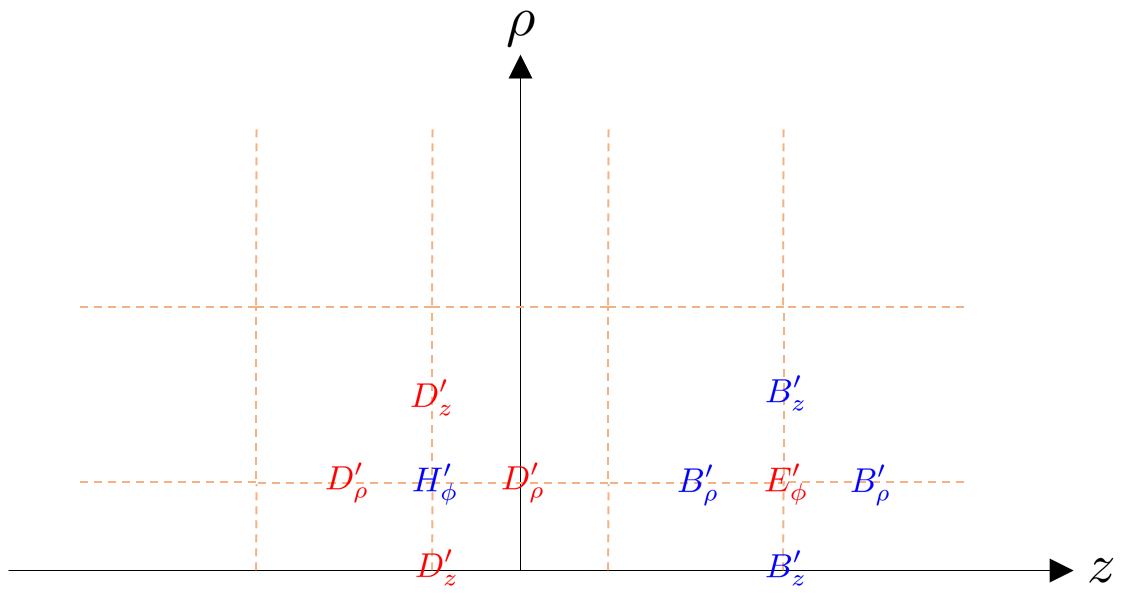}}
\caption{
The regular-mesh counterpart of the field discretization on an irregular mesh.
}
\label{fig:primal_dual_meshes_regular_mesh}
\end{figure}
Fig. \ref{fig:primal_dual_meshes_regular_mesh} illustrates the regular-mesh counterpart for the better understanding of the field discretization on an irregular mesh.
In the same primal mesh, we discretize the transverse part of $\mathbf{E}'$ and the azimuthal part of $\mathbf{B}'$ for TE\textsuperscript{$\phi$} polarized fields and the azimuthal part of $\mathbf{D}'$ and the transverse part of $\mathbf{H}'$ for TM\textsuperscript{$\phi$} polarized fields.

\subsection{TM\textsuperscript{$\phi$} field solver} \label{Sec.TM_Phi_FETD}
As mentioned earlier, for TM\textsuperscript{$\phi$} polarized fields, we regard electric flux density $\mathbf{D}'_{\perp}$ and magnetic field intensity $\mathbf{H}'_{\parallel}$ as primal quantities instead so that $\mathbf{D}'_{\perp}$ and $\mathbf{H}'_{\parallel}$ are discretized on the primal mesh.

Since $\mathbf{D}_{\perp}$ is a 2-form and $\mathbf{H}_{\parallel}$ is a 1-form~\cite{Teixeira1999Lattice,He2007Differential}, they can be expanded by Whitney 2- and 1-forms respectively, i.e.,
\begin{flalign} \label{eq:exp3}
\mathbf{D}'_{\perp}(\mathbf{r}_{\parallel},t)
&=
\sum_{i=1}^{N_{2}}
d_{\perp,i}(t)
\mathbf{W}^{(2)}_{i}(\mathbf{r}_{\parallel}),
\end{flalign}
\begin{flalign} \label{eq:exp4}
\mathbf{H}'_{\parallel}(\mathbf{r}_{\parallel},t)
&=
\sum_{i=1}^{N_{1}}
h_{\parallel,i}(t)
\mathbf{W}^{(1)}_{i}(\mathbf{r}_{\parallel}),
\end{flalign}
which are dual to the TE\textsuperscript{$\phi$} polarization case~\cite{He2007Differential}.

\begin{figure}
\centering
\includegraphics[width=.75\linewidth]
{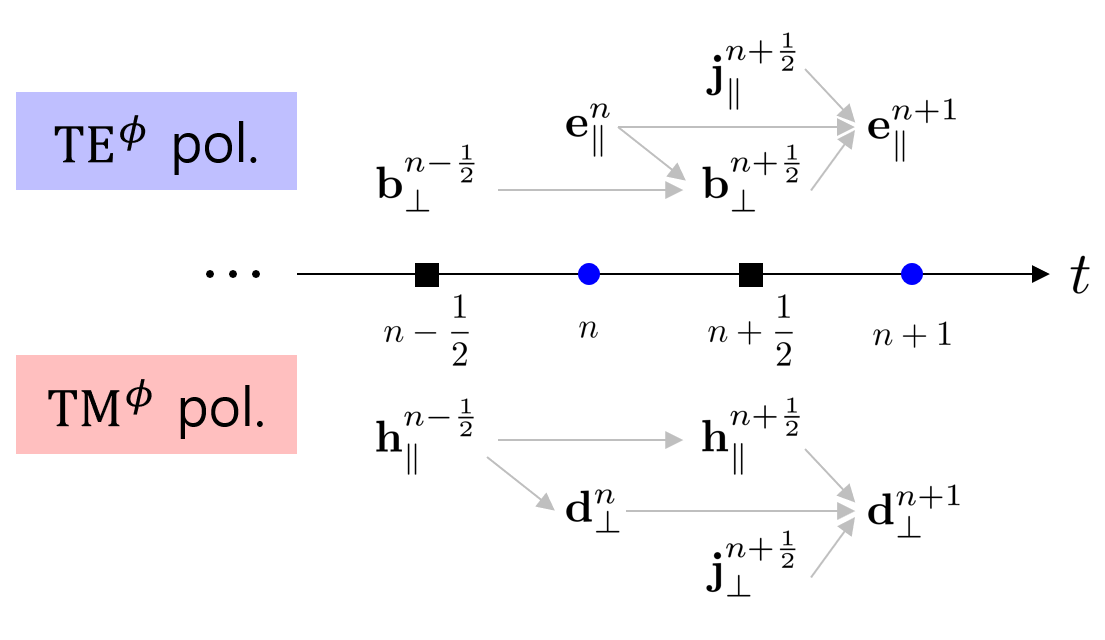}
\caption{Causal field dependency after time discretization, for both polarizations.}
\label{fig:time_discretization}
\end{figure}

Following analogous steps as those shown in the TE\textsuperscript{$\phi$} case, we obtain the following field update equations in the TM\textsuperscript{$\phi$} case :
\begin{flalign}
\label{eqn:TMupdate_1}
\left[\star_{\mu}\right] \cdot \mathbf{h}^{n+\frac{1}{2}}_{\parallel}
&= \left[\star_{\mu}\right] \cdot \mathbf{h}^{n-\frac{1}{2}}_{\parallel}
- \Delta t
\bar{\mathbf{C}}^{T}\cdot
\left[\star_{\epsilon^{-1}}\right]\cdot\mathbf{d}^{n}_{\perp},
\\
\label{eqn:TMupdate_2}
\mathbf{d}^{n+1}_{\perp}
&= \mathbf{d}^{n}_{\perp}
+ \Delta t\left(\bar{\mathbf{C}}\cdot\mathbf{h}^{n+\frac{1}{2}}_{\parallel}
- \mathbf{j}^{n+\frac{1}{2}}_{\perp}\right),
\end{flalign}
where boldface lowercases $\mathbf{d}_{\perp}^{n}$, $\mathbf{h}_{\parallel}^{n+\frac{1}{2}}$, and  $\mathbf{j}_{\perp}^{n+\frac{1}{2}}$ are column vectors collecting DoFs of $\mathbf{D}'_{\perp}$, $\mathbf{H}'_{\parallel}$, and $\mathbf{J}'_{\perp}$, respectively. The new discrete Hodge matrices $\left[\star_{\epsilon^{-1}}\right]$ and $\left[\star_{\mu}\right]$ in the TM\textsuperscript{$\phi$} case are given by
\begin{flalign}
\left[\star_{\epsilon^{-1}}\right]_{i,j}
&=
\iint_{\Omega}
\mathbf{W}^{(2)}_{i}
\cdot
\bar{\bar{\epsilon}}'^{-1} \cdot
\mathbf{W}^{(2)}_{j}
\nonumber \\
&=
\iint_{\Omega}  \,
\left(\hat{\phi}\cdot\mathbf{W}^{(2)}_{i}\right)
{\bar{\bar{\epsilon}}'}^{-1}_{\phi,\phi}
\left(\hat{\phi}\cdot\mathbf{W}^{(2)}_{j}\right) \, dA,
\label{eqn:TM_Hodge_eps}
\, dA,
\\
\left[\star_{\mu}\right]_{i,j}
&=
\iint_{\Omega}
\mathbf{W}^{(1)}_{i}
\cdot
\bar{\bar{\mu}}' \cdot
\mathbf{W}^{(1)}_{j}
\, dA
\nonumber \\
&=
\iint_{\Omega}\,
\left[
\begin{matrix}
\hat{\rho}\cdot\mathbf{W}^{(1)}_{i} \\ \hat{z}\cdot\mathbf{W}^{(1)}_{i}
\end{matrix}
\right]^{T}
\cdot
\left[
\begin{matrix}
\bar{\bar{\mu}}'_{\rho,\rho} & \bar{\bar{\mu}}'_{\rho,z} \\
\bar{\bar{\mu}}'_{z,\rho} & \bar{\bar{\mu}}'_{z,z} \\
\end{matrix}
\right]
\cdot
\left[
\begin{matrix}
\hat{\rho}\cdot\mathbf{W}^{(1)}_{j} \\ \hat{z}\cdot\mathbf{W}^{(1)}_{j}
\end{matrix}
\right] \, dA.
\label{eqn:TM_Hodge_mu}
\end{flalign}
Similarly, the discrete Hodge matrices for TM\textsuperscript{$\phi$} fields can be also well-defined since they do not present any singularity at $\rho=0$, as observed in \eqref{eqn:TM_Hodge_mu} and \eqref{eqn:TM_Hodge_eps} with \eqref{eqn:eff_permittivity} and \eqref{eqn:eff_permeability}.
The same observations made for eq.~(\ref{eqn:update_TE_2}) in regards to the SPAI solution can be made for eq.~(\ref{eqn:TMupdate_1}) as well.
 Note that a similar finite-element discretization based on $\mathbf{D}$ and $\mathbf{H}$ was considered before in~\cite{Na2019Polynomial}.
Fig. \ref{fig:time_discretization} depicts the resulting field dependency in the time-stepping algorithm, for both polarizations.

\subsection{Field boundary conditions on the symmetry axis}
For BOR problems where the line $\rho=0$ (symmetry axis) is part of the solution domain, it becomes necessary to treat the field there by means of appropriate boundary conditions in the $z \rho$ plane.
The boundary conditions at $\rho=0$ are mode-dependent (eigenmodes along azimuth indexed by $m$) and should reflect the cylindrical coordinate system singularity and the related degeneracy of the cylindrical unit vectors  $\hat{\rho}$ and $\hat{\phi}$ there.
When $m = 0$, there is no field variation along azimuth and, in the absence of charges at $\rho=0$, both azimuthal and radial field components are zero at $\rho=0$.
As a result, when $m = 0$, the boundary $\rho=0$ can be represented as a perfect magnetic conductor (PMC) for the TE\textsuperscript{$\phi$} field since $E'_{z}\neq 0, E'_{\rho}={B'}_{\phi}=0$.\footnote{Note that the perfect electric conductor (PEC) boundary condition at $\rho=0$ physically models a long PEC thin wire along the $z$ axis in real 3D space. In this case, tangential components on the surface of the PEC rod are along z and $\phi$ directions while the $\rho$ component is perpendicular to the surface.
Since tangential electric field components vanish on a PEC surface, $E_z=E_\phi=0$. However, the normal component $E_\rho$ does not have to be zero when charged particles are present on the $z$-axis. This scenario is reversed when considering perfect magnetic conductor (PMC) boundary conditions.}
On the other hand, the radial boundary $\rho=0$ can be represented as a homogeneous Dirichlet boundary condition for the $E_{\phi}$ field. A homogeneous Neumann boundary condition for the electric field can be used to represent the perfect magnetic conductor case and a homogeneous Dirichlet boundary condition for the perfect electric conductor case.
More precisely, when edges in the E-B primal mesh (where $E_{\rho}$ and $E_z$ are defined) at $\rho=0$ are treated as unknowns, it corresponds to PMC or Neumann boundary conditions.
On the other hand, when the edges in the E-B primal mesh at $\rho=0$ are discarded, it corresponds to PEC or Dirichlet boundary conditions.
Similarly, when edges in the D-H primal mesh (where $H_\rho$ and $H_z$ are defined) at $\rho=0$ are treated as unknowns, it corresponds to PEC or Neumann boundary conditions.
Finally, when the edges in the D-H primal mesh at $\rho=0$ are discarded, it corresponds to PMC or Dirichlet boundary conditions.

The present BORPIC++ algorithm mainly considers linear superposition of eigenmodes with $m=0$.
Nevertheless, one can model the linear superposition of high-order azimuthal eigenmodes by introducing azimuthal variations in the field expansion, as explained in \cite{Na2019Finite}.
If $m\neq 0$, then $E_z=0$, $E_\rho\neq 0$, and $B_\phi \neq 0$, which can be implemented by assuming PEC boundary conditions on the axis, whereas $H_z=0$, $H_\rho\neq 0$, and $D_\phi\neq 0$, which can be implemented by assuming PMC boundary conditions.

Periodic boundary conditions (PBC) are used as lateral boundary conditions, i.e., at $z=L_{min}$ and $z=L_{max}$ where $L_{min}$ and $L_{max}$ are $z$ coordinates of the two ends of the problem domain. We employ a perfectly matched layers (PML) in the outward radial direction to model an open domain. The PML is considered later in this work  in a separate section.

\section{Gather Stage}
Discrete solutions of electromagnetic fields should be interpolated at the location of $p$-th charged particle $\mathbf{r}_{\parallel,p}^{n}$ (at $t=n\Delta t$) to evaluate the Lorentz force. Using the expansions in eqs. (\ref{eqn:exp1}) and (\ref{eqn:exp2}), the interpolated fields at the particle positions are obtained as
\begin{flalign}
{\boldsymbol{\mathcal{E}}}^{n}_{p}
&\equiv
\mathbf{E}'(\mathbf{r}_{\parallel,p}^{n},n\Delta t)
\nonumber \\
&=
\sum_{i=1}^{N_{1}}e_{\parallel,i}^{n}\mathbf{W}^{(1)}_{i}(\mathbf{r}_{\parallel,p}^{n})
+
\frac{1}{\epsilon_0}\sum_{i=1}^{N_{2}}d_{\perp,i}^{n}\mathbf{W}^{(2)}_{i}(\mathbf{r}_{\parallel,p}^{n}),
\\
{\boldsymbol{\mathcal{B}}}^{n}_{p}
&\equiv
\mathbf{B}'(\mathbf{r}_{\parallel,p}^{n},n\Delta t)
\nonumber \\
&=
\sum_{i=1}^{N_{2}}
\left(
\frac{b_{\perp,i}^{n-\frac{1}{2}}+b_{\perp,i}^{n+\frac{1}{2}}}{2}
\right)
\mathbf{W}^{(2)}_{i}(\mathbf{r}_{\parallel,p}^{n})
+
\mu_0\sum_{i=1}^{N_{1}}
\left(
\frac{h_{\parallel,i}^{n-\frac{1}{2}}+h_{\parallel,i}^{n+\frac{1}{2}}}{2}
\right)
\mathbf{W}^{(2)}_{i}(\mathbf{r}_{\parallel,p}^{n}).
\end{flalign}

\section{Pusher Stage}
Unlike a standard PIC algorithm,  BORPIC++ considers the kinematics of charge {\it rings}.
Each charge ring is represented as a computational point particle on $z \rho$-plane and the algorithm simulates of the motion of such point particles in the $z \rho$-plane.

We allow the charge rings to expand or compress radially, and to move forward and backward along the $z$ direction, while maintaining the same total charge\footnote{As a result, the {\it line charge density} on each ring varies according to the radial coordinate of the associated superparticle in the $z \rho$-plane.}. We also assume that the charge rings may rotate, producing an uniform current along the $\phi$ direction.
In the Pusher stage, the velocity and position of charged rings are updated in reaction to the Lorentz force obtained from the Gather stage.
In what follows we consider the non-relativistic case for the sake of simplicity\footnote{Implementation of a relativistic pusher in BORPIC++ does not pose any fundamental challenges~\cite{Na2018Relativistic}.}.

We express the position and velocity of the $p$-th superparticle in the $z \rho$-plane as
\begin{flalign}
\mathbf{r}_{\parallel,p}^{n} &= \hat{z} z_{p}^{n} + \hat{\rho} \rho_{p}^{n},\\
\mathbf{v}_{p}^{n+\frac{1}{2}}
&=
\underbrace{
\hat{z} {v_{p}}_{z}^{n+\frac{1}{2}}
+
\hat{\rho} {v_{p}}_{\rho}^{n+\frac{1}{2}}
}_{\mathbf{v}_{\parallel,p}^{n+\frac{1}{2}}}
+
\underbrace{
\hat{\phi} {v_{p}}_{\phi}^{n+\frac{1}{2}}
}_{\mathbf{v}_{\perp,p}^{n+\frac{1}{2}}}.
\end{flalign}

\subsection{Velocity update}
The velocity of $p$-th particle is updated using to the Lorentz force law as~\cite{Moon2015Exact}
\begin{flalign}
\left[
\begin{matrix}
v_{p,z}^{n+\frac{1}{2}}\\
v_{p,\rho}^{n+\frac{1}{2}}\\
v_{p,\phi}^{n+\frac{1}{2}}
\end{matrix}
\right]
=
\bar{\mathbf{N}}^{-1}
\cdot
\Biggl(
\bar{\mathbf{N}}^{T}
\cdot
\left[
\begin{matrix}
v_{p,z}^{n-\frac{1}{2}}\\
v_{p,\rho}^{n-\frac{1}{2}}\\
v_{p,\phi}^{n-\frac{1}{2}}
\end{matrix}
\right]
+
s
\left[
\begin{matrix}
\mathcal{E}_{p,z}^{n}\\
\mathcal{E}_{p,\rho}^{n}\\
\mathcal{E}_{p,\phi}^{n}
\end{matrix}
\right]
\Biggr)
\end{flalign}
where $s=\frac{Q_{p}\Delta t}{m_{p}}$, with $Q_{p}$ and $m_{p}$ being the charge and mass for $p$-th superparticle, respectively, and
\begin{flalign}
\bar{\mathbf{N}} &=
\frac{1}{2}
\left[
\begin{matrix}
2 & -s \mathcal{B}_{p,\phi}^{n} & s \mathcal{B}_{p,\rho}^{n} \\
s \mathcal{B}_{p,\phi}^{n} & 2 & -s \mathcal{B}_{p,z}^{n} \\
-s \mathcal{B}_{p,\rho}^{n} & s \mathcal{B}_{p,z}^{n} & 2  \\
\end{matrix}
\right].
\end{flalign}

\subsection{Particle position update}
After the velocity update, the position of $p$-th particle in the $(\rho,z)$ is updated as
\begin{flalign}
\mathbf{r}_{\parallel,p}^{n+1} = \mathbf{r}_{\parallel,p}^{n} + \Delta t \mathbf{v}_{\parallel,p}^{n+\frac{1}{2}}.
\end{flalign}
Note that we do not need to update the azimuth coordinate since we assume axisymmetry.

\section{Scatter Stage}

We next describe the scatter stage that maps the movement of the particles to a current distribution on the mesh. We consider the TE\textsuperscript{$\phi$} and TM\textsuperscript{$\phi$} cases separately in sequence. Importantly, the scatter stage described below ensures charge-conservation on an unstructured mesh~\cite{Moon2015Exact}, i.e. it ensures satisfaction of (the discrete analogue of) the charge continuity equation
$d\varrho/dt + \nabla \cdot \varrho \mathbf{v} = 0$ on the mesh, where $\varrho$ is the volumetric charge density. A similar strategy has been
considered in~\cite{pinto2014charge,9599150}.

\subsection{TE\textsuperscript{$\phi$} case}\label{Sec:Parallel_Scatter}
Assuming that the charged point particles are represented by delta functions on the $z \rho$-plane, the current density $\mathbf{J}'_{\parallel}$ produced by the collective motion of the superparticles in the plasma medium is given by\footnote{As noted before, a point charge in the $z \rho$ plane represents a charge ring in three-dimensional space. $Q_p$ is the total charge in such ring. The $\rho$- and $z$-components of the current density $\mathbf{J}_{\parallel}$ that arise from the movement of the associated superparticle in the $z \rho$ plane are {\it line} current densities (along the ring) in the three-dimensional space. These densities scale with $1/\rho$, where $\rho$ is the radius of the ring (or, equivalently, the radial coordinate of the superparticle). In eq.~(\ref{eqn:Jcollective}), such $1/\rho$ factor gets canceled out by the factor $\rho$ that arises in the mapping of the $\rho$ and $z$ components from $\mathbf{J}$  to $\mathbf{J}'$  in eq.~(\ref{eqn:transfJ}).}

\begin{align}
\label{eqn:Jcollective}
\mathbf{J}'_{\parallel}(\mathbf{r}_{\parallel},t)=\sum_{p=1}^{N_{p}}
\mathbf{v}_{p}(t) \, Q_{p} \, \delta\left(\mathbf{r}_{\parallel}-\mathbf{r}_{\parallel,p}(t)\right),
\end{align}
where $N_{p}$ denotes the total number of superparticles and we assume the shape factor of each particle in the $z \rho$-plane to be a Dirac delta function, i.e., a point particle.

From eq.~(\ref{eqn:Jproj}), the current scattered on edge $i$ due to the motion of the charged particles in the $z \rho$-plane is expressed as the Galerkin projection of the average current during a time interval $\Delta t$ onto each edge $i$ through an inner product with the associated edge element (Whitney basis) function $\mathbf{W}^{(1)}_{i}$, i.e.,
\begin{flalign}
\label{eq:proj1}
j^{n+\frac{1}{2}}_{\parallel,i}
=\frac{1}{\Delta t}
\int_{n\Delta t}^{(n+1)\Delta t} \left[
\iint_{\cup_i \mathcal{T}}
\mathbf{J}'_{\parallel}(\mathbf{r}_{\parallel},t) \cdot
\mathbf{W}^{(1)}_{i}(\mathbf{r}_{\parallel}) \, dA \right] \, dt
\end{flalign}
where $\cup_i \mathcal{T}$ denotes the union of the two area elements that touch edge $i$. Substituting eq.~(\ref{eqn:Jcollective}) in eq.~(\ref{eq:proj1}),
and using the sifting property of the delta function:
\begin{flalign}
\iint_{\cup_i \mathcal{T}}\delta\left(\mathbf{r}_{\parallel}-\mathbf{r}_{\parallel,p}(t)\right)
\mathbf{W}^{(1)}_{i}(\mathbf{r}_{\parallel}) \, dA = \mathbf{W}^{(1)}_{i}  (\mathbf{r}_{\parallel,p}(t) ),
\end{flalign}
allows us to write the following expression for the current scattered on edge $i$ due to the motion of the charged particles in the $z \rho$-plane~\cite{Teixeira2013Differential,Moon2015Exact,Na2016Local,Na2017Axisymmetric}:
\begin{flalign}
j^{n+\frac{1}{2}}_{\parallel,i}
=
\sum_{p=1}^{N_{p}}
\frac{Q_{p}}{\Delta t}
\int_{\mathbf{r}^{n}_{\parallel,p}}^{\mathbf{r}^{n+1}_{\parallel,p}}
\mathbf{W}^{(1)}_{i}(\mathbf{r}_{\parallel})\cdot d\mathbf{l},
\end{flalign}
where we used  fact that,
assuming a constant $\mathbf{v}(t)$ in the interval $\Delta t$, $\mathbf{v}(t) dt = d \mathbf{l}$ where $d \mathbf{l}$ is an infinitesimal displacement.

\subsection{TM\textsuperscript{$\phi$} case}\label{Sec:Perp_Scatter}

The azimuthal current density $\mathbf{J}'_{\perp}$ due to rotation of the charge rings is expressed as
\begin{flalign}
\mathbf{J}'_{\perp}(\mathbf{r}_{\parallel},t)
=
\sum_{p=1}^{N_{p}}
\hat{\phi}
\left(\hat{\phi}\cdot\mathbf{v}_{p}(t)\right) \, Q_{p} \,
S\left(\mathbf{r}_{\parallel}-\mathbf{r}_{\parallel,p}(t)\right),
\end{flalign}
where $S(\cdot)$ is a particle shape function to be discussed below and $N_{p}$ is the total number of superparticles.
Using again a Galerkin projection but now based on $\mathbf{W}^{(2)}_{i}$,
the associated scatter current is given by
\begin{flalign}
\label{eq:proj2}
j^{n+\frac{1}{2}}_{\perp,i}
=
\iint_{\mathcal{T}_{i}}
\mathbf{J}'_{\perp}\Bigl(\mathbf{r}_{\parallel},(n+\frac{1}{2})\Delta t\Bigr) \cdot
\mathbf{W}^{(2)}_{i}(\mathbf{r}_{\parallel})
~
dA
\end{flalign}
on each $i$-th cell of the primal mesh.
To numerically approximate this integral, the particle position can be evaluated as $\mathbf{r}_{\parallel,p}^{n+\frac{1}{2}}=(\mathbf{r}_{\parallel,p}^{n+1}+\mathbf{r}_{\parallel,p}^{n})/2$.
It is important to note that the total charge of each ring (represented as a point charge in the $z\rho$-plane) remains constant as the ring expands or contracts over time.
In addition, since it is assumed that each charged ring moves at a constant speed during $\Delta t$, it always forms a uniform current.
Hence, azimuthal currents are charge-conserving.

Among the various possible shape functions, here we employ the following:
\begin{flalign}
S(\mathbf{r}-\mathbf{r}_{p})
=
P^{(m)}_{z}(z-z_p)
P^{(m)}_{\rho}(\rho-\rho_p)
\end{flalign}
with symmetric polynomial factors of the form
\begin{flalign}
P^{(m)}_{\zeta}(\zeta)
=
\begin{cases}
\left[1-(\zeta / H_{m})^2\right]^{m}/\alpha,
& \left|\zeta / H_{m}\right| \leq 1 \\
0,
& \text{elsewhere}
\end{cases}
\end{flalign}
where $H_{m}=\alpha h_{m}$ is determined so that the integral of $P^{(m)}_{\zeta}(\zeta)$ is unity and $\alpha$ is an arbitrary scale factor to control the size of the superparticle.
In particular, $h_{0}=0.5$, $h_{1}=0.75$, $h_{2}=0.9375$, and $h_{3}=1.039475$.

{\it Remark 3}: The choice of shape functions (finite-sized for ${\mathbf{J}'}_{\perp}$ and delta for ${\mathbf{J}'}_{\parallel}$) is dictated by the degree of continuity in the remaining functions present in the integrands of~(\ref{eq:proj1}) and~(\ref{eq:proj2}) and by the need to avoid discontinuities of the scatter current as the particles cross from one mesh element to another. In the $z \rho$-plane, the vector function $\mathbf{W}^{(1)}_{i}$ is tangentially continuous across mesh elements while $\mathbf{W}^{(2)}_{i}$ is discontinuous across elements~\cite{10.1145/1141911.1141991}. The use of a delta function in~(\ref{eq:proj2}) together with (discontinuous) $\mathbf{W}^{(2)}_{i}$ would cause excessive numerical noise because of the sudden current turn-on or -off by particles crossing from one element to another\footnote{This can be interpreted as a spurious numerical radiation originating from the interaction between grid dispersion and space charge dispersion.
The shape function can interpreted as a low pass filter imposed on a point particle so that spurious radiation is reduced.}.
On the other hand, as the spatial size (and the polynomial order) of the shape function increases, the computational load associated with the numerical evaluation of the scatter current integrals
increases. Hence, a fundamental trade-off exists that is managed by the above shape function choices.

\begin{figure}
\centering
\includegraphics[width=.5\linewidth]
{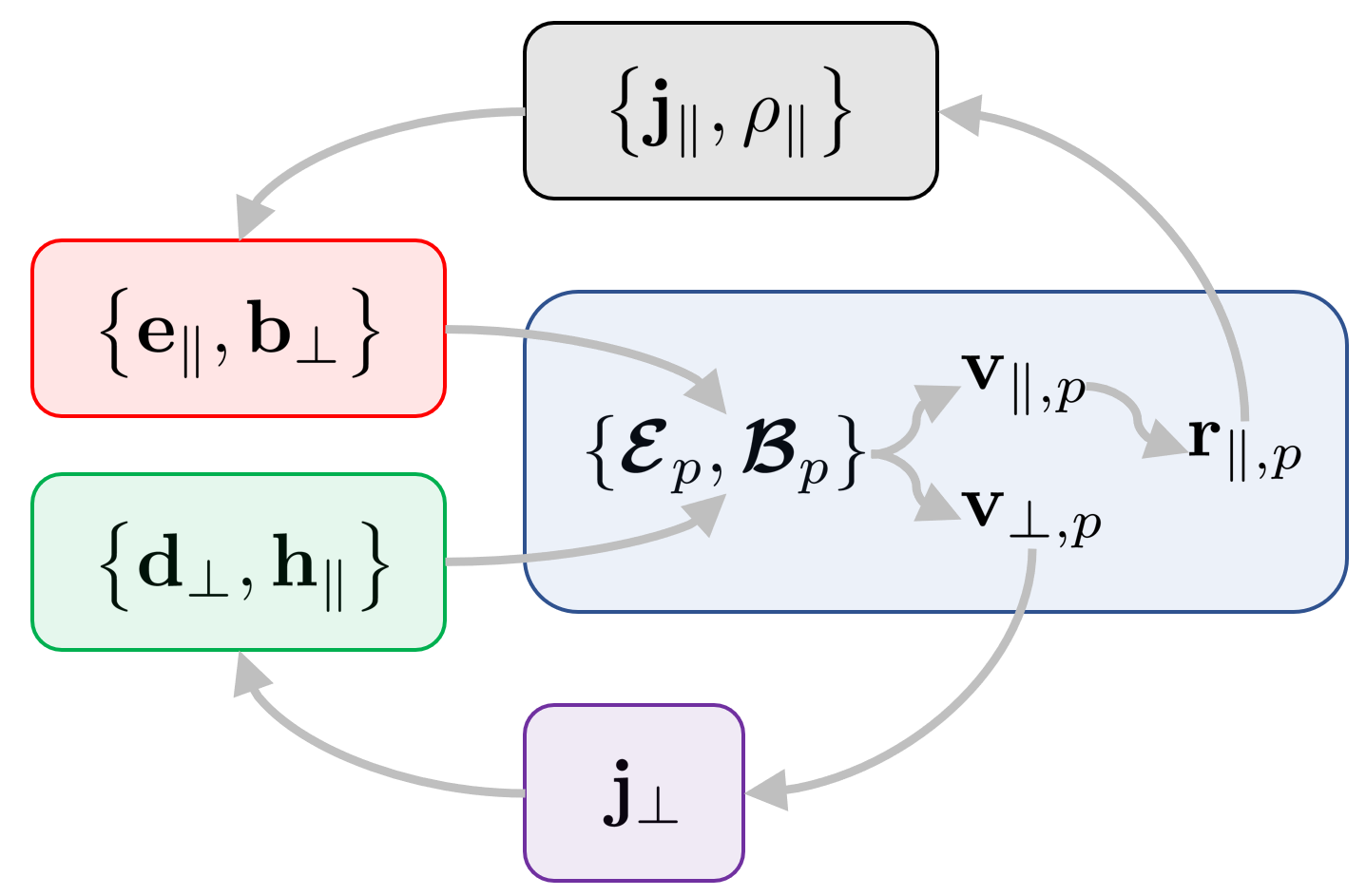}
\caption{The overall update schematic of the proposed BORPIC algorithm.}
\label{fig:overall_update_scheme}
\end{figure}

Fig. \ref{fig:overall_update_scheme} summarizes the sequential update of the dynamic variables during each time step of BORPIC++.

\section{Radial PML implementation}
In order to study the generation of lower hybrid and whistler waves and other radiation by ion velocity rings and diverse types of beams in an open space, it is necessary to incorporate an absorbing boundary condition (ABC) at the boundary of the computational domain. Otherwise, radiated waves would reflect back to the computational domain from the mesh termination, leading to erroneous results. BORPIC++ implements a perfectly matched layer (PML) absorbing boundary condition along the radial direction.
The radial PML implemented in BORPIC++ absorbs outgoing electromagnetic waves very effectively and enables accurate simulations of long-term plasma behavior.

One can refer to \ref{app:radialPML} for a comprehensive description of the radial PML implementation into the TE\textsuperscript{$\phi$} and TM\textsuperscript{$\phi$} field solvers.

\section{Numerical results}
To verify the accuracy of our algorithm implementation, we have conducted a number of numerical tests discussed in detail below.

\subsection{Impressed surface currents}
\begin{figure}
\centering
\subfloat[TE\textsuperscript{$\phi$} configuration.]
{\includegraphics[width=.6\linewidth]{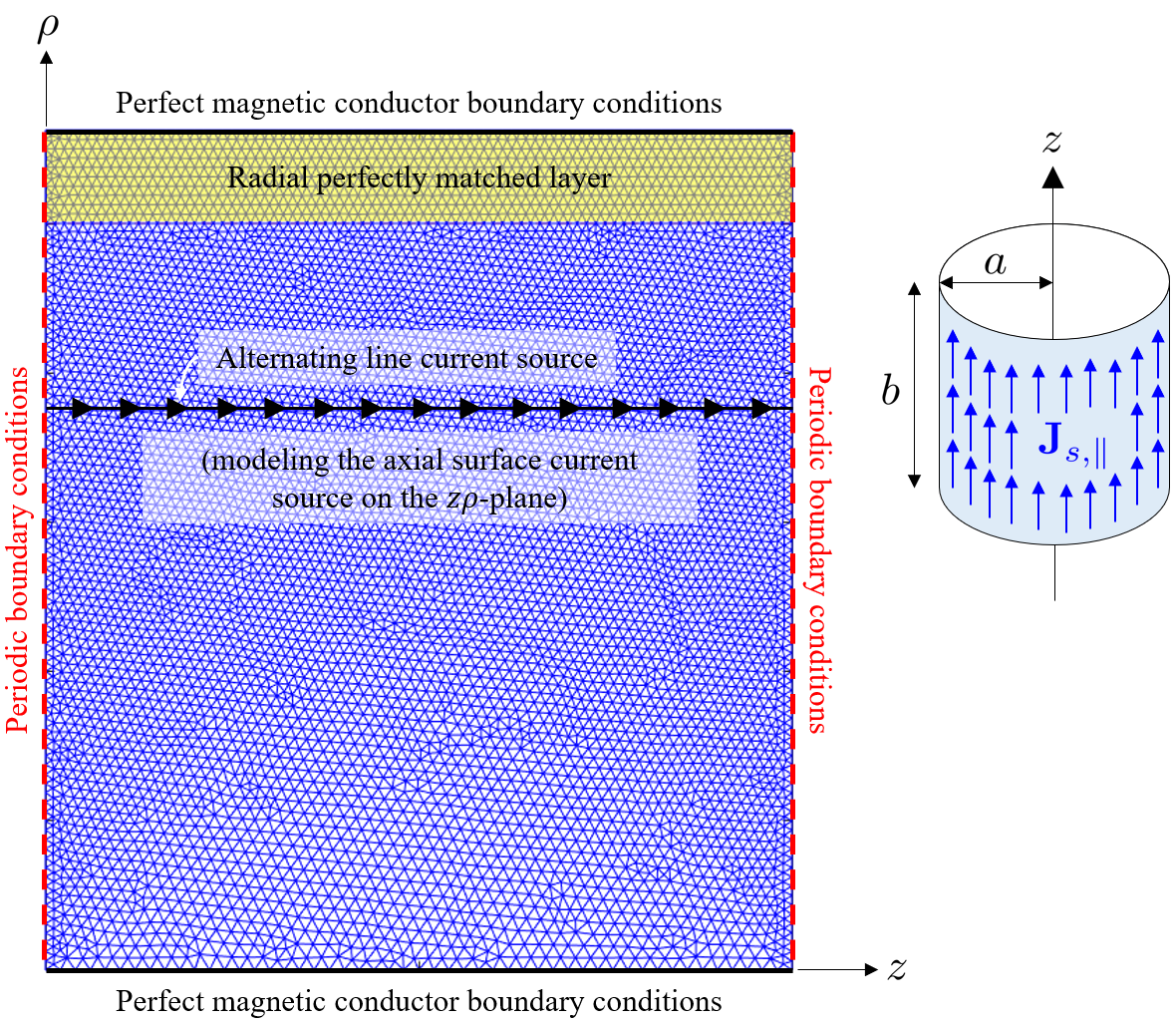}}
\quad\quad
\\
\subfloat[TM\textsuperscript{$\phi$} configuration.]
{\includegraphics[width=.6\linewidth]{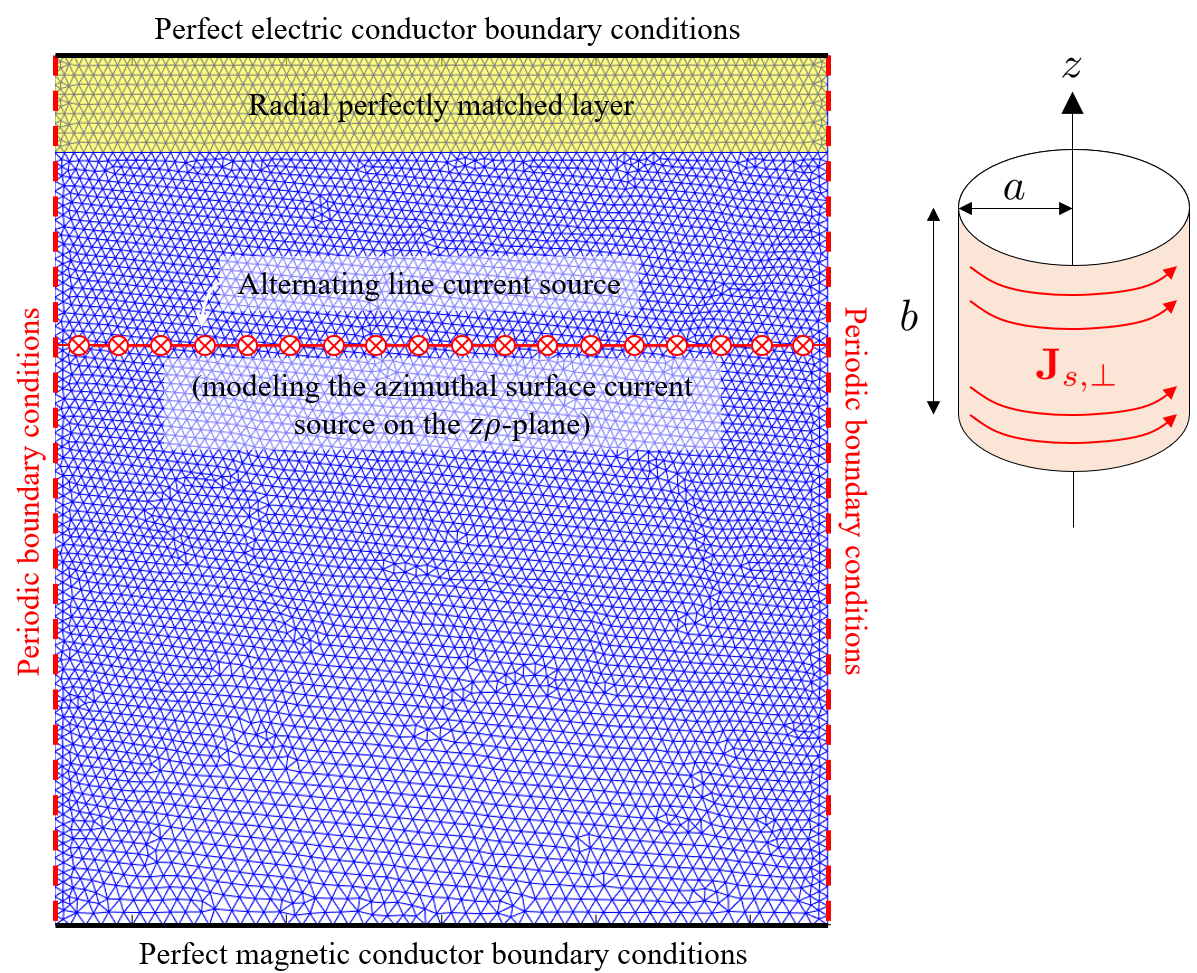}}
\caption{
Problem description, geometry, and finite-element mesh for analyzing the radial PML performance.
}
\label{fig:problem_geometry_mesh}
\end{figure}

\begin{figure}
\centering
\includegraphics[width=.4\linewidth]
{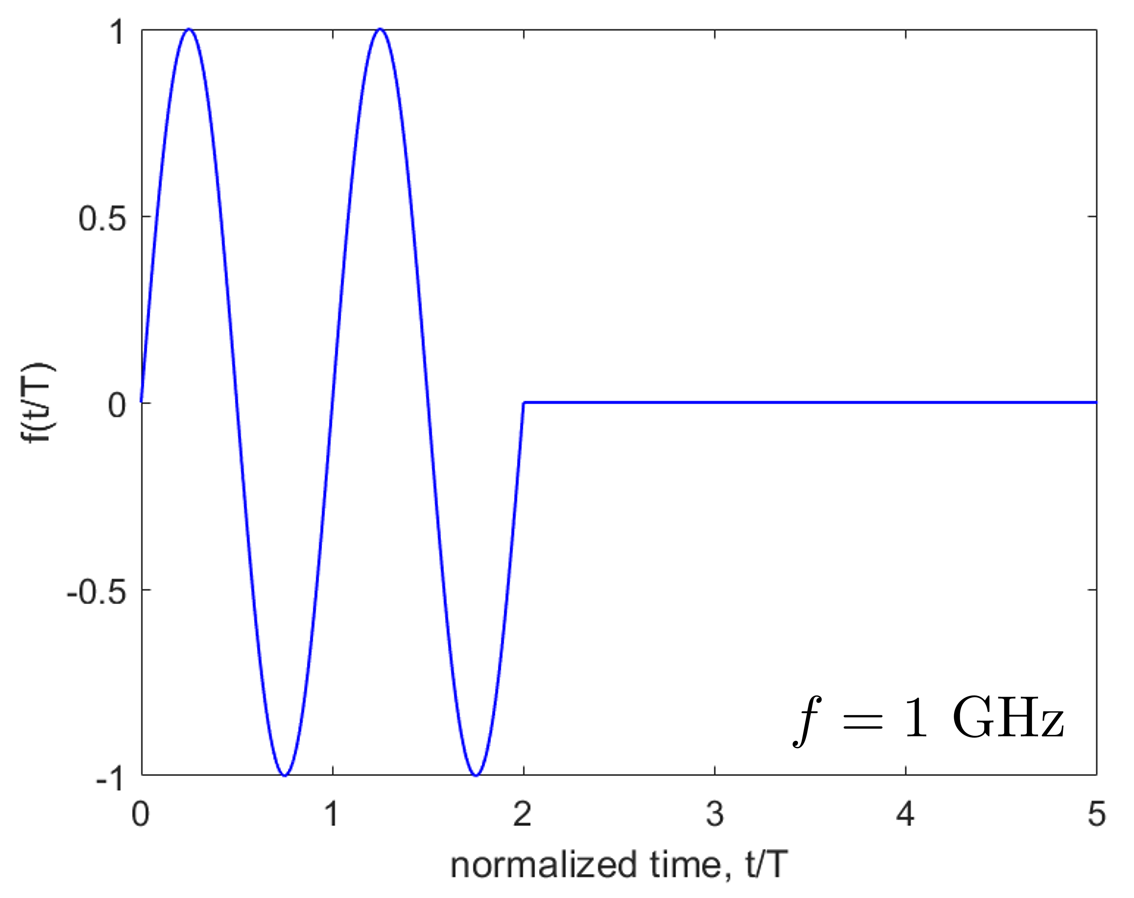}
\caption{Temporal pulse profile of the impressed surface current densities.}
\label{fig:pulse}
\end{figure}

To analyze the performance of the radial PML, we first perform simulations with electromagnetic fields only (no particles).
Specifically, we observe the radiation fields due to two types of current densities impressed on a cylinder surface of radius $a$ and height $b$, as depicted in Fig.~\ref{fig:problem_geometry_mesh}.
One current (denoted by $\mathbf{J}_{s,\parallel}$) flows along the $z$ direction and the other current (denoted by $\mathbf{J}_{s,\perp}$) flows along the $\phi$ direction.
Thus, $\mathbf{J}_{s,\parallel}$ produces a TE\textsuperscript{$\phi$}  whereas $\mathbf{J}_{s,\perp}$ produces a TM\textsuperscript{$\phi$} field.
The two surface currents are excited with the same temporal pulse profile, a truncated sinusoidal signal with $f=1$ GHz , as illustrated in Fig.~\ref{fig:pulse}.
We use an unstructured mesh consisting of 7,041 nodes, 20,810 edges, and 13,770 triangular elements.
The physical domain (excluding the PML region) consists of the region $-0.5\text{~[m]~}\leq z \leq 0.5$ [m] and $0 \leq \rho \leq 1$ [m].

\begin{figure}
\centering
\includegraphics[width=\linewidth]
{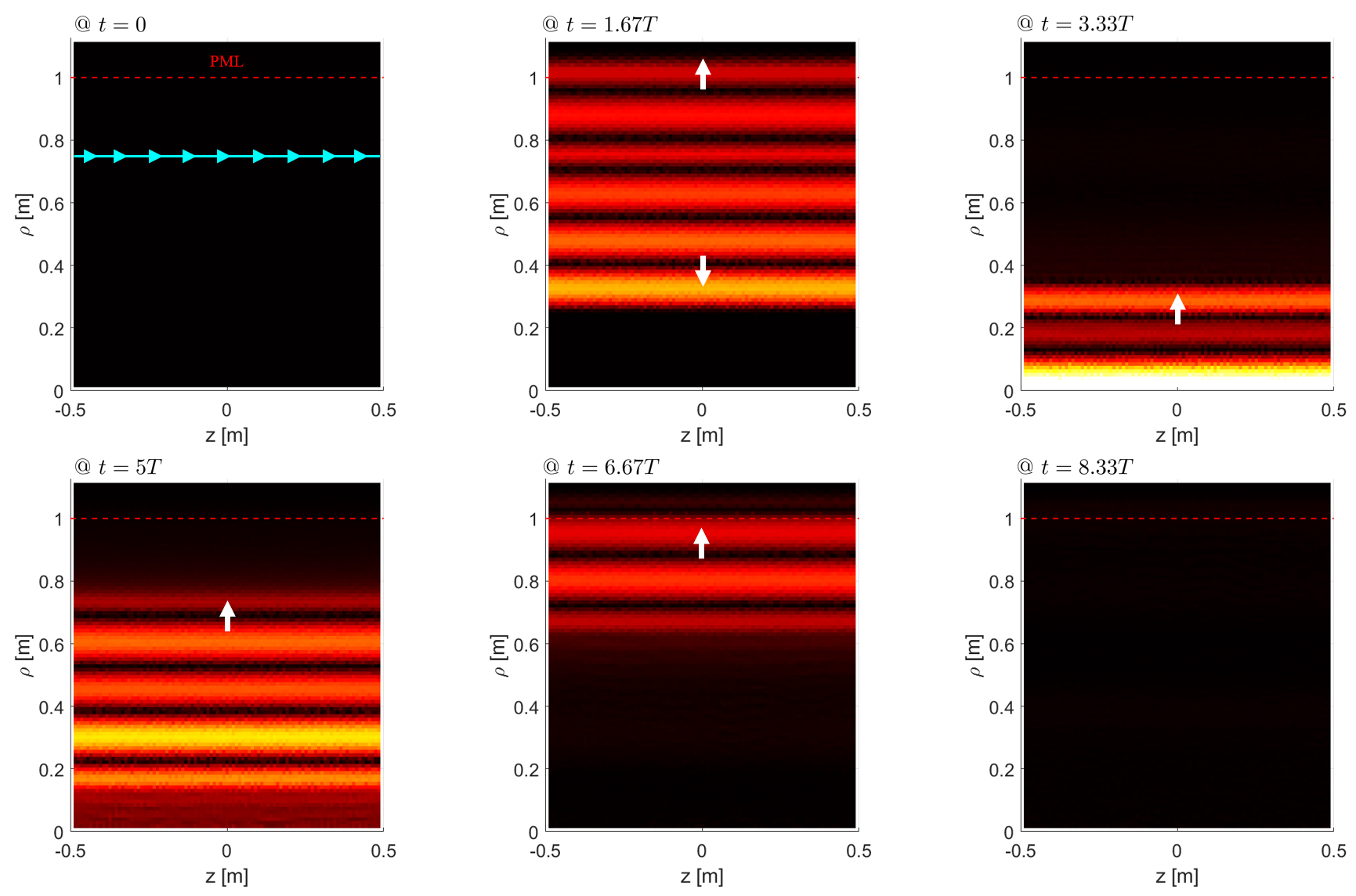}
\caption{$\left|\mathbf{E}_{\parallel}(\mathbf{r}_{\parallel},t)\right|$ due to $\mathbf{K}_{\parallel}$ at different time steps.}
\label{fig:TE_Phi_E_Field}
\end{figure}

\begin{figure}
\centering
\includegraphics[width=\linewidth]
{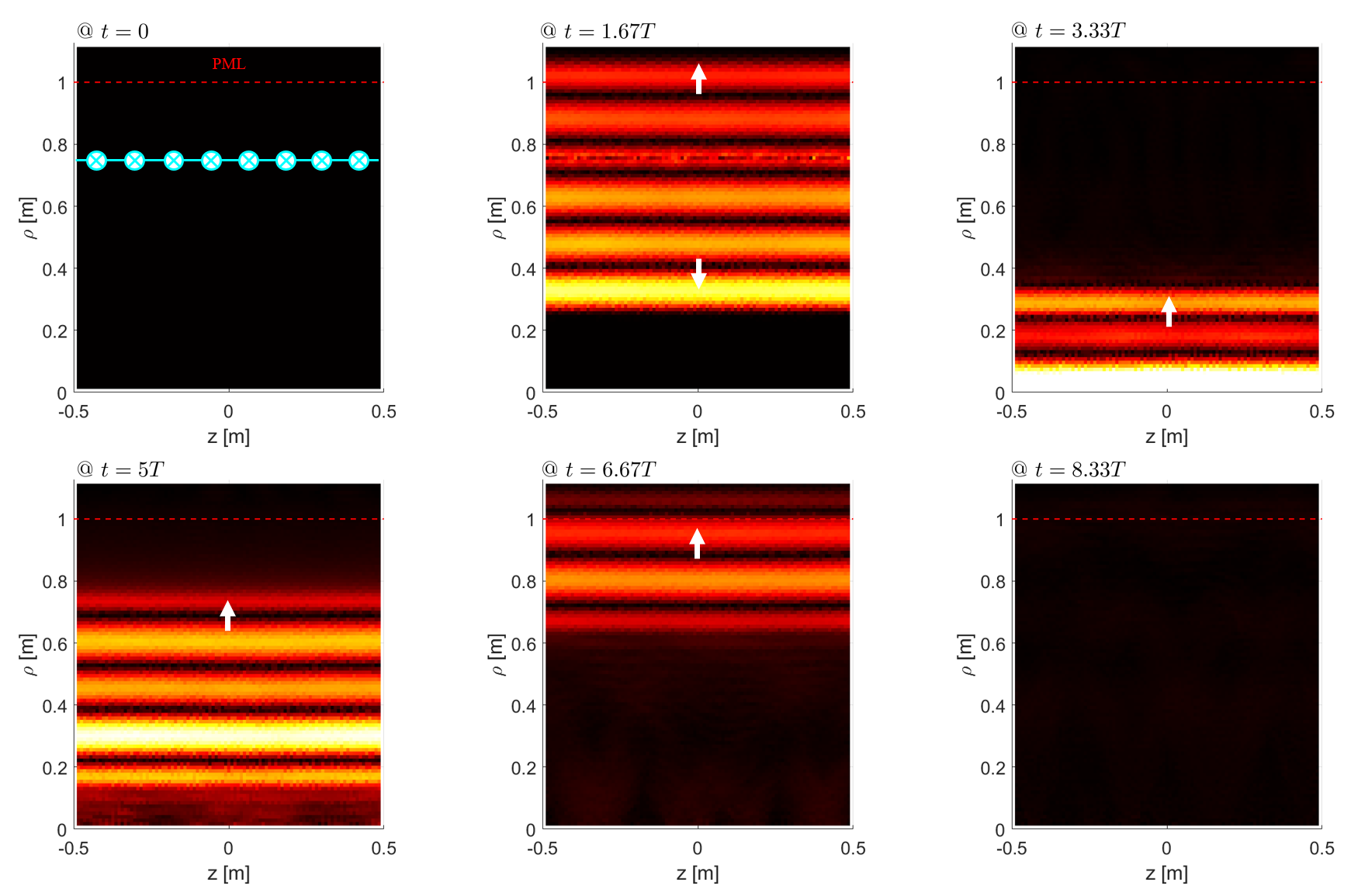}
\caption{$\left|\mathbf{H}_{\parallel}(\mathbf{r}_{\parallel},t)\right|$ due to $\mathbf{K}_{\perp}$ at different time steps.}
\label{fig:TM_Phi_H_Field}
\end{figure}

Fig.~\ref{fig:TE_Phi_E_Field} shows snapshots of the $\left|\mathbf{E}_{\parallel}(\mathbf{r}_{\parallel},t)\right|$ due to $\mathbf{J}_{s,\parallel}$ at different time steps.
Similarly, Fig. \ref{fig:TM_Phi_H_Field} shows snapshots of $\left|\mathbf{H}_{\parallel}(\mathbf{r}_{\parallel},t)\right|$ due to $\mathbf{J}_{s,\perp}$ at different time steps.
In both simulations, it is observed that the outgoing waves are absorbed very effectively as they enter in the PML region.

\subsection{Charge ring motion}
We next consider the motion of a charge ring, which is a combination of its poloidal ($\rho$ and $z$) and azimuthal ($\phi$) motions.
We also analyze the resulting electromagnetic fields.

A charged superparticle ($10^6$ electrons) is initialized at $(\rho,z)=(0.75,-0.15)$ [m] and launched with initial velocity $\mathbf{v}_{p}=(0.0129c,0.0233c,0.0025c)$ where $c$ denotes the speed of light in the free space. The time increment $\Delta t$ is set to be $\Delta l/c$ [s] where $\Delta l=0.001$ [m].
First we consider the poloidal $\mathbf{J}_{\parallel}$ due to the axial and radial motions of the charge ring only.
The PML extends away from its front end surface positioned at $\rho = 1.0$ [m].
We tested particle reflector surfaces placed at three different radial locations: inside the PML at $\rho=1.1$ [m] (Test Case 1), coinciding with the PML front end at $\rho=1.0$ [m] (Test Case 2), and at some distance off from the PML at $\rho=0.9$ [m] (Test Case 3).
These choices are made to examine the interaction between the near field produced by the superparticle and the PML.
Note that, ideally, the PML should absorb any electromagnetic radiation out of the computational domain and, at the same time, avoid any spurious interactions with the near-field (nonradiative) produced by the particles.
We used an irregular mesh consisting of $7,041$ nodes, $20,810$ edges, and $13,770$ cells.

\begin{figure}
\centering
\includegraphics[width=\linewidth]
{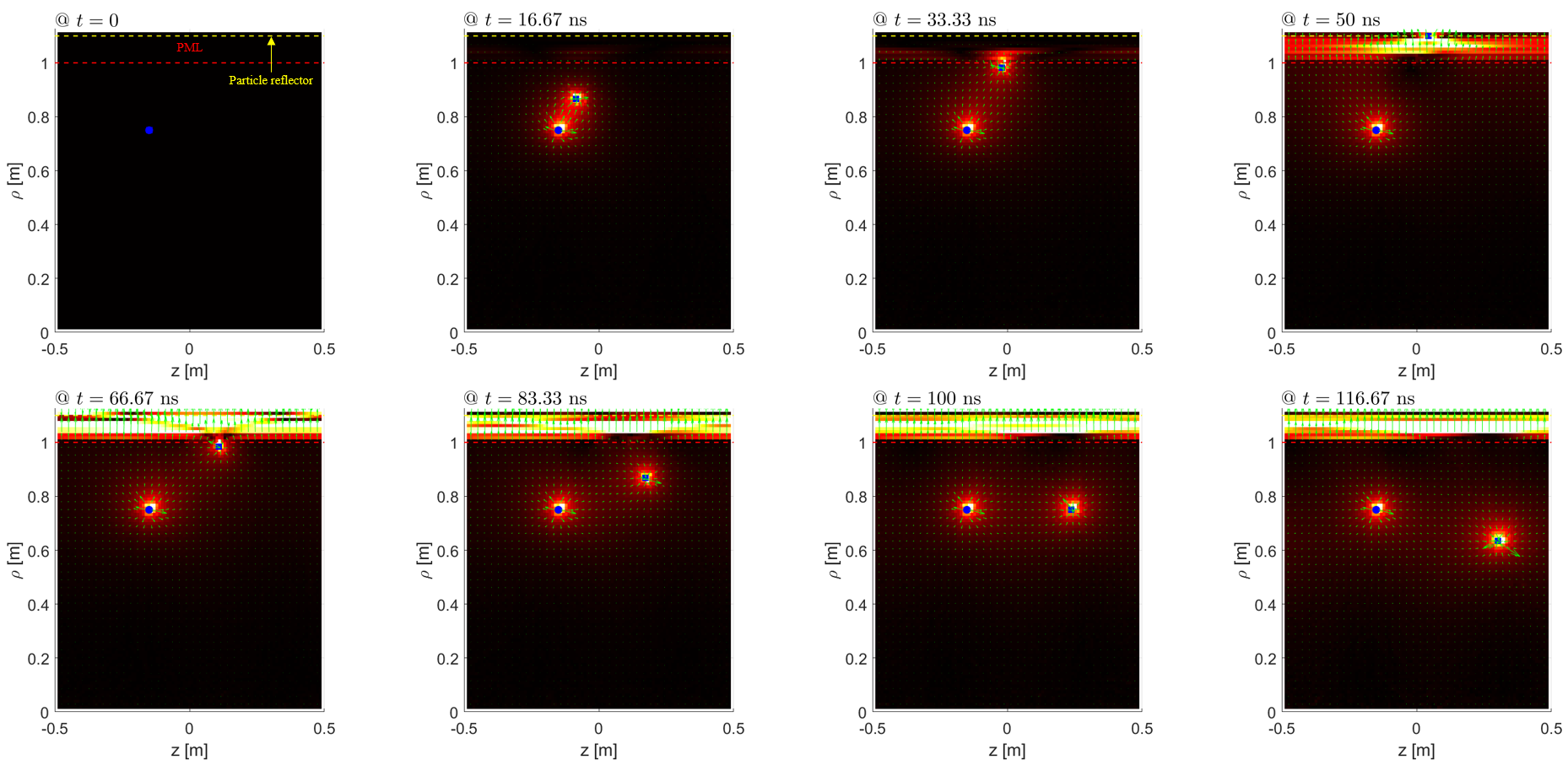}
\caption{$\left|\mathbf{E}_{\parallel}(\mathbf{r}_{\parallel},t)\right|$ at different time steps for Test Case 1.}
\label{fig:PIC_TE_Phi_E_Case_1}
\end{figure}

\begin{figure}
\centering
\includegraphics[width=\linewidth]
{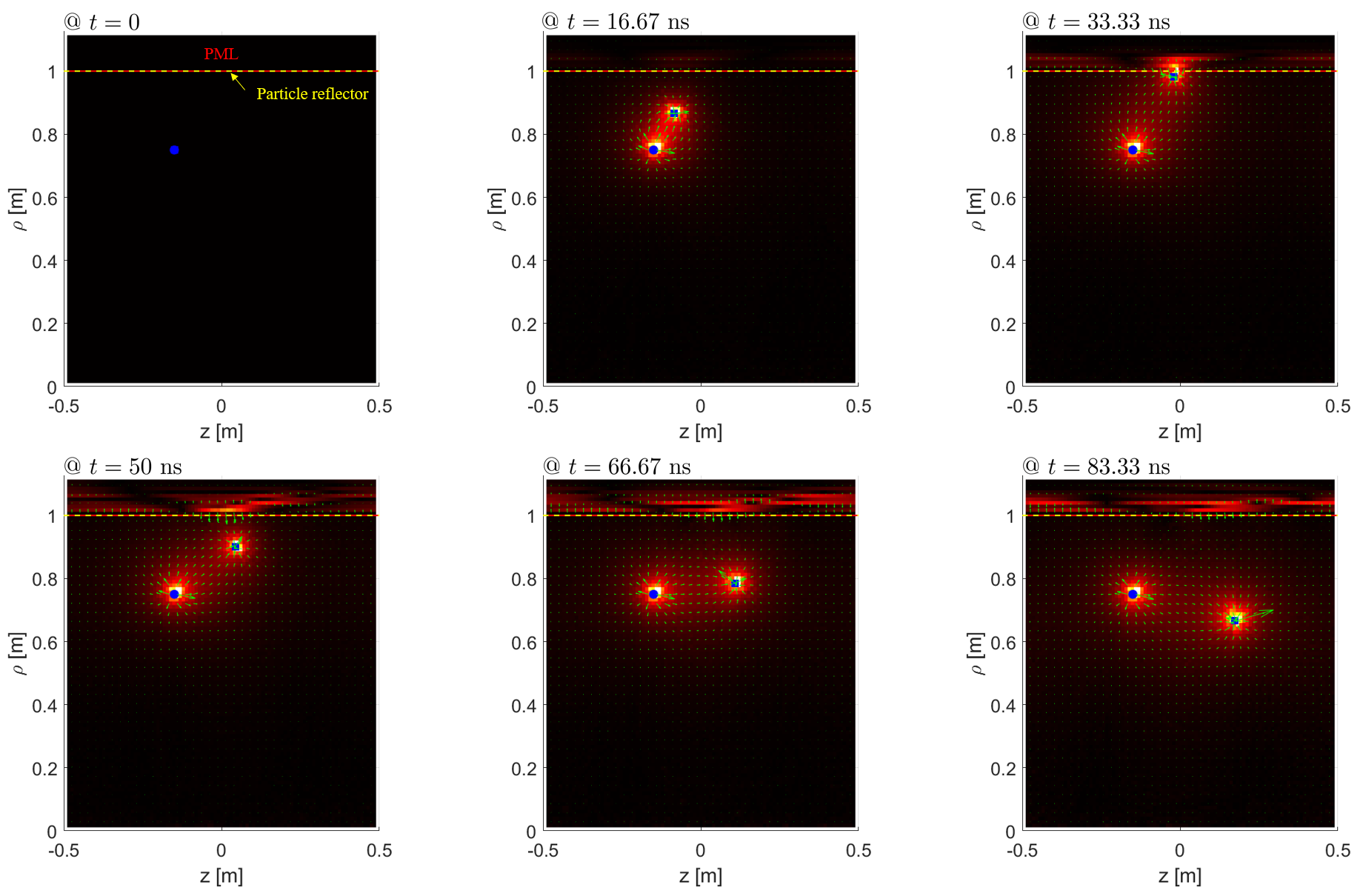}
\caption{$\left|\mathbf{E}_{\parallel}(\mathbf{r}_{\parallel},t)\right|$ at different time steps for Test Case 2.}
\label{fig:PIC_TE_Phi_E_Case_2}
\end{figure}

\begin{figure}
\centering
\includegraphics[width=\linewidth]
{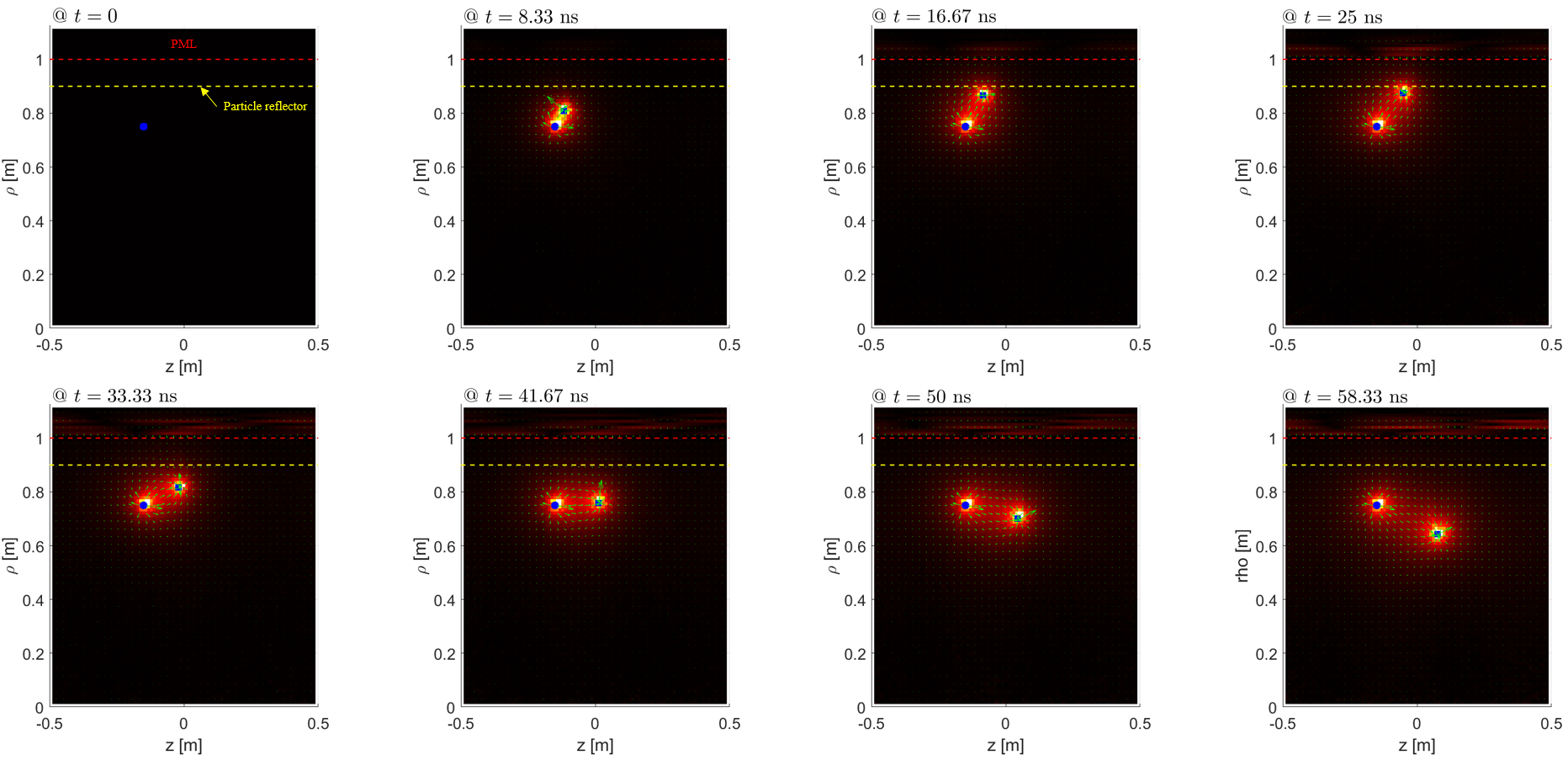}
\caption{$\left|\mathbf{E}_{\parallel}(\mathbf{r}_{\parallel},t)\right|$ at different time steps for Test Case 3.}
\label{fig:PIC_TE_Phi_E_Case_3}
\end{figure}

Figs. \ref{fig:PIC_TE_Phi_E_Case_1},  \ref{fig:PIC_TE_Phi_E_Case_2}, and \ref{fig:PIC_TE_Phi_E_Case_3} show $\left|\mathbf{E}_{\parallel}(\mathbf{r}_{\parallel},t)\right|$ at different time steps for Test Cases 1, 2, and 3, respectively.
These figures indicate that, as expected, the placement of a particle reflector surface within the PML or coinciding with the front-end of the PML is not desirable because the PML indeed perturbs the near field of the particle. In addition, a spurious field is clearly induced inside the PML region (from the interaction with the nonradiative near-field of the particle. On the other hand, Fig.~(\ref{fig:PIC_TE_Phi_E_Case_3}) shows that, once the reflecting surface is placed at a sufficient distance away from the PML, only a very small residual field is present inside the PML region and the field perturbation in the physical region is negligible.

\begin{figure}
\centering
\includegraphics[width=\linewidth]
{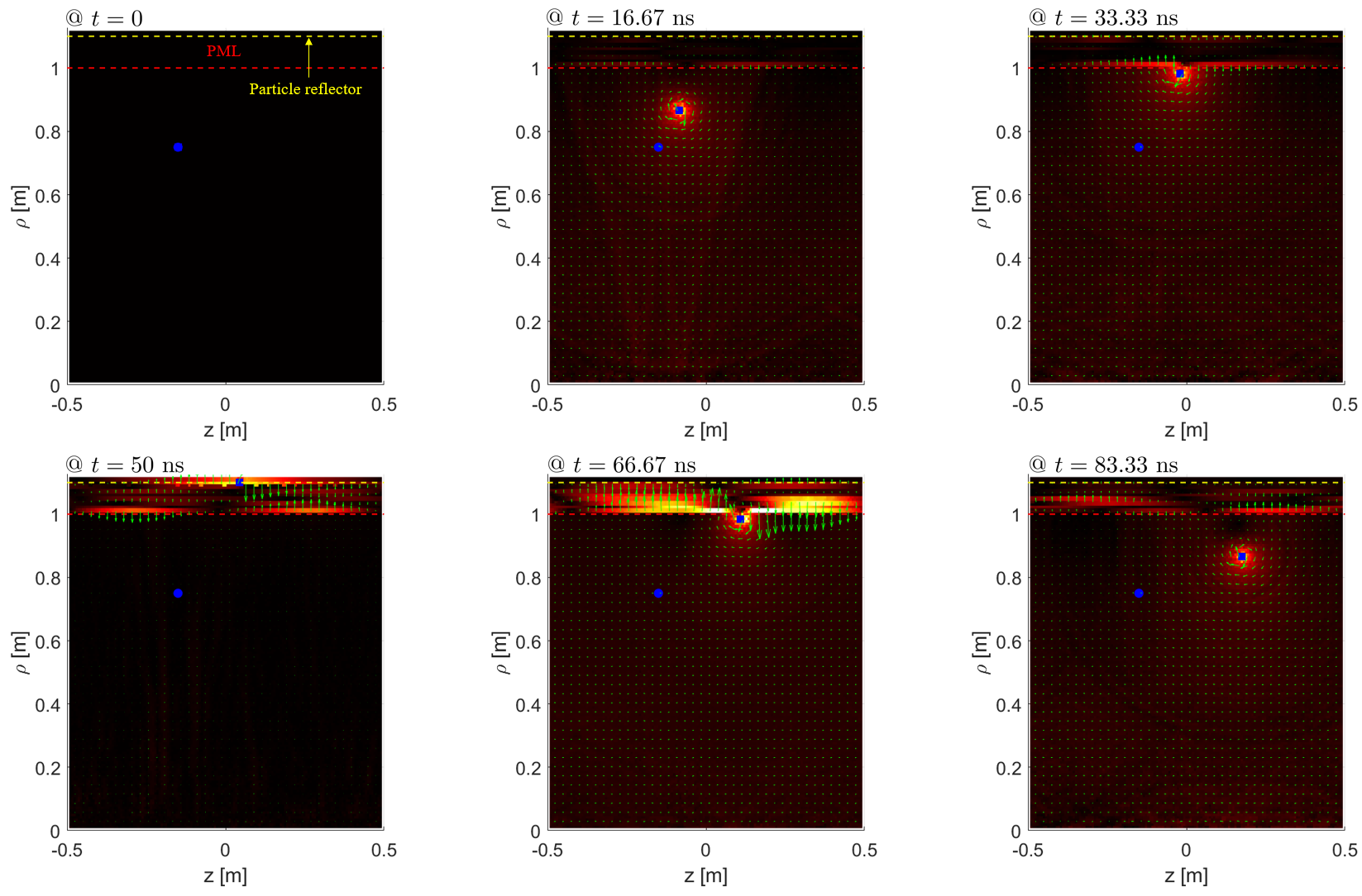}
\caption{$\left|\mathbf{H}_{\parallel}(\mathbf{r}_{\parallel},t)\right|$ at different time steps for Case 1.}
\label{fig:PIC_TM_Phi_H_Case_1}
\end{figure}

\begin{figure}
\centering
\includegraphics[width=\linewidth]
{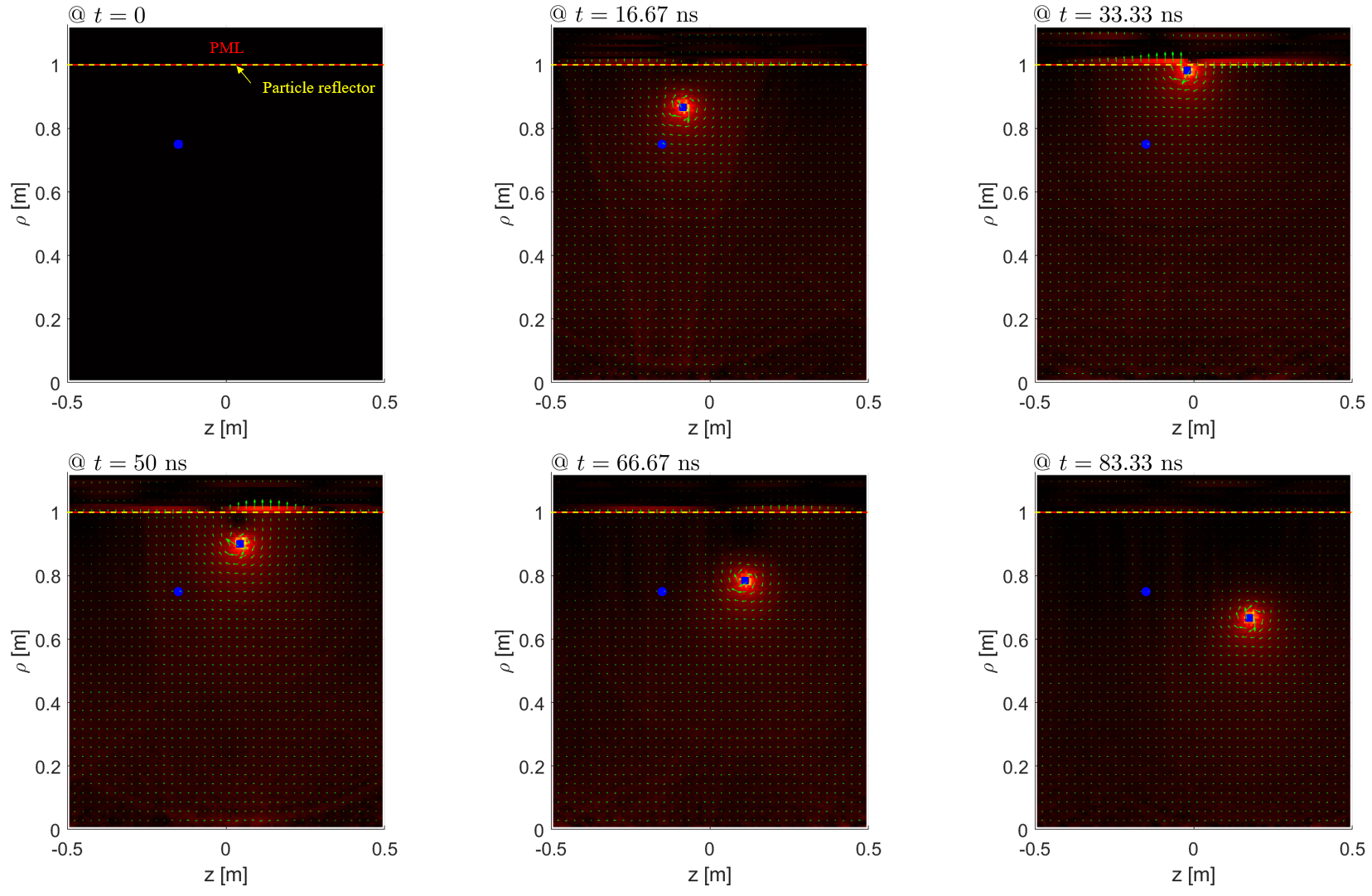}
\caption{$\left|\mathbf{H}_{\parallel}(\mathbf{r}_{\parallel},t)\right|$ at different time steps for Test Case 2.}
\label{fig:PIC_TM_Phi_H_Case_2}
\end{figure}

\begin{figure}
\centering
\includegraphics[width=\linewidth]
{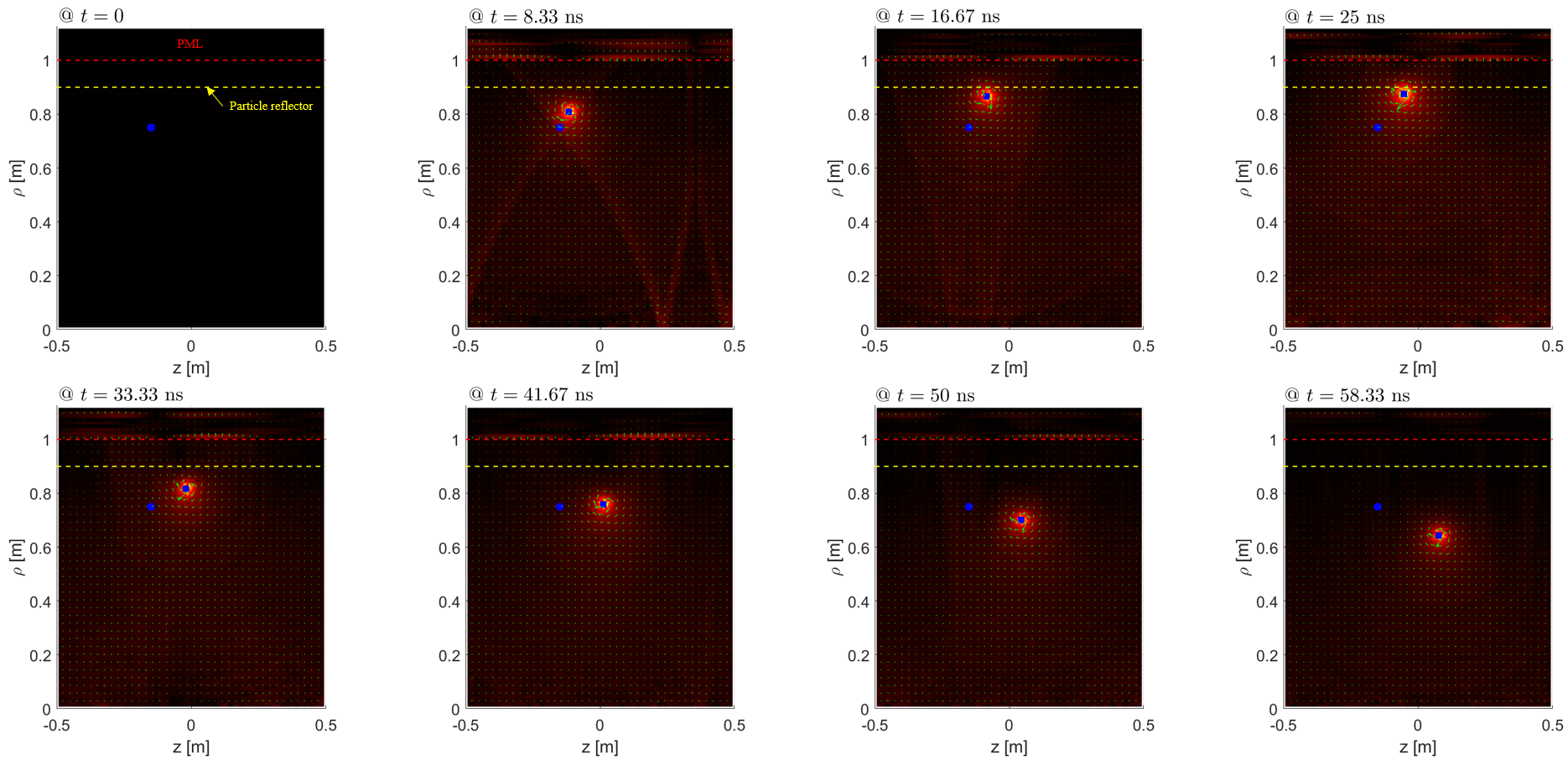}
\caption{$\left|\mathbf{H}_{\parallel}(\mathbf{r}_{\parallel},t)\right|$ at different time steps for Test Case 3.}
\label{fig:PIC_TM_Phi_H_Case_3}
\end{figure}

The next example considers a charge ring with current $\mathbf{J}_{\perp}$ due to azimuthal rotation, with other parameters staying the same. Figs.~\ref{fig:PIC_TM_Phi_H_Case_1}, \ref{fig:PIC_TM_Phi_H_Case_2}, and \ref{fig:PIC_TM_Phi_H_Case_3} show snapshots of $\left|\mathbf{H}_{\parallel}(\mathbf{r}_{\parallel},t)\right|$ at different time steps for Test Cases 1, 2, and 3, respectively. Again, it can be observed that, in order to minimize spurious fields inside the PML and particle perturbations, a buffer region should be placed between the particle reflecting surface and the front-end of the PML.  In Fig.~\ref{fig:PIC_TM_Phi_H_Case_3}, a radiation wavefront from the sudden initialization of the superparticle at $t=0$ is clearly visible, particularly at the early time $t=8.33$ ns (such type of wavefront is also present  in Figs.~\ref{fig:PIC_TM_Phi_H_Case_1} and \ref{fig:PIC_TM_Phi_H_Case_2},
 but it is less visible there because of those sequences of plots refer to later times). Since we assume zero field  initial conditions, the initialization and movement of the superparticle induces a stationary charge (infinitely massive ion) of opposite value at the initial location of the superparticle.  We thus have the sudden appearance at $t=0$ of a charge dipole, which produces the said radiation. This radiation is reflected from $\rho=0$ and the periodic boundary conditions along $z$, but it is well absorbed by the PML.

\subsection{Gyromotion of a charge particle in $z \rho$-plane}
\begin{figure}
\centering
\subfloat[With radial PML.]
{\includegraphics[width=.45\linewidth]{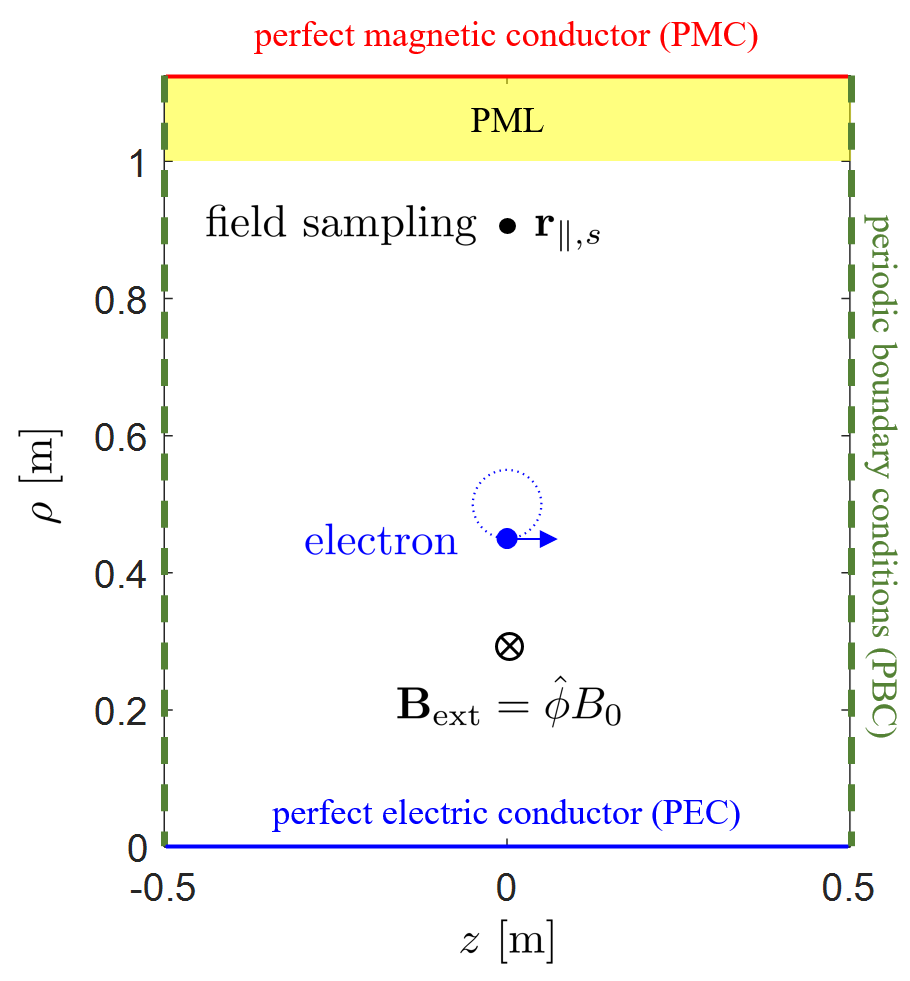}}
\quad\quad
\subfloat[Without radial PML.]
{\includegraphics[width=.45\linewidth]{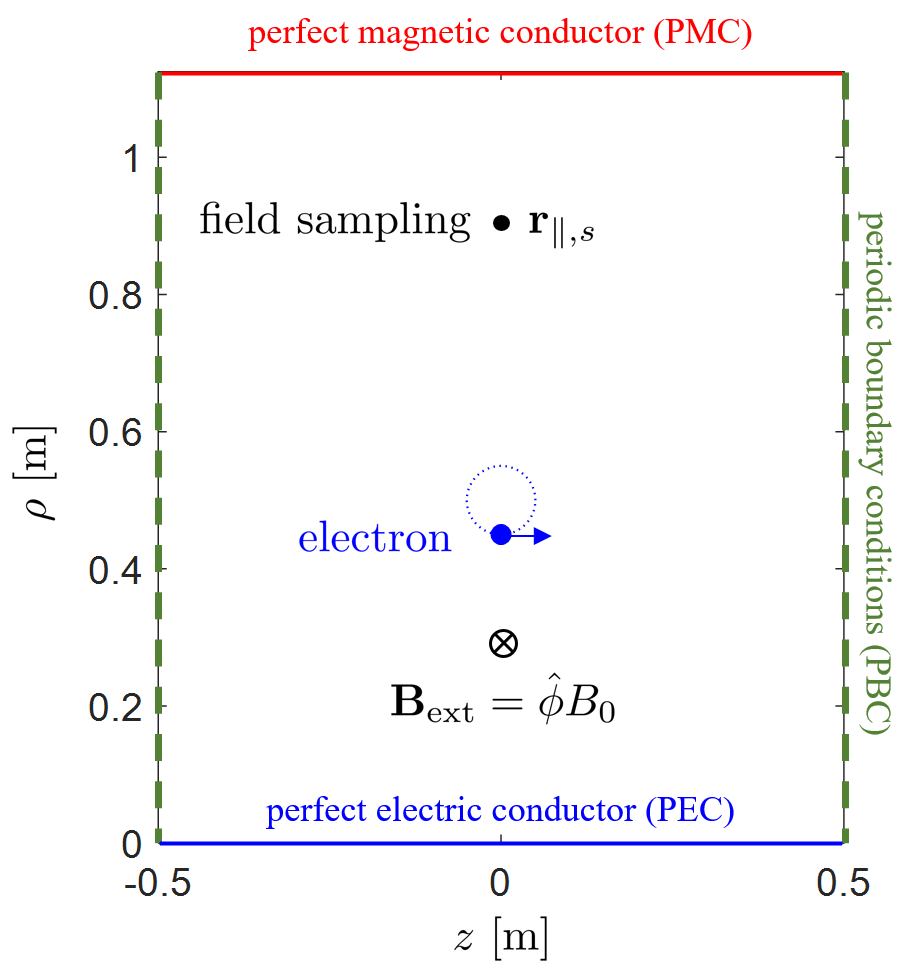}}
\caption{
Geometry for observing gyromotion (a) with and (b) without the radial PML.
}
\label{fig:gyromotion_problem_geometry}
\end{figure}

We next consider a gyromotion of a charged superparticle ($10^{6}$ electrons) in the $z \rho$-plane due to an external azimuthal magnetic field of constant magnitude, $\mathbf{B}'_{ext}=\hat{\phi} \, 8.53\times 10^{-4}$ [T], as illustrated in Fig.~\ref{fig:gyromotion_problem_geometry}.
The superparticle is initialized at $(\rho,z)=(0,0.45)$ [m] with initial velocity $\hat{z}0.025c$ .
The superparticle rotates with Larmor radius $r_{L}=0.05$ [m] and Larmor frequency $f_{L}=1.5 \times 10^{8}$ [Hz].
Electric and magnetic fields are sampled at $(\rho,z)=(0.9,0)$ [m] during ten periodic gyromotions, and their behaviors in time and frequency domains are observed.
Two cases are compared: one has a radial PML whereas the other does not.
Both cases use the same irregular mesh consisting of $7,041$ nodes, $20,810$ edges, and $13,770$ cells.
The case without the radial PML is equivalent to a cylindrical cavity terminated by a lateral PMC boundary condition.

\begin{figure}
\centering
\subfloat[]
{\includegraphics[width=.5\linewidth]{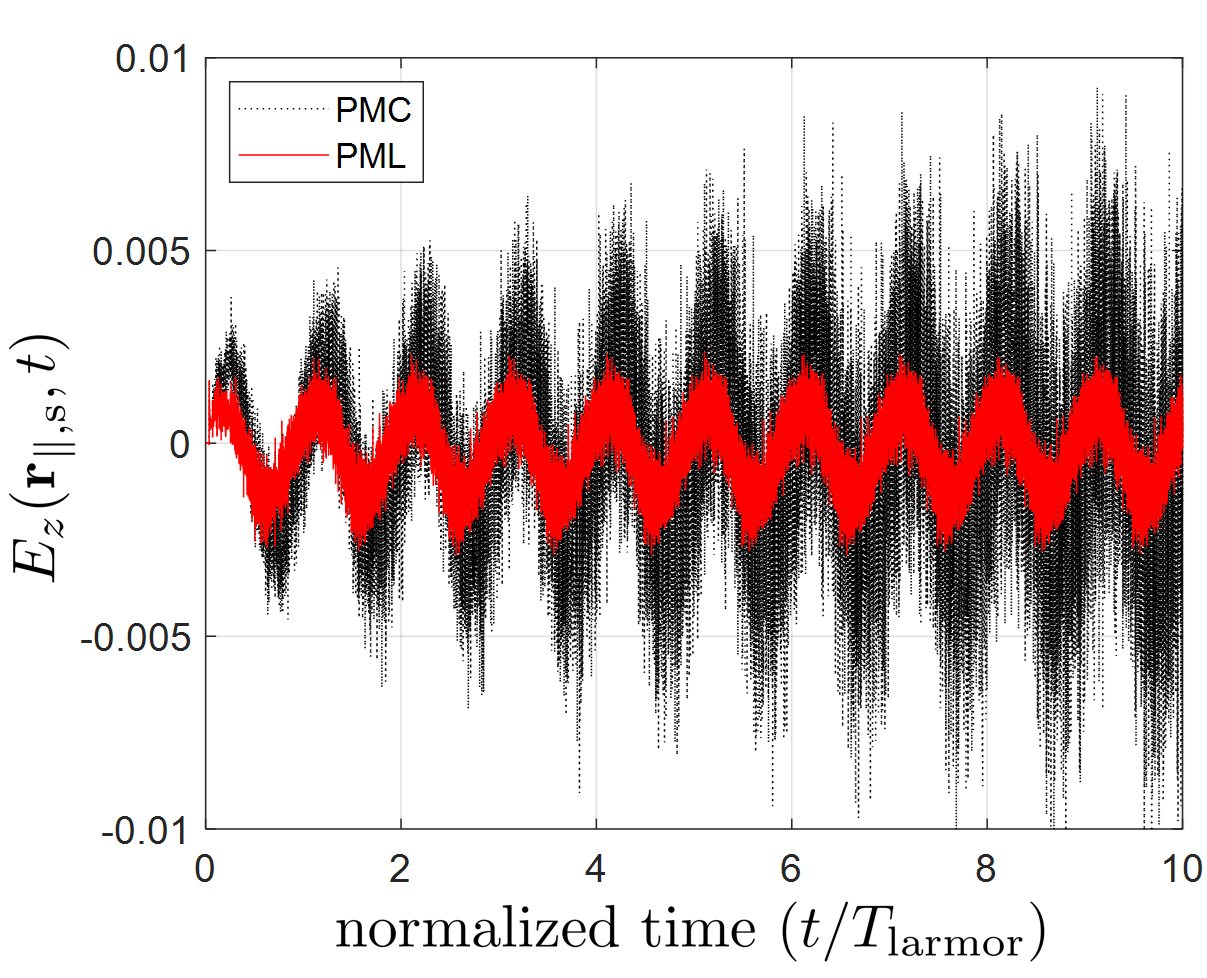}}
\\
\subfloat[]
{\includegraphics[width=.5\linewidth]{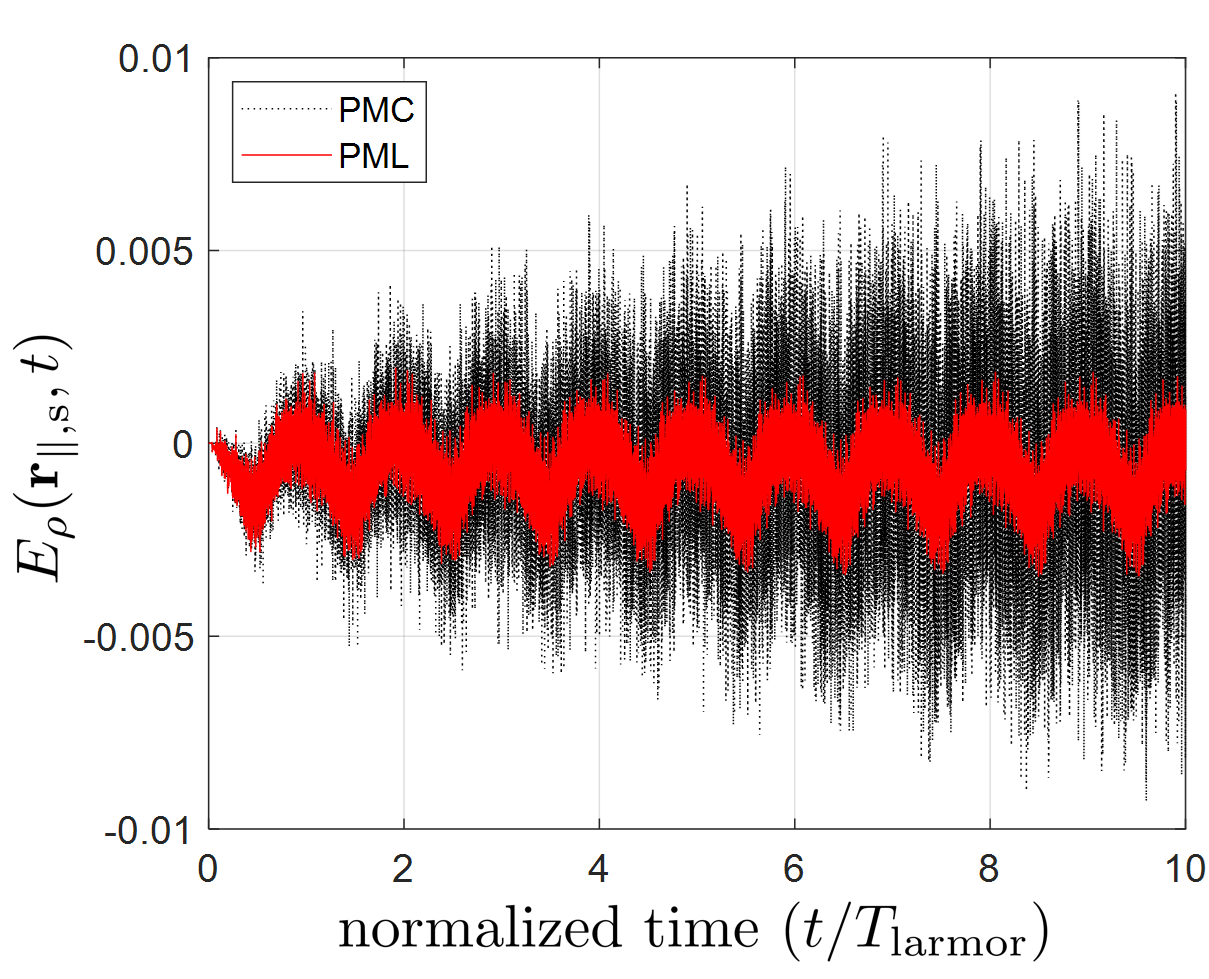}}
\\
\subfloat[]
{\includegraphics[width=.5\linewidth]{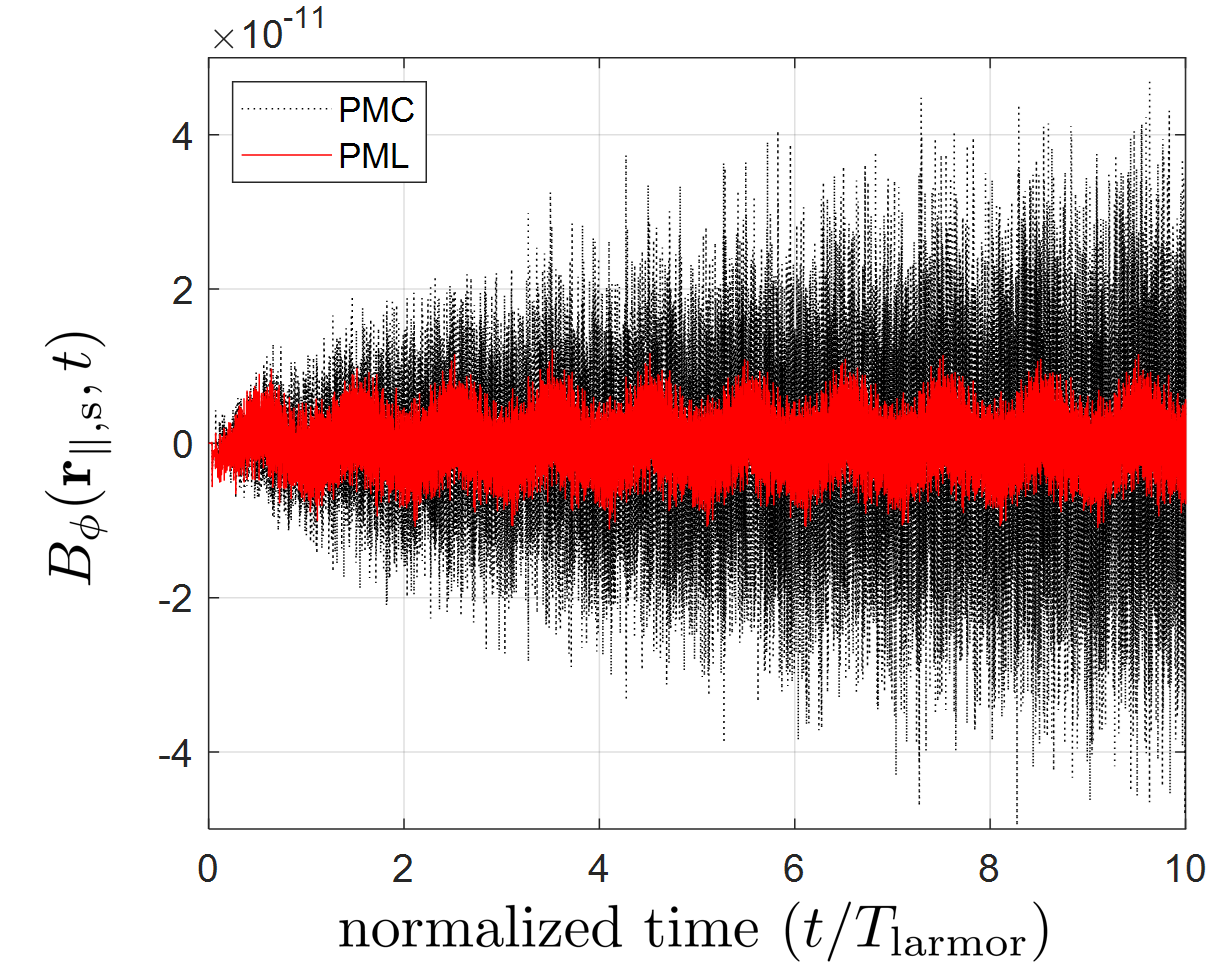}}
\caption{
Temporal evolution of $(E_{z},E_{\rho},B_{\phi})$ field values sampled at $(\rho,z)=(0.9,0)$.
}
\label{fig:time_domain_two_case}
\end{figure}

\begin{figure}
\centering
\subfloat[]
{\includegraphics[width=.5\linewidth]{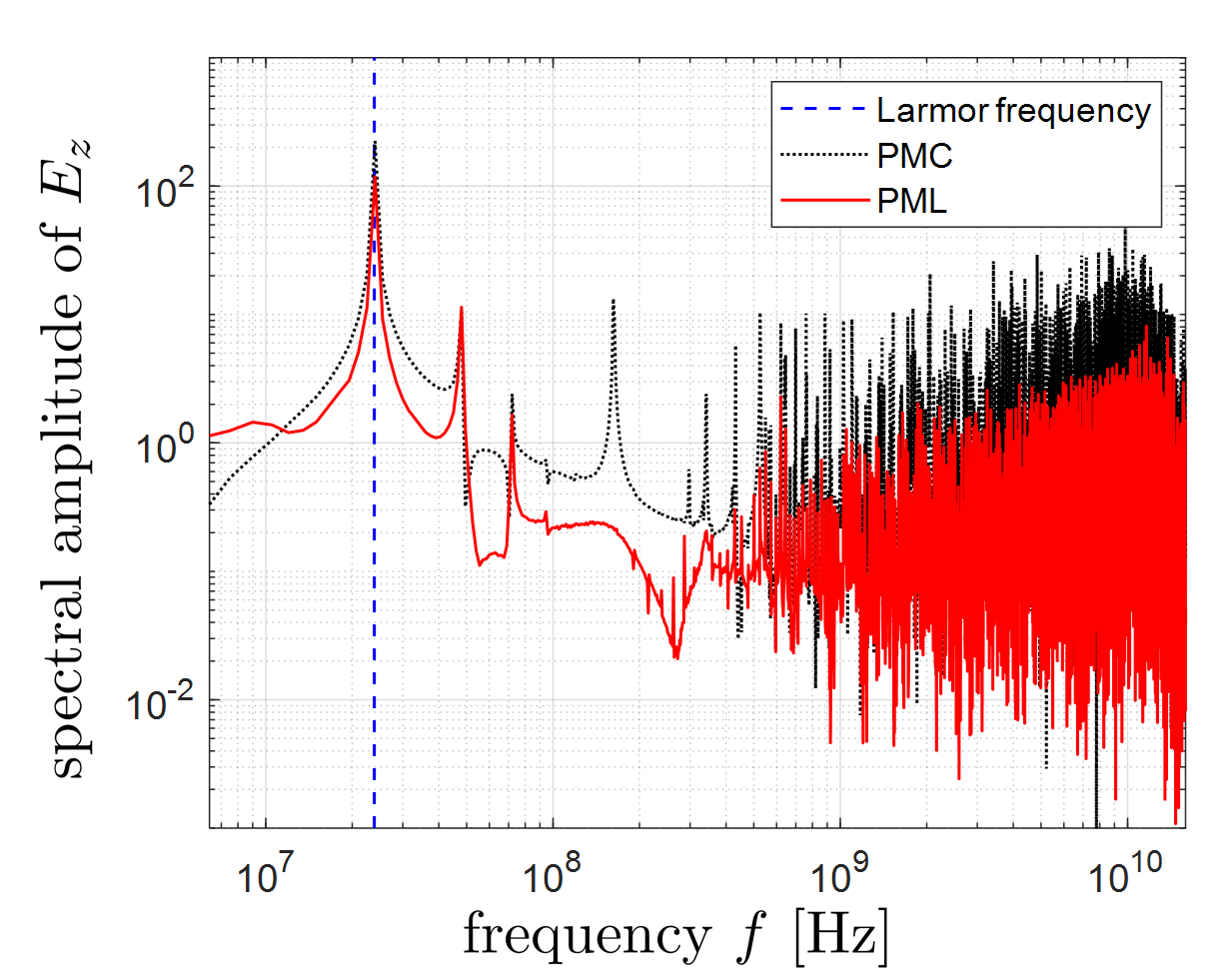}}
\\
\subfloat[]
{\includegraphics[width=.5\linewidth]{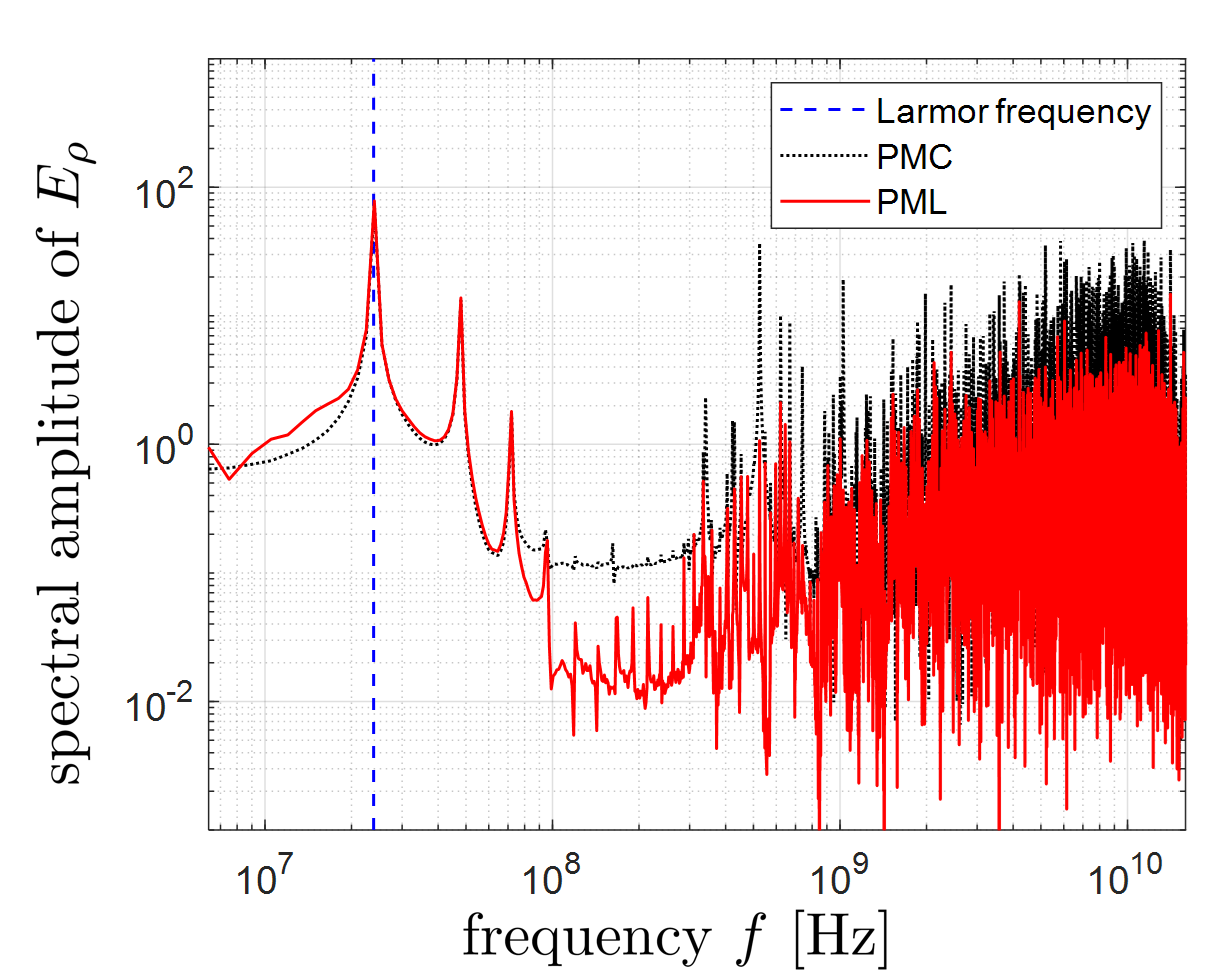}}
\\
\subfloat[]
{\includegraphics[width=.5\linewidth]{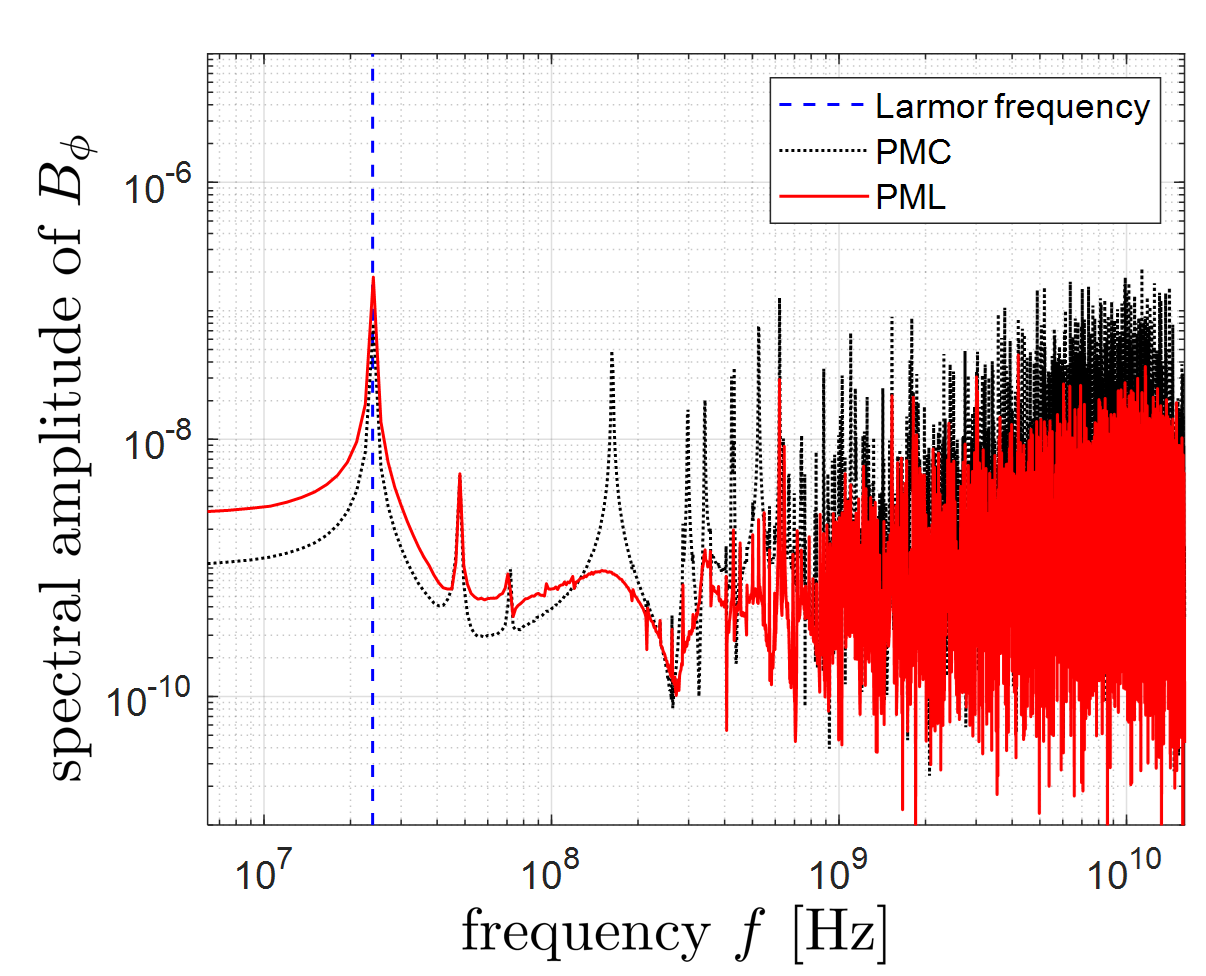}}
\caption{
Spectral amplitudes of the $(E_{z},E_{\rho},B_{\phi})$ fields sampled at  $(\rho,z)=(0.9,0)$.
}
\label{fig:spectral_amp_two_case}
\end{figure}

Fig. \ref{fig:time_domain_two_case} shows sampled values of the field components $E_z$, $E_{\rho}$, and $B_{\phi}$ at a point inside the computational domain as indicated in Fig.~\ref{fig:gyromotion_problem_geometry}, as a function of time for the two cases.
Oscillations in both the electric and magnetic fields are observed, which are a consequence of the gyromotion of the charged particle.
It is interesting to observe from Fig.~\ref{fig:time_domain_two_case} that in the simulation without PML (denoted as PMC), the field amplitude exhibits a secular grow in time since the electromagnetic radiation is not absorbed by the PML and remains trapped in simulation domain cavity.
On the other hand, the peak signal amplitude in the simulation with PML remains steady, which indicates that the radiation is well absorbed by the PML.
The corresponding spectral amplitudes are depicted in Fig. \ref{fig:spectral_amp_two_case}.
The fundamental frequency in the spectrum is in excellent agreement with the theoretical Larmor frequency.
In the case without PML, a large resonant peak appears at approximately $1.62\times 10^{8}$ [Hz], which corresponds to a resonant mode of the cylindrical cavity.
The high-frequency noise can be attributed to particle noise and can be reduced by using higher-order shape functions~\cite{pinto2016electromagnetic,campos2015towards}.

\section{Conclusions}
We have developed a novel electromagnetic particle-in-cell code, BORPIC++, for conducting fully kinetic plasma simulations in complex body-of-revolution geometries with radial perfectly matched layers, axial periodic conditions, and dual field polarizations.
This algorithm employs a combination of field scalings and coordinate mapping techniques to simplify the Maxwell field problem in a cylindrical system.
By doing so, we were able to transform it into a Cartesian finite element Maxwell solver operating in the meridian plane, which enables greater adaptability in handling complex geometries.
This approach effectively eliminates the coordinate singularity at the axis of symmetry and enables easy adaptation of the code to handle different polarizations.
Furthermore, we exploited symmetries inherent in Maxwell's equations to decompose the problem into two dual polarizations.
These states possess isomorphic representations that facilitate code reuse and simplify implementation.
Importantly, we incorporated a cylindrical perfectly matched layer as a boundary condition in the radial direction.
This allows us to simulate open space problems effectively.
We provided several numerical examples to demonstrate the initial application of the algorithm and showcase its capabilities in handling complex geometries and polarizations.

Our next objective is to employ the BORPIC++ code for performing collective electromagnetic plasma simulations. These simulations will focus on ion velocity rings and various types of beams in an open space. The main goal is to accurately model the nonlinear interactions between charged plasma particles and waves as well as capture the time evolution of these phenomena over long durations.
To effectively address large-scale problems, we can employ a parallel computing method in BORPIC++ simulations. Although the update scheme of the standard finite-element time-domain (FETD) method is inherently implicit, making it challenging to adapt it to parallel computing in general, our FETD scheme takes advantage of the local sparse approximate inverse (SPAI), which significantly facilitates the algorithm parallelization.
Furthermore, the previous research has already demonstrated that obtaining the SPAI mass matrix can be done in an embarrassingly parallel manner \cite{Kim2011Parallel}; thus, efficient parallelization of BORPIC++ is expected to be relatively straightforward.

We have implemented the BORPIC++ algorithm in MATLAB, which is available on GitHub \cite{Github_BORPIC_PP}.

\section*{Acknowledgments}
The authors thank the reviewers for their valuable feedback. This work was funded in part by the POSCO Holdings Inc. Grant 20233355, the U.S. Department of Energy Grant DE-SC0022982 through the NSF/DOE Partnership in Basic Plasma Science and Engineering, and the Ohio Supercomputer Center Grant PAS-0061.

\appendix
\section{Equivalence between Maxwell's equations in the original and rescaled cylindrical coordinate systems}\label{equivalence}
\subsection{Faraday's law}\label{app:FL}
Let us first examine the Faraday's law.
The original Faraday's law in the cylindrical coordinate system is written by
\begin{flalign}
\nabla\times\mathbf{E}=-\frac{\partial}{\partial t}\mathbf{B}
\end{flalign}
where
\begin{flalign}
\text{LHS}
&=
\hat{\rho}\frac{1}{\rho}\left(\frac{\partial E_{z}}{\partial \phi} - \rho \frac{\partial E_{\phi}}{\partial z}\right)
+
\hat{\phi}\left(\frac{\partial E_{\rho}}{\partial z} - \frac{\partial E_{z}}{\partial \rho}\right)
+
\hat{z}\frac{1}{\rho}\left(E_\phi + \rho \frac{\partial E_{\phi}}{\partial \rho} - \frac{\partial E_{\rho}}{\partial \phi}\right),
\\
\text{RHS}
&=
-
\hat{\rho}\frac{\partial}{\partial t}B_{\rho}
-
\hat{\phi}\frac{\partial}{\partial t}B_{\phi}
-
\hat{z}\frac{\partial}{\partial t}B_{z}.
\end{flalign}
Thus, we have three component equations as follows:
\begin{flalign}
\frac{1}{\rho}\left(\frac{\partial E_{z}}{\partial \phi} - \rho \frac{\partial E_{\phi}}{\partial z}\right) &= -\frac{\partial}{\partial t}B_{\rho}, \label{eqn:FL_z_comp_org}\\
\left(\frac{\partial E_{\rho}}{\partial z} - \frac{\partial E_{z}}{\partial \rho}\right) &= -\frac{\partial}{\partial t}B_{\phi},
\label{eqn:FL_rho_comp_org}\\
\frac{1}{\rho}\left(E_\phi + \rho \frac{\partial E_{\phi}}{\partial \rho} - \frac{\partial E_{\rho}}{\partial \phi}\right) &= -\frac{\partial}{\partial t}B_{z}.
\label{eqn:FL_phi_comp_org}
\end{flalign}
The transformed Maxwell's equation is written by
\begin{flalign}
\nabla'\times\mathbf{E}'=-\frac{\partial}{\partial t}\mathbf{B}'
\end{flalign}
where
\begin{flalign}
\text{LHS}
&=
\hat{\rho}\left(\frac{\partial {E'}_{z}}{\partial \phi} - \frac{\partial {E'}_{\phi}}{\partial z}\right)
+
\hat{\phi}\left(\frac{\partial {E'}_{\rho}}{\partial z} - \frac{\partial {E'}_{z}}{\partial \rho}\right)
+
\hat{z}\left(\frac{\partial {E'}_{\phi}}{\partial \rho} - \frac{\partial {E'}_{\rho}}{\partial \phi}\right),
\\
\text{RHS}
&=
-
\hat{\rho}\frac{\partial}{\partial t}{B'}_{\rho}
-
\hat{\phi}\frac{\partial}{\partial t}{B'}_{\phi}
-
\hat{z}\frac{\partial}{\partial t}{B'}_{z}.
\end{flalign}
Similarly, we have three component equations
\begin{flalign}
\frac{\partial {E'}_{z}}{\partial \phi} - \frac{\partial {E'}_{\phi}}{\partial z}
&= -\frac{\partial}{\partial t}{B'}_{\rho}, \\
\frac{\partial {E'}_{\rho}}{\partial z} - \frac{\partial {E'}_{z}}{\partial \rho}
&= -\frac{\partial}{\partial t}{B'}_{\phi},\\
\frac{\partial {E'}_{\phi}}{\partial \rho} - \frac{\partial {E'}_{\rho}}{\partial \phi}
&= -\frac{\partial}{\partial t}{B'}_{z}.
\end{flalign}
When using the relation in \eqref{eqn:transfE} and \eqref{eqn:transfB}, the above three equations can be rewritten in terms of original Maxwellian field variables as
\begin{flalign}
\frac{\partial {E}_{z}}{\partial \phi} - \frac{\partial \left(\rho{E}_{\phi}\right)}{\partial z}
&= -\frac{\partial}{\partial t}(\rho{B}_{\rho}),
\label{eqn:FL_z_comp_trs}
\\
\frac{\partial {E}_{\rho}}{\partial z} - \frac{\partial {E}_{z}}{\partial \rho}
&= -\frac{\partial}{\partial t}{B}_{\phi},
\label{eqn:FL_rho_comp_trs}\\
\frac{\partial \left(\rho{E}_{\phi}\right)}{\partial \rho} - \frac{\partial E_{\rho}}{\partial \phi}
&= -\frac{\partial}{\partial t}(\rho{B}_{z}),
\label{eqn:FL_phi_comp_trs}
\end{flalign}
which are identical to eqs. \eqref{eqn:FL_z_comp_org}, \eqref{eqn:FL_rho_comp_org}, and \eqref{eqn:FL_phi_comp_org}, respectively.

Since the transformed Faraday's law can be compactly written by
\begin{flalign}
\nabla'\times \bar{\bar{U}}\mathbf{E} = -\frac{\partial}{\partial t}\bar{\bar{V}}\cdot\mathbf{B},
\end{flalign}
it can be generalized that any equation taking the form of the above, i.e.,
\begin{flalign}
\nabla'\times \bar{\bar{U}}\mathbf{P'} = -\frac{\partial}{\partial t}\bar{\bar{V}}\cdot\mathbf{Q'},
\end{flalign}
for arbitrary rescaled vector fields $\mathbf{P'}$ and $\mathbf{Q'}$ is identical to
\begin{flalign}
\nabla \times \mathbf{P} = -\frac{\partial}{\partial t} \mathbf{Q}.
\end{flalign}

\subsection{Ampere's law}
The original Ampere's law in the cylindrical coordinate system is written by
\begin{flalign}
\nabla\times\mathbf{H}=\frac{\partial}{\partial t}\mathbf{D}+\mathbf{J}.
\label{eqn:org_AL}
\end{flalign}
And the rescaled Ampere's law given by
\begin{flalign}
\nabla'\times\mathbf{H}'=\frac{\partial}{\partial t}\mathbf{D}'+\mathbf{J}'
\label{eqn:res_AL}
\end{flalign}
can be rewritten with the use of \eqref{eqn:transfH}, \eqref{eqn:transfD}, and \eqref{eqn:transfJ} by
\begin{flalign}
\nabla'\times \bar{\bar{U}}\cdot\mathbf{H} = \frac{\partial}{\partial t}\bar{\bar{V}}\cdot\mathbf{D}+\bar{\bar{V}}\cdot\mathbf{J}.
\end{flalign}
Hence, \eqref{eqn:org_AL} and \eqref{eqn:res_AL} are identical to each other in accordance with the conclusion shown in \ref{app:FL}.

We further examine the equivalence of the Ampere's law written in term of electric field intensity and magnetic flux density.
For the simplicity, we assume the background is the vacuum.
In the original Maxwell's equations, the Ampere's law can be written by
\begin{flalign}
\nabla\times \mu_{0}^{-1}\mathbf{B} = \frac{\partial}{\partial t} \epsilon_{0} \mathbf{E} + \mathbf{J},
\label{eqn:AL_E_B}
\end{flalign}
and the rescaled version of it is given by
\begin{flalign}
\nabla'\times \bar{\bar{\mu}}'^{-1}\cdot\mathbf{B'} = \frac{\partial}{\partial t} {\bar{\bar{\epsilon}}'}\cdot \mathbf{E'} + \mathbf{J'}.
\end{flalign}
Using \eqref{eqn:eff_permittivity}, \eqref{eqn:eff_permeability}, \eqref{eqn:transfH}, \eqref{eqn:transfD}, and \eqref{eqn:transfJ}, one can rewrite the above as
\begin{flalign}
\nabla'\times \mu_0^{-1}\bar{\bar{R}}^{-1}\cdot\bar{\bar{V}}\cdot\mathbf{B} = \frac{\partial}{\partial t} \epsilon_0 \bar{\bar{R}}\cdot\bar{\bar{U}}\cdot \mathbf{E} + {\bar{\bar{V}}}\cdot\mathbf{J}.
\end{flalign}
Since $\bar{\bar{R}}^{-1}\cdot\bar{\bar{V}}=\bar{\bar{U}}$ and $\bar{\bar{R}}\cdot\bar{\bar{U}}=\bar{\bar{V}}$, we finally have the equation below:
\begin{flalign}
\nabla'\times \mu_0^{-1}\bar{\bar{U}}\cdot\mathbf{B} = \frac{\partial}{\partial t} \epsilon_0 {\bar{\bar{V}}}\cdot \mathbf{E} + {\bar{\bar{V}}}\cdot\mathbf{J}.
\end{flalign}
Thus, the above equation is the same as the original Ampere's law written in term of electric field intensity and magnetic flux density in \eqref{eqn:AL_E_B}.

\section{Singularity-free discrete Hodge matrices}\label{app:singularity}
When the background is vacuum, one can deduce from \eqref{eqn:eff_permittivity}, \eqref{eqn:eff_permeability}, \eqref{eqn:rescale_GL}, and \eqref{eqn:rescale_GLM}:
\begin{flalign}
{D'}_{\rho} = \epsilon_0 \rho {E'}_{\rho}, \quad
{D'}_{\phi} = \epsilon_0 \rho^{-1} {E'}_{\phi}, \quad
{D'}_{z} = \epsilon_0 \rho {E'}_{z}
\end{flalign}
and
\begin{flalign}
{B'}_{\rho} = \mu_0 \rho {H'}_{\rho}, \quad
{B'}_{\phi} = \mu_0 \rho^{-1} {H'}_{\phi}, \quad
{B'}_{z} = \mu_0 \rho {H'}_{z}.
\end{flalign}
From the above relations, it is clear that TE\textsuperscript{$\phi$} polarized field components have no singularity in the mapping from $({E'}_{z},{E'}_{\rho},{B'}_{\phi})$, defined in a primal mesh, to $({D'}_{z},{D'}_{\rho},{H'}_{\phi})$, defined in a dual mesh, since $({E'}_{z},{E'}_{\rho},{B'}_{\phi})$ are always scaled by $\rho$, i.e.,
\begin{flalign}
\left[
\begin{matrix}
{D'}_{z}  & {D'}_{\rho}
\end{matrix}
\right]
=
\epsilon_0
\left[
\begin{matrix}
\rho & 0 \\ 0 & \rho
\end{matrix}
\right]
\cdot
\left[
\begin{matrix}
{E'}_{z}  & {E'}_{\rho}
\end{matrix}
\right],
\quad
{H'}_{\phi}
=
\left(\mu_0^{-1}\rho\right){B'}_{\phi}.
\end{flalign}
While encoding the above effective medium, the discrete Hodge matrix for the TE\textsuperscript{$\phi$} polarized field components can be defined as
\begin{flalign}
\left[\star_{\epsilon}\right]_{i,j}
&\equiv
\sum_{n=1}^{N_{2}}
\iint_{T_{n}}\,
\mathbf{W}^{(1)}_{i}
\cdot
\bar{\bar{\epsilon}}' \cdot
\mathbf{W}^{(1)}_{j} \, dA,
\nonumber \\
&=
\sum_{n=1}^{N_{2}}
\iint_{T_{n}}
\,
\left[
\begin{matrix}
\hat{\rho}\cdot\mathbf{W}^{(1)}_{i} \\ \hat{z}\cdot\mathbf{W}^{(1)}_{i}
\end{matrix}
\right]^{T}
\cdot
\epsilon_0
\left[
\begin{matrix}
\rho & 0 \\
0 & \rho \\
\end{matrix}
\right]
\cdot
\left[
\begin{matrix}
\hat{\rho}\cdot\mathbf{W}^{(1)}_{j} \\ \hat{z}\cdot\mathbf{W}^{(1)}_{j}
\end{matrix}
\right]
\, dA,
\label{app_eqn:TE_Hodge_eps}
\end{flalign}
\begin{flalign}
 \left[\star_{\mu^{-1}}\right]_{i,j} &\equiv
\sum_{n=1}^{N_{2}}
\iint_{T_{n}}
 \,
 \mathbf{W}^{(2)}_{i}  \cdot
 \bar{\bar{\mu}}'^{-1} \cdot
\mathbf{W}^{(2)}_{j} \, dA
\nonumber \\
&=
\sum_{n=1}^{N_{2}}
\iint_{T_{n}}
\,
\left(\hat{\phi}\cdot\mathbf{W}^{(2)}_{i}\right)
\mu_0^{-1}\rho
\left(\hat{\phi}\cdot\mathbf{W}^{(2)}_{j}\right) \, dA,
\label{app_eqn:TE_Hodge_mu}
\end{flalign}
where $T_{n}$ is $n$-th triangle of the primal mesh.
As observed in the expression, even in triangles that takes their nodes located on $z$-axis ($\rho=0$), the discrete Hodge matrices have no singularity.

However, this is not the case for TM\textsuperscript{$\phi$} polarized field components, which can be easily identified by doing the similar procedure.
Specifically, when discretizing ${E'}_\phi$, ${B'}_z$, ${B'}_\rho$ on a primal mesh and transforming from ${E'}_\phi$, ${B'}_z$, ${B'}_\rho$ to ${D'}_\phi$, ${H'}_z$, ${H'}_\rho$, one can show that ${E'}_\phi$, ${B'}_z$, ${B'}_\rho$ are always scaled by $\rho^{-1}$, which produces a singularity at $\rho=0$, i.e.,
\begin{flalign}
\left[
\begin{matrix}
{H'}_{z}  & {H'}_{\rho}
\end{matrix}
\right]
=
\mu_0^{-1}
\left[
\begin{matrix}
\rho^{-1} & 0 \\ 0 & \rho^{-1}
\end{matrix}
\right]
\cdot
\left[
\begin{matrix}
{B'}_{z}  & {B'}_{\rho}
\end{matrix}
\right],
\quad
{D'}_{\phi}
=
\left(\epsilon_0\rho^{-1}\right){E'}_{\phi}.
\end{flalign}
Thus, the discrete Hodge matrix cannot be well-defined since the effective medium contains a $\rho^{-1}$ factor such that elements of the discrete Hodge matrix contributed from triangles touching the symmetric axis ($\rho=0$) diverge.
To resolve the singularity issue when dealing with TM\textsuperscript{$\phi$} polarized field components, we treat ${D'}_\phi$, ${H'}_z$, ${H'}_\rho$ as primal quantities and ${E'}_\phi$, ${B'}_z$, ${B'}_\rho$ as dual quantities.
With this treatment, we now map from ${D'}_\phi$, ${H'}_z$, ${H'}_\rho$ to ${E'}_\phi$, ${B'}_z$, ${B'}_\rho$, which is devoid of the singularity since
\begin{flalign}
\left[
\begin{matrix}
{B'}_{z}  & {B'}_{\rho}
\end{matrix}
\right]
=
\mu_0
\left[
\begin{matrix}
\rho & 0 \\ 0 & \rho
\end{matrix}
\right]
\cdot
\left[
\begin{matrix}
{H'}_{z}  & {H'}_{\rho}
\end{matrix}
\right],
\quad
{E'}_{\phi}
=
\left(\epsilon_0^{-1}\rho\right){D'}_{\phi}.
\end{flalign}
Thus, the resulting Hodge matrix can be also well-defined by
\begin{flalign}
\left[\star_{\epsilon^{-1}}\right]_{i,j}
&=
\sum_{n=1}^{N_{2}}
\iint_{T_{n}}
\mathbf{W}^{(2)}_{i}
\cdot
\bar{\bar{\epsilon}}'^{-1} \cdot
\mathbf{W}^{(2)}_{j}
\nonumber \\
&=
\sum_{n=1}^{N_{2}}
\iint_{T_{n}}  \,
\left(\hat{\phi}\cdot\mathbf{W}^{(2)}_{i}\right)
\epsilon_0^{-1}
\rho
\left(\hat{\phi}\cdot\mathbf{W}^{(2)}_{j}\right) \, dA,
\label{app_eqn:TM_Hodge_eps}
\, dA,
\\
\left[\star_{\mu}\right]_{i,j}
&=
\sum_{n=1}^{N_{2}}
\iint_{T_{n}}
\mathbf{W}^{(1)}_{i}
\cdot
\bar{\bar{\mu}}' \cdot
\mathbf{W}^{(1)}_{j}
\, dA
\nonumber \\
&=
\sum_{n=1}^{N_{2}}
\iint_{T_{n}}\,
\left[
\begin{matrix}
\hat{\rho}\cdot\mathbf{W}^{(1)}_{i} \\ \hat{z}\cdot\mathbf{W}^{(1)}_{i}
\end{matrix}
\right]^{T}
\cdot
\mu_0
\left[
\begin{matrix}
\rho & 0 \\
0 & \rho \\
\end{matrix}
\right]
\cdot
\left[
\begin{matrix}
\hat{\rho}\cdot\mathbf{W}^{(1)}_{j} \\ \hat{z}\cdot\mathbf{W}^{(1)}_{j}
\end{matrix}
\right] \, dA.
\label{app_eqn:TM_Hodge_mu}
\end{flalign}
Another key advantage of the above treatment is the fact that the discrete Hodge matrices for TM\textsuperscript{$\phi$} polarized field components are dual to those for TE\textsuperscript{$\phi$} polarized field components, i.e.,
\begin{flalign}
\frac{1}{\mu_0}\left[\star_{\mu}\right] = \frac{1}{\epsilon_0}\left[\star_{\epsilon}\right], \quad
\epsilon_0\left[\star_{\epsilon^{-1}}\right] = \mu_0 \left[\star_{\mu^{-1}}\right].
\end{flalign}
Consequently, it is not necessary to perform a separate implementation for the discrete Hodge matrices for TM\textsuperscript{$\phi$} polarized field components once those for TE\textsuperscript{$\phi$} polarized field components are complete.

\section{Radial PML implementation}\label{app:radialPML}

We first consider the TE\textsuperscript{$\phi$} case.
Assuming a background medium with constitutive tensors ${\bar{\bar{\epsilon}}'}$ and ${\bar{\bar{\mu}}'}$ from \eqref{eqn:eff_permittivity} and \eqref{eqn:eff_permeability}, respectively, the radial PML can be implemented by assuming constitutive tensors ${\bar{\bar{\epsilon}}'}^{(\text{PML})}$ and ${\bar{\bar{\mu}}'}^{(\text{PML})}$ defined over the outer annular region encircling the cylindrical domain of interest~\cite{teixeira1997systematic} and given by
\begin{flalign}
\label{eqn:PMLtensors}
{\bar{\bar{\epsilon}}'}^{(\text{PML})}
=
\left(\text{det}\bar{\bar{S}}\right)^{-1}
\bar{\bar{S}}
\cdot
{\bar{\bar{\epsilon}}'}
\cdot
\bar{\bar{S}}
,
\quad
{\bar{\bar{\mu}}'}^{(\text{PML})}
=
\left(\text{det}\bar{\bar{S}}\right)^{-1}
\bar{\bar{S}}
\cdot
{\bar{\bar{\mu}}'}
\cdot
\bar{\bar{S}}
\end{flalign}
where the tensor $\bar{\bar{S}}$ is expressed in the frequency domain ($e^{j \omega t}$ convention) as
\begin{flalign}
\label{eqn:Spml}
\bar{\bar{S}} = \hat{z}\hat{z}+\hat{\rho}\hat{\rho}\left(\frac{1}{s_{\rho}}\right)+\hat{\phi}\hat{\phi}\left(\frac{\rho}{{\tilde \rho}}\right),
\end{flalign}
with
\begin{flalign}
s_{\rho}(\rho) = 1+\frac{\sigma_{\rho}(\rho)}{j\omega\epsilon_0},
\\
{\tilde \rho} = \int_0^{\rho} s_{\rho}(\rho) d\rho,
\end{flalign}
where $\sigma_{\rho}(\rho) \geq 0$ is an artificial conductivity factor to produce wave absorption within the PML. Outside the PML region (physical domain), $\sigma_{\rho}=0$ recovers the original medium. Inside the PML region, $\sigma_{\rho}(\rho)$ is typically chosen to have a polynomial taper growth.
 The PML tensors in eq.~$(\ref{eqn:PMLtensors})$ provide reflection absorption of the waves in the continuum limit. For the range of parameters considered here,
 $|{\tilde \rho} - \rho|/\rho \ll 1$, so we will let $\rho/{\tilde \rho} \rightarrow 1$ in eq.~(\ref{eqn:Spml}). This approximation simplifies the PML implementation considerably while having negligible effects on performance. Consequently,
\begin{flalign}
{\bar{\bar{\epsilon}}'}^{(\text{PML})}
=
\epsilon_0
\left[
\begin{matrix}
\rho s_{\rho} & 0 & 0 \\
0 & \frac{\rho}{s_{\rho}} & 0 \\
0 & 0& \frac{s_{\rho}}{\rho}\\
\end{matrix}
\right],
\quad
{\bar{\bar{\mu}}'}^{(\text{PML})}
=
\mu_0
\left[
\begin{matrix}
\frac{s_{\rho}}{\rho} & 0 & 0 \\
0 & \frac{1}{\rho s_{\rho}} & 0 \\
0 & 0& \rho s_{\rho}\\
\end{matrix}
\right].
\end{flalign}
where
\begin{flalign}
\left[{\bar{\bar{\epsilon}}'}^{(\text{PML})}\right]_{zz}
&=
\epsilon_0 \rho s_{\rho}
=
\epsilon_0
\rho
\left(1+\frac{\sigma_{\rho}}{j\omega\epsilon_0}\right),
\\
\left[{\bar{\bar{\epsilon}}'}^{(\text{PML})}\right]_{\rho\rho}
&=
\frac{\epsilon_0 \rho}{s_{\rho}}
=
\epsilon_0
\rho
\frac{j\omega\epsilon_0}{j\omega\epsilon_0+\sigma_{\rho}},
\\
\left[{\bar{\bar{\epsilon}}'}^{(\text{PML})}\right]_{\phi\phi}
&=
\frac{\epsilon_0 s_{\rho}}{\rho}
=
\frac{\epsilon_0}{\rho}
\left(1+\frac{\sigma_{\rho}}{j\omega\epsilon_0}\right),
\\
\left[{\bar{\bar{\mu}}'}^{(\text{PML})}\right]_{zz}
&=
\frac{\mu_0 s_{\rho}}{\rho}
=
\frac{\mu_0}{\rho}
\left(1+\frac{\sigma_{\rho}}{j\omega\epsilon_0}\right),
\\
\left[{\bar{\bar{\mu}}'}^{(\text{PML})}\right]_{\rho\rho}
&=
\frac{\mu_0\rho}{s_{\rho}}
=
\frac{\mu_0}{\rho}
\frac{j\omega\epsilon_0}{j\omega\epsilon_0+\sigma_{\rho}},
\\
\left[{\bar{\bar{\mu}}'}^{(\text{PML})}\right]_{\phi\phi}
&=
\mu_0\rho s_{\rho}
=
\mu_0 \rho
\left(1+\frac{\sigma_{\rho}}{j\omega\epsilon_0}\right).
\end{flalign}

For implementation in finite element algorithm, we follow similar steps as described in~\cite{Donderici2008Mixed}.
First, we assume that the above permittivity and permeability are element-wise constant; hence, the above can be rewritten in the $k$-th mesh element ${\mathcal{T}_{k}}$ as
\begin{flalign}
\left[{\bar{\bar{\epsilon}}'}^{(\text{PML})}(k)\right]_{zz}
&=
\epsilon_0
\rho(k)
\left(1+\frac{\sigma_{\rho}(k)}{j\omega\epsilon_0}\right),
\\
\left[{\bar{\bar{\epsilon}}'}^{(\text{PML})}(k)\right]_{\rho\rho}
&=
\epsilon_0
\rho
\frac{j\omega\epsilon_0}{j\omega\epsilon_0+\sigma_{\rho}(k)},
\\
\left[{\bar{\bar{\epsilon}}'}^{(\text{PML})}(k)\right]_{\phi\phi}
&=
\frac{\epsilon_0}{\rho(k)}
\left(1+\frac{\sigma_{\rho}(k)}{j\omega\epsilon_0}\right),
\\
\left[{\bar{\bar{\mu}}'}^{(\text{PML})}(k)\right]_{zz}
&=
\frac{\mu_0}{\rho(k)}
\left(1+\frac{\sigma_{\rho}(k)}{j\omega\epsilon_0}\right),
\\
\left[{\bar{\bar{\mu}}'}^{(\text{PML})}(k)\right]_{\rho\rho}
&=
\frac{\mu_0}{\rho(k)}
\frac{j\omega\epsilon_0}{j\omega\epsilon_0+\sigma_{\rho}(k)},
\\
\left[{\bar{\bar{\mu}}'}^{(\text{PML})}(k)\right]_{\phi\phi}
&=
\mu_0 \rho(k)
\left(1+\frac{\sigma_{\rho}(k)}{j\omega\epsilon_0}\right),
\end{flalign}
where $\rho(k)$ denotes the centroid of radial coordinates of $k$-th triangular element, and $\sigma_{\rho}(k)$ is the value of the PML conductivity assigned to the $k$-th element.

We next build the discrete Hodge matrix
\begin{flalign}
\left[\star_{\epsilon}\right]
=
\sum_{k=1}^{N_{2}}
\left[{\bar{\bar{\epsilon}}'}^{(\text{PML})}(k)\right]_{z,z}
\left[\bar{\bar{\mathcal{L}}}_{z}(k)\right]
+
\sum_{k=1}^{N_{2}}
\left[{\bar{\bar{\epsilon}}'}^{(\text{PML})}(k)\right]_{\rho,\rho}
\left[\bar{\bar{\mathcal{L}}}_{\rho}(k)\right]
\end{flalign}
with
\begin{flalign}
\left[\bar{\bar{\mathcal{L}}}_{\zeta}(k)\right]_{i,j}
=
\iint_{\mathcal{T}_{k}}
\left(\hat{\zeta}\cdot\mathbf{W}^{(1)}_{i}(z,\rho)\right)
\left(\hat{\zeta}\cdot\mathbf{W}^{(1)}_{j}(z,\rho)\right)
d A
\end{flalign}
for $\zeta=z,\rho$.
\begin{figure}
\centering
{\includegraphics[width=.5\linewidth]{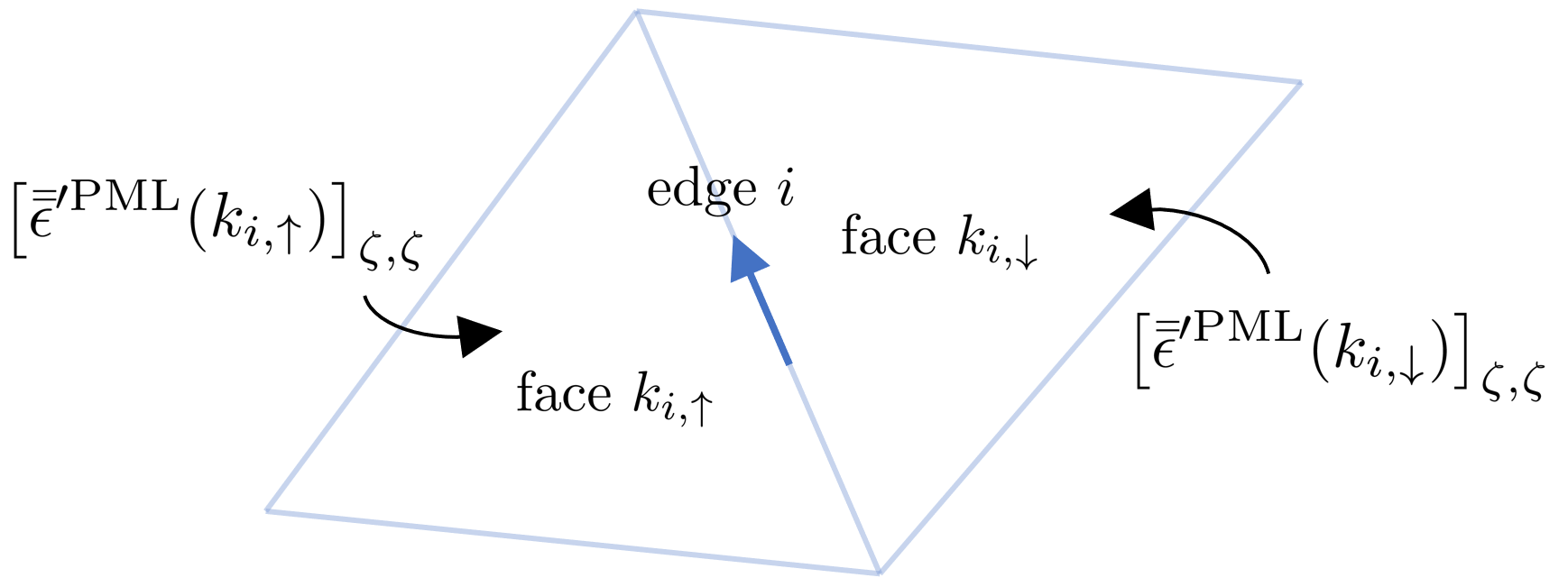}}
\caption{
Indices of two cells that touch $i$-th edge, denoted by $k_{i,\uparrow}$ and $k_{i,\downarrow}$, respectively.
}
\label{fig:PML_face_edge_rel}
\end{figure}
We relate the DoFs of the electric field with those of the electric flux density through the field splitting below~\cite{Donderici2008Mixed}
\begin{flalign}
\left[\tilde{\mathbf{d}}_{\parallel}\right]_{i}
=
\sum_{\zeta=\left\{z,\rho\right\}}
\left(
\left[{\bar{\bar{\epsilon}}'}^{(\text{PML})}(k_{i,\uparrow})\right]_{\zeta,\zeta}
\left[\tilde{\mathbf{e}}_{\parallel}\right]_{i,\zeta,\uparrow}
+
\left[{\bar{\bar{\epsilon}}'}^{(\text{PML})}(k_{i,\downarrow})\right]_{\zeta,\zeta}
\left[\tilde{\mathbf{e}}_{\parallel}\right]_{i,\zeta,\downarrow}
\right)
\label{eqn:d_e_rel_freq}
\end{flalign}
for $\zeta=z,\rho$, where $k_{i,\uparrow}$ and $k_{i,\downarrow}$ are indices of the two mesh elements that touch the $i$-th edge, as illustrated in Fig.~\ref{fig:PML_face_edge_rel}, and
\begin{flalign}
\left[\tilde{\mathbf{e}}_{\parallel}\right]_{i,\zeta,\uparrow}
\equiv
\sum_{j=1}^{N_{1}}
\left[\bar{\bar{\mathcal{L}}}_{\zeta}(k_{i,\uparrow})\right]_{i,j}\left[\tilde{\mathbf{e}}_{\parallel}\right]_{j},
\quad
\left[\tilde{\mathbf{e}}_{\parallel}\right]_{i,\zeta,\downarrow}
\equiv
\sum_{j=1}^{N_{1}}
\left[\bar{\bar{\mathcal{L}}}_{\zeta}(k_{i,\downarrow})\right]_{i,j}\left[\tilde{\mathbf{e}}_{\parallel}\right]_{j}.
\end{flalign}
for $\zeta=z,\rho$.
The tilde on the frequency-domain $\tilde{\mathbf{e}}_{\parallel}$ and $\tilde{\mathbf{d}}_{\parallel}$ is used to distinguish them from the time-domain counterparts.
Expanding both sides of eq.~\eqref{eqn:d_e_rel_freq} and converting the resulting equation to the time domain, we obtain a differential equation of the form
\begin{flalign}
\sum_{p=0}^{N_{p}} r_{\zeta}^{(p)}(k_{i,\xi})
\left[{\mathbf{d}}_{\parallel}^{(p)}\right]_{i,\zeta,\xi}
=
\sum_{p=0}^{N_{p}} q_{\zeta}^{(p)}(k_{i,\xi})
\left[{\mathbf{e}}_{\parallel}^{(p)}\right]_{i,\zeta,\xi}
\label{eqn:d_e_time_domain_eqn}
\end{flalign}
for $\zeta=z,\rho$ and $\xi=\uparrow,\downarrow$. The superscript $p$ indicates a time derivative of order $p$. In the present case,
$N_{p}=1$ with the following coefficients:
\begin{flalign}
r_{z}^{(0)}(k_{i,\uparrow})
& = 0,
\quad
r_{z}^{(1)}(k_{i,\uparrow})
 = 1,
\\
r_{\rho}^{(0)}(k_{i,\uparrow})
& = \frac{\sigma_{\rho}}{\epsilon_{0}},
\quad
r_{\rho}^{(1)}(k_{i,\uparrow})
 = 1,
\\
q_{z}^{(0)}(k_{i,\uparrow})
& = \rho \sigma_{\rho},
\quad
q_{z}^{(1)}(k_{i,\uparrow})
 = \epsilon_0 \rho,
\\
q_{\rho}^{(0)}(k_{i,\uparrow})
& = 0,
\quad
q_{\rho}^{(1)}(k_{i,\uparrow})
 = \epsilon_0 \rho.
\end{flalign}
We discretize \eqref{eqn:d_e_time_domain_eqn} in time by approximating the time derivatives by central differences in time in a nested fashion,
\begin{flalign}
\left[{\mathbf{e}}_{\parallel}^{(p)}\right]_{i,\zeta,\xi}^{n+1}
\approx
-
\left[{\mathbf{e}}_{\parallel}^{(p)}\right]_{i,\zeta,\xi}^{n}
+
\frac{2}{\Delta t}
\left(
\left[{\mathbf{e}}_{\parallel}^{(p-1)}\right]_{i,\zeta,\xi}^{n+1}
-
\left[{\mathbf{e}}_{\parallel}^{(p-1)}\right]_{i,\zeta,\xi}^{n}
\right),
\end{flalign}
and similarly for $\left[{\mathbf{d}}_{\parallel}^{(p)}\right]_{i,\zeta,\xi}^{n+1}$.
Consequently, we arrive at the following relation
\begin{flalign}
\label{eqn:etod}
\left[{\mathbf{d}}_{\parallel}\right]_{i,\zeta,\xi}^{n+1}
=
w_{\zeta}(k_{i,\xi})
\left[{\mathbf{e}}_{\parallel}\right]_{i,\zeta,\xi}^{n+1}
+
\left[{\mathbf{g}}_{\parallel}\right]_{i,\zeta,\xi}^{n}
\end{flalign}
with the auxiliary field ${\mathbf{g}}_{\parallel}$ governed by
\begin{flalign}
\left[{\mathbf{g}}_{\parallel}\right]_{i,\zeta,\xi}^{n}
=
\sum_{p=0}^{N_{p}}
u_{\zeta}^{(p)}(k_{i,\xi})
\left[{\mathbf{d}}_{\parallel}^{(p)}\right]_{i,\zeta,\xi}^{n}
-
\sum_{p=0}^{N_{p}}
v_{\zeta}^{(p)}(k_{i,\xi})
\left[{\mathbf{e}}_{\parallel}^{(p)}\right]_{i,\zeta,\xi}^{n},
\end{flalign}
\begin{flalign}
w_{\zeta}(k_{i,\xi})
&=
\Lambda^{-1}(k_{i,\xi})\left(q_{\zeta}^{(0)}(k_{i,\xi})+
q_{\zeta}^{(1)}(k_{i,\xi})\tau\right),
\\
u_{\zeta}^{(0)}(k_{i,\xi})
&=
\Lambda^{-1}(k_{i,\xi}) r_{\zeta}^{(1)}(k_{i,\xi})\tau,
\quad
u_{\zeta}^{(1)}(k_{i,\xi})
=
\Lambda^{-1}(k_{i,\xi}) r_{\zeta}^{(1)}(k_{i,\xi}),
\\
v_{\zeta}^{(0)}(k_{i,\xi})
&=
\Lambda^{-1}(k_{i,\xi}) q_{\zeta}^{(1)}(k_{i,\xi})\tau,
\quad
v_{\zeta}^{(1)}(k_{i,\xi})
=
\Lambda^{-1}(k_{i,\xi}) q_{\zeta}^{(1)}(k_{i,\xi}),
\label{eqn:w_u_v_coefficients}
\end{flalign}
where $\Lambda(k_{i,\xi}) = r_{\zeta}^{(0)}(k_{i,\xi})+r_{\zeta}^{(1)}(k_{i,\xi})\tau$, and $\tau = \frac{2}{\Delta t}$.
Eq. (\ref{eqn:etod}) can be rewritten in matrix form as
\begin{flalign}
\label{eqn:g_update}
\mathbf{d}^{n+1}_{\parallel,\zeta,\xi}
=
\bar{\bar{\mathcal{W}}}_{\zeta,\xi}
\cdot
\bar{\bar{\mathcal{L}}}_{\zeta,\xi}
\cdot
\mathbf{e}^{n+1}_{\parallel,\zeta,\xi} + \mathbf{g}^{n+1}_{\parallel,\zeta,\xi}
\end{flalign}
with $\bar{\bar{\mathcal{W}}}_{\zeta,\xi} = \text{diag}\left\{
w_{\zeta}(k_{1,\xi}),w_{\zeta}(k_{2,\xi}),\cdots,w_{\zeta}(k_{N_{1},\xi}) \right\}$
for $\zeta=z,\rho$ and $\xi=\uparrow,\downarrow$.
Solving the above update equation separately with respect to $\zeta=z,\rho$ and $\xi=\uparrow,\downarrow$, we obtain
\begin{flalign}
\mathbf{d}_{\parallel}^{n+1}
=
\sum_{\xi=\left\{\uparrow,\downarrow\right\}}
\sum_{\zeta=\left\{z,\rho\right\}}
\mathbf{d}_{\parallel,\zeta,\xi} ^{n+1}.
\end{flalign}
Once $\mathbf{d}_{\parallel}^{n+1}$ and $\mathbf{g}_{\parallel}^{n}$ are obtained using the steps above,
$\mathbf{e}_{\parallel}^{n+1}$ can be finally obtained through
\begin{flalign}
\label{eqn:update_e}
\bar{\bar{\mathcal{A}}}_{\parallel}\cdot \mathbf{e}_{\parallel}^{n+1}
=
\mathbf{d}_{\parallel}^{n+1}
-
\mathbf{g}_{\parallel}^{n}
\end{flalign}
where
\begin{flalign}
\bar{\bar{\mathcal{A}}}_{\parallel}
=
\sum_{\xi=\left\{\uparrow,\downarrow\right\}}
\sum_{\zeta=\left\{z,\rho\right\}}
\bar{\bar{\mathcal{W}}}_{\zeta,\xi}\cdot\bar{\bar{\mathcal{L}}}_{\zeta,\xi}.
\end{flalign}

Before discussing the final time-update equations in the PML region, we need to discuss next the magnetic field relations in the PML.
The magnetic Hodge matrix in the PML is given by
\begin{flalign}
\left[\star_{\mu^{-1}}\right] = \sum_{k=1}^{N_{2}}
\left[{\bar{\bar{\mu}}'}^{(\text{PML})}(k)\right]_{\phi,\phi}
\bar{\bar{\mathcal{L}}}_{\phi}(k)
\end{flalign}
where
\begin{flalign}
\left[\bar{\bar{\mathcal{L}}}_{\phi}(k)\right]_{i,j}
=
\iint_{\mathcal{T}_{k}}
\left(\hat{\phi}\cdot\mathbf{W}^{(2)}_{i}(z,\rho)\right)
\left(\hat{\phi}\cdot\mathbf{W}^{(2)}_{j}(z,\rho)\right)
d A
=
\left\{
\begin{matrix}
\Omega_{k}, & \text{if~} (i,j)=(k,k) \\
0, & \text{otherwise}
\end{matrix}
\right.
\end{flalign}
and $\Omega_{k}$ denotes the area of $k$-th triangle of a 2D mesh.
Consequently, we obtain the following equation for the magnetic field:
\begin{flalign}
\left[\mathbf{h}_{\perp}\right]_{k}
=
\Omega_{k}^{-1}
\left(\left[{\bar{\bar{\mu}}'}^{(\text{PML})}(k)\right]_{\phi\phi}\right)^{-1}
\left[\mathbf{b}_{\perp}\right]_{k}
\end{flalign}
for $k=1,2,\cdots,N_{2}$.
Note that while the spatial support associated with each edge element consists of two triangular cells, the support of
each face element consists of one single cell (in a 2D mesh). Therefore,
there is no need for field splitting in this case.
Furthermore,
\begin{flalign}
\label{eqn:update_h}
\bar{\bar{\mathcal{A}}}_{\perp}
\cdot
\mathbf{h}_{\perp}^{n+\frac{1}{2}}
=
\mathbf{b}_{\perp}^{n+\frac{1}{2}}
-
\mathbf{g}_{\perp}^{n-\frac{1}{2}}
\end{flalign}
where $\bar{\bar{\mathcal{A}}}_{\perp}$ and $\mathbf{g}_{\perp}^{n-\frac{1}{2}}$ are analogous to $\bar{\bar{\mathcal{A}}}_{\parallel}$ and $\mathbf{g}_{\parallel}^{n}$, respectively, but now for the permeability.
For calculating $\mathbf{g}_{\perp}^{n-\frac{1}{2}}$, coefficients such as
in \label{eqn:d_e_time_domain_eqn} can be evaluated to give
\begin{flalign}
r_{\phi}^{(0)}(k_{\uparrow i})
& = \mu_0\rho,
\quad
r_{\phi}^{(1)}(k_{\uparrow i})
 = \frac{\mu_0}{\epsilon_0}\rho\sigma_\rho,
\\
q_{\phi}^{(0)}(k_{\uparrow i})
& = 0,
\quad
q_{\phi}^{(1)}(k_{\uparrow i})
= \Omega_{k}^{-1}.
\end{flalign}
The matrix $\bar{\bar{\mathcal{A}}}_{\perp}$ is given by the product of two diagonal matrices
\begin{flalign}
\bar{\bar{\mathcal{A}}}_{\perp}
=
\bar{\bar{\mathcal{W}}}_{\phi}
\cdot
\bar{\bar{\mathcal{L}}}_{\phi}.
\end{flalign}

The update equations for $\mathbf{b}_{\perp}$ and  $\mathbf{d}_{\parallel}$ directly follow from
 Eqs. (\ref{eqn:update_TE_1}) and (\ref{eqn:update_TE_2}) and write as
\begin{flalign}
\label{eqn:update_TEpml_1}
\mathbf{b}^{n+\frac{1}{2}}_{\perp}
&= \mathbf{b}^{n-\frac{1}{2}}_{\perp}
- \Delta t\bar{\mathbf{C}}\cdot\mathbf{e}^{n}_{\parallel},
\\
\label{eqn:update_TEpml_2}
\mathbf{d}^{n+1}_{\parallel}
&= \mathbf{d}^{n}_{\parallel}
+ \Delta t
\cdot
\left(\bar{\mathbf{C}}^{T}\cdot\mathbf{h}^{n+\frac{1}{2}}_{\perp}
- \mathbf{j}^{n+\frac{1}{2}}_{\parallel}\right).
\end{flalign}
The final update equations for the fields in the PML region consist on the cyclical application of  Eqs.~(\ref{eqn:update_TEpml_1}), (\ref{eqn:update_h}),  (\ref{eqn:update_TEpml_2}), and
(\ref{eqn:update_e}), in this order.

One can immediately carry out the implementation of the PML in the TM\textsuperscript{$\phi$} case
by reusing the TE\textsuperscript{$\phi$} results above and invoking the duality relations expressed by (\ref{eqn:duality}).

\bibliography{mybibfile}

\begin{thebibliography}{10}
\expandafter\ifx\csname url\endcsname\relax
  \def\url#1{\texttt{#1}}\fi
\expandafter\ifx\csname urlprefix\endcsname\relax\def\urlprefix{URL }\fi
\expandafter\ifx\csname href\endcsname\relax
  \def\href#1#2{#2} \def\path#1{#1}\fi

\bibitem{Na2017Axisymmetric}
D.-Y. Na, Y.~A. Omelchenko, H.~Moon, B.-H.~V. Borges, F.~L. Teixeira,
  Axisymmetric charge-conservative electromagnetic particle simulation
  algorithm on unstructured grids: Application to microwave vacuum electronic
  devices, J. Comp. Phys. 346 (2017) 295 -- 317.

\bibitem{Na2019Finite}
D.-Y. Na, B.-H.~V. Borges, F.~L. Teixeira, Finite element time-domain
  body-of-revolution {M}axwell solver based on discrete exterior calculus,
  Journal of Computational Physics 376 (2019) 249--275.

\bibitem{liu2011dynamic}
W.~Liu, M.~Fujimoto, The {D}ynamic {M}agnetosphere, Vol.~3, Springer, 2011.

\bibitem{ganguli2019understanding}
G.~Ganguli, C.~Crabtree, A.~C. Fletcher, L.~Rudakov, A.~S. Richardson, J.~Huba,
  C.~Siefring, W.~Amatucci, C.~D. Lewis, Understanding and harnessing the dual
  electrostatic/electromagnetic character of plasma turbulence in the
  near-earth space environment, Journal of Geophysical Research: Space Physics
  124~(12) (2019) 10365--10375.

\bibitem{carlsten2019radiation}
B.~E. Carlsten, P.~L. Colestock, G.~S. Cunningham, G.~L. Delzanno, E.~E. Dors,
  M.~A. Holloway, C.~A. Jeffery, J.~W. Lewellen, Q.~R. Marksteiner, D.~C.
  Nguyen, et~al., Radiation-belt remediation using space-based antennas and
  electron beams, IEEE Transactions on Plasma Science 47~(5) (2019) 2045--2063.

\bibitem{winske2012generation}
D.~Winske, W.~Daughton, Generation of lower hybrid and whistler waves by an ion
  velocity ring distribution, Physics of Plasmas 19~(7) (2012) 072109.

\bibitem{omelchenko2021rate}
Y.~A. Omelchenko, L.~Rudakov, J.~Ng, C.~Crabtree, G.~Ganguli, On the rate of
  energy deposition by an ion ring velocity beam, Physics of Plasmas 28~(5)
  (2021) 052102.

\bibitem{main2018excitation}
D.~Main, V.~Sotnikov, J.~Caplinger, D.~Rose, Excitation of electromagnetic
  whistler waves due to a parametric interaction between magnetosonic and lower
  oblique resonance modes in a cold, magnetized plasma, Physics of Plasmas
  25~(6) (2018) 062310.

\bibitem{delzanno2019high}
G.~L. Delzanno, V.~Roytershteyn, High-frequency plasma waves and pitch angle
  scattering induced by pulsed electron beams, Journal of Geophysical Research:
  Space Physics 124~(9) (2019) 7543--7552.

\bibitem{VINCENTI201822}
H.~Vincenti, J.-L. Vay, Ultrahigh-order {M}axwell solver with extreme
  scalability for electromagnetic {PIC} simulations of plasmas, Computer
  Physics Communications 228 (2018) 22--29.

\bibitem{vincenti2017pic}
H.~Vincenti, M.~Lobet, R.~Lehe, J.-L. Vay, J.~Deslippe, {PIC} codes on the road
  to exascale architectures, in: Exascale Scientific Applications, Chapman and
  Hall/CRC, 2017, pp. 375--408.

\bibitem{10.1063/5.0046842}
S.~O'Connor, Z.~D. Crawford, O.~H. Ramachandran, J.~Luginsland, B.~Shanker,
  {Time integrator agnostic charge conserving finite element {PIC}}, Physics of
  Plasmas 28~(9) (09 2021).

\bibitem{Na2018Relativistic}
D.-Y. Na, H.~Moon, Y.~A. Omelchenko, F.~L. Teixeira, Relativistic extension of
  a charge-conservative finite element solver for time-dependent
  {M}axwell-{V}lasov equations, Physics of Plasmas 25~(1) (2018) 013109.

\bibitem{10.1063/1.1384387}
D.~L. Bruhwiler, R.~Giacone, J.~R. Cary, J.~P. Verboncoeur, P.~Mardahl,
  E.~Esarey, W.~Leemans, {Modeling beam-driven and laser-driven plasma
  wakefield accelerators with {XOOPIC}}, AIP Conference Proceedings 569~(1)
  (2001) 591--604.

\bibitem{mahalingam2010particle}
S.~Mahalingam, J.~A. Menart, Particle-based plasma simulations for an ion
  engine discharge chamber, Journal of Propulsion and Power 26~(4) (2010)
  673--688.

\bibitem{LEHE201666}
R.~Lehe, M.~Kirchen, I.~A. Andriyash, B.~B. Godfrey, J.-L. Vay, A spectral,
  quasi-cylindrical and dispersion-free particle-in-cell algorithm, Computer
  Physics Communications 203 (2016) 66--82.

\bibitem{MASSIMO2016841}
F.~Massimo, S.~Atzeni, A.~Marocchino, Comparisons of time explicit hybrid
  kinetic-fluid code architect for plasma wakefield acceleration with a full
  {PIC} code, Journal of Computational Physics 327 (2016) 841--850.

\bibitem{Teixeira1999Lattice}
F.~L. Teixeira, W.~C. Chew, Lattice electromagnetic theory from a topological
  viewpoint, J. Math. Phys. 40 (1999) 169--187.

\bibitem{Teixeira1999Differential}
F.~Teixeira, W.~Chew, Differential forms, metrics, and the reflectionless
  absorption of electromagnetic waves, J. Electromagn. Waves Appl. 13 (1999)
  665--686.

\bibitem{Pendry2006Controlling}
J.~B. Pendry, D.~Schurig, D.~R. Smith, Controlling electromagnetic fields,
  Science 312 (2006) 1780--1782.

\bibitem{HE20061}
B.~He, F.~Teixeira, Geometric finite element discretization of maxwell
  equations in primal and dual spaces, Physics Letters A 349~(1) (2006) 1--14.

\bibitem{doi:10.1137/08073901X}
M.~E. Rognes, R.~C. Kirby, A.~Logg, Efficient assembly of {H}(div) and
  {H}(curl) conforming finite elements, SIAM Journal on Scientific Computing
  31~(6) (2010) 4130--4151.

\bibitem{monkFEM}
P.~Monk, {F}inite {E}lement {M}ethods for {M}axwell's {E}quations, Oxford
  University press, 2003.

\bibitem{campos2016constructing}
M.~Campos~Pinto, Constructing exact sequences on non-conforming discrete
  spaces, Comptes Rendus. Math{\'e}matique 354~(7) (2016) 691--696.

\bibitem{pinto2022semi}
M.~C. Pinto, V.~Pag{\`e}s, A semi-implicit electromagnetic {FEM-PIC} scheme
  with exact energy and charge conservation, Journal of Computational Physics
  453 (2022) 110912.

\bibitem{Ramachandran2023Review}
O.~H. Ramachandran, L.~C. Kempel, J.~P. Verboncoeur, B.~Shanker, A necessarily
  incomplete review of electromagnetic finite element particle-in-cell methods,
  IEEE Transactions on Plasma Science (2023) 1--11\href
  {https://doi.org/10.1109/TPS.2023.3257165}
  {\path{doi:10.1109/TPS.2023.3257165}}.

\bibitem{Bossavit1988Whitney}
A.~Bossavit, Whitney forms: A class of finite elements for three-dimensional
  computations in electromagnetism, IEE Proc., Part A: Phys. Sci., Meas.
  Instrum., Manage. Educ. 135 (1988) 493--500.

\bibitem{GILLETTE20111213}
A.~Gillette, C.~Bajaj, Dual formulations of mixed finite element methods with
  applications, Computer-Aided Design 43~(10) (2011) 1213--1221, solid and
  Physical Modeling 2010.

\bibitem{LOHI2021113520}
J.~Lohi, L.~Kettunen, Whitney forms and their extensions, Journal of
  Computational and Applied Mathematics 393 (2021) 113520.

\bibitem{10.1007/978-3-319-01601-6}
A.~Bossavit, F.~Rapetti, Whitney forms, from manifolds to fields, in:
  M.~Aza{\"i}ez, H.~El~Fekih, J.~S. Hesthaven (Eds.), Spectral and High Order
  Methods for Partial Differential Equations - ICOSAHOM 2012, Springer
  International Publishing, Cham, 2014, pp. 179--189.

\bibitem{pinto2016electromagnetic}
M.~C. Pinto, M.~Lutz, M.~Mounier, Electromagnetic {PIC} simulations with smooth
  particles: a numerical study, ESAIM: Proceedings and Surveys 53 (2016)
  133--148.

\bibitem{gross2004electromagnetic}
P.~W. Gross, P.~R. Kotiuga, Electromagnetic theory and computation: a
  topological approach, Vol.~48, Cambridge University Press, 2004.

\bibitem{Donderici2008Mixed}
B.~Donderici, F.~L. Teixeira, Mixed finite-element time-domain method for
  transient {M}axwell equations in doubly dispersive media, IEEE Transactions
  on Microwave Theory and Techniques 56~(1) (2008) 113--120.

\bibitem{Kim2011Parallel}
J.~Kim, F.~L. Teixeira, Parallel and explicit finite-element time-domain method
  for {M}axwell's equations, IEEE Trans. Antennas Propag. 59 (2011) 2350--2356.

\bibitem{He2007Differential}
B.~He, F.~L. Teixeira, Differential forms, {G}alerkin duality, and sparse
  inverse approximations in finite element solutions of {M}axwell equations,
  IEEE Trans. Antennas Propag. 55 (2007) 1359--1368.

\bibitem{Teixeira2013Differential}
F.~L. Teixeira, Differential forms in lattice field theories: An overview, ISRN
  Math. Phys. (2013) 487270.

\bibitem{Na2016Local}
D.-Y. Na, H.~Moon, Y.~A. Omelchenko, F.~L. Teixeira, Local, explicit, and
  charge-conserving electromagnetic particle-in-cell algorithm on unstructured
  grids, IEEE Trans. Plasma Sci. 44 (2016) 1353--1362.

\bibitem{bochev2006principles}
P.~B. Bochev, J.~M. Hyman, Principles of mimetic discretizations of
  differential operators, in: Compatible spatial discretizations, Springer,
  2006, pp. 89--119.

\bibitem{arnold2010finite}
D.~Arnold, R.~Falk, R.~Winther, Finite element exterior calculus: from {H}odge
  theory to numerical stability, Bulletin of the American mathematical society
  47~(2) (2010) 281--354.

\bibitem{PIER2001}
F.~L. Teixeira~(ed.), {G}eometric {M}ethods in {C}omputational
  {E}lectromagnetics, Progress In Electromagnetics Research 32, EMW Publishing,
  Cambridge, Mass., 2001.

\bibitem{hirani2003discrete}
A.~N. Hirani, Discrete exterior calculus, Ph.D. thesis, California Institute of
  Technology, 2003.

\bibitem{Teixeira2014Lattice}
F.~L. Teixeira, Lattice {M}axwell's equations, Progress In Electromagnetics
  Research 148 (2014) 113--128.

\bibitem{arnold2018finite}
D.~N. Arnold, Finite element exterior calculus, SIAM, 2018.

\bibitem{9599150}
Z.~D. Crawford, S.~O’Connor, J.~Luginsland, B.~Shanker, Rubrics for charge
  conserving current mapping in finite element electromagnetic particle in cell
  methods, IEEE Transactions on Plasma Science 49~(11) (2021) 3719--3732.
\newblock \href {https://doi.org/10.1109/TPS.2021.3122410}
  {\path{doi:10.1109/TPS.2021.3122410}}.

\bibitem{squire2012geometric}
J.~Squire, H.~Qin, W.~M. Tang, Geometric integration of the {V}lasov-{M}axwell
  system with a variational particle-in-cell scheme, Physics of Plasmas 19~(8)
  (2012) 084501.

\bibitem{PhysRevE.61.3174}
S.~Sen, S.~Sen, J.~C. Sexton, D.~H. Adams, Geometric discretization scheme
  applied to the {A}belian {C}hern-{S}imons theory, Phys. Rev. E 61 (2000)
  3174--3185.

\bibitem{Moon2015Exact}
H.~Moon, F.~L. Teixeira, Y.~A. Omelchenko, Exact charge-conserving
  scatter–gather algorithm for particle-in-cell simulations on unstructured
  grids: {A} geometric perspective, Comput. Phys. Commun. 194 (2015) 43--53.

\bibitem{osti_950065}
M.~L. Stowell, \href{https://www.osti.gov/biblio/950065}{Discretizing transient
  current densities in the maxwell equations} (11 2008).
\newline\urlprefix\url{https://www.osti.gov/biblio/950065}

\bibitem{campospinto:hal-01303852}
M.~Campos~Pinto, E.~Sonnendr{\"u}cker,
  \href{https://hal.science/hal-01303852}{{Compatible Maxwell solvers with
  particles I: conforming and non-conforming 2D schemes with a strong Ampere
  law}}, {SMAI Journal of Computational Mathematics} 3 (2017) 53--89.
\newblock \href {https://doi.org/10.5802/smai-jcm.20}
  {\path{doi:10.5802/smai-jcm.20}}.
\newline\urlprefix\url{https://hal.science/hal-01303852}

\bibitem{Na2019Polynomial}
D.-Y. Na, F.~L. Teixeira, W.~C. Chew, Polynomial finite-size shape functions
  for electromagnetic particle-in-cell algorithms based on unstructured meshes,
  IEEE Journal on Multiscale and Multiphysics Computational Techniques 4 (2019)
  317--328.

\bibitem{pinto2014charge}
M.~C. Pinto, S.~Jund, S.~Salmon, E.~Sonnendr{\"u}cker, Charge-conserving
  {FEM--PIC} schemes on general grids, Comptes Rendus Mecanique 342~(10-11)
  (2014) 570--582.

\bibitem{10.1145/1141911.1141991}
K.~Wang, Weiwei, Y.~Tong, M.~Desbrun, P.~Schr\"{o}der, Edge subdivision schemes
  and the construction of smooth vector fields, ACM Trans. Graph. 25~(3) (2006)
  1041–1048.

\bibitem{campos2015towards}
M.~Campos~Pinto, Towards smooth particle methods without smoothing, Journal of
  Scientific Computing 65~(1) (2015) 54--82.

\bibitem{Github_BORPIC_PP}
F.~L.~T. Dong-Yeop~Na, Y.~A. Omelchenko, Borpic-pp-matlab-code,
  \url{https://github.com/DYNPOSTECH/BORPIC-PP-MATLAB-CODE} (2023).

\bibitem{teixeira1997systematic}
F.~Teixeira, W.~C. Chew, Systematic derivation of anisotropic {PML} absorbing
  media in cylindrical and spherical coordinates, IEEE microwave and guided
  wave letters 7~(11) (1997) 371--373.

\end{thebibliography}

\end{document}